\documentclass[a4paper,10pt,twoside]{article}
\usepackage{german,amssymb,amsmath,latexsym,epsfig,graphics,graphicx,mathrsfs,color,amsthm}
\usepackage[latin1]{inputenc}
\usepackage{fancyhdr} 
\newtheorem{theorem}{Theorem}[section]

\newtheorem{definition}[theorem]{Definition}
\newtheorem{lemma}[theorem]{Lemma}

\newtheorem{folgerung}[theorem]{Folgerung}
\newtheorem{beispiel}[theorem]{Beispiel}
\newtheorem{satz}[theorem]{Satz}
\newtheorem{bemerkung}[theorem]{Bemerkung}

\numberwithin{equation}{section}
\newcommand{\R}{{\mathbb R}}
\newcommand{\N}{{\mathbb N}}

\begin{document}
\thispagestyle{empty}
\rule{1.0\textwidth}{1.0pt}

\vspace*{5mm}
{\bf{\large Das Pontrjaginsche Maximumprinzip f"ur nichtlineare \\[2mm] Steuerungsprobleme mit unendlichem Zeithorizont \\[2mm]
            unter Zustandsbeschr"ankungen} \\[10mm]
Nico Tauchnitz} \\[25mm]
{\bf Vorwort} \\[2mm]
In der vorliegenden Ausarbeitung stelle ich notwendige und hinreichende Optimalit"atsbedingungen f"ur
starke lokale Extrema in Aufgaben mit unendlichem Zeithorizont vor.
Im Vergleich zu den
klassischen Steuerungsproblemen besitzen Aufgaben mit unendlichem Zeithorizont ihren eigenen Charakter,
da durch das unbeschr"ankte Zeitintervall die Aufgabenstellung eine Singularit"at beinhaltet.
Die L"osungsmethoden, die f"ur die klassischen Aufgaben entwickelt wurden,
k"onnen f"ur die Situation des unendlichen Horizontes nicht einfach "ubernommen werden. \\[2mm]
Mit der vorliegenden Arbeit pr"asentiere ich Ergebnisse,
die die Verfahren von Dubovitskii \& Milyutin und Ioffe \& Tichomirov an die Probleme der Optimalen Steuerung mit
unendlichem Zeithorizont anpassen.
Als Herausforderungen erweisen sich neben dem unbeschr"ankten Zeitintervall
die Einf"uhrung mehrfacher Nadelvariationen und die Behandlung einer lokal unbeschr"ankten Verteilungsfunktion,
wie z.\;B. eine Weibull-Verteilung, im Integranden.
Zur Untersuchung von Zustandsbeschr"ankungen f"uhre ich den Raum der stetigen Funktionen, die im Unendlichen einen Grenzwert besitzen, ein.
Die erzielten Resultate weisen eine direkte Verwandtschaft zu den Ergebnissen zu den  Steuerungsproblemen mit endlichem Zeithorizont auf.
Dadurch wird eine Einordnung der Aufgabenklassen zueinander m"oglich. \\[2mm]
Die vorliegende "uberarbeitete Version behandelt Steuerungsprobleme unter restriktiven Annahmen an die Daten der Aufgabe.
Einen alternativen Zugang zu schwachen lokalen Minimalstellen in Steuerungsproblemen mit unendlichem Zeithorizont im Rahmen gewichteter Funktionenr"aume
sind in der Arbeit \cite{TauchnitzWMPIHOC} angegeben. 
Eine umfassendere Darstellung kann man \cite{TauchnitzOC} entnehmen.
Bei n"aherer Betrachtung der Ergebnisse in \cite{TauchnitzWMPIHOC} zeigen sich aber die Schwierigkeiten und Pathologien,
die sich im Rahmen gewichteter Funktionenr"aume ergeben k"onnen.
Diesbez"uglich habe ich ausf"uhrliche Bemerkungen in \cite{TauchnitzWMPIHOC} aufgef"uhrt. \\[5mm]
Juli 2018

\newpage
\lhead[\thepage \, Inhaltsverzeichnis]{Optimale Steuerung mit unendlichem Zeithorizont}
\rhead[Optimale Steuerung mit unendlichem Zeithorizont]{Inhaltsverzeichnis \thepage}
\tableofcontents

\newpage
\lhead[\thepage \hspace*{1mm} Starkes lokales Minimum]{}
\rhead[]{Starkes lokales Minimum \hspace*{1mm} \thepage}
\section{Starkes lokales Minimum "uber unendlichem Zeithorizont} \label{KapitelStrong}
Die Optimale Steuerung mit unendlichem Zeithorizont liefert die wesentliche Grundlage zur Formulierung und Untersuchung
von Aufgaben in der "Okonomischen Wachstumstheorie.
Im Rahmen der "Okonomischen Wachstumstheorie werden z.\,B. die Interaktionen sich "uberschneidender Generationen oder die Determinanten 
des wirtschaftlichen Wachstums,
insbesondere unter sich "andernden Umweltbedingungen wie globaler Erw"armung oder ersch"opfenden nat"urlichen Ressourcen, untersucht.
Aufgrund der Langlebigkeit der wirtschafts- und sozialpolitischen Entscheidungen muss dabei die Frage nach einem geeigneten 
Planungszeitraum aufgeworfen werden:
Jeder endliche Zeithorizont stellt die Forderung nach einer ad"aquaten Ausgangslage f"ur die nachfolgenden Generationen.
Um die Beachtung aller nachfolgenden Generationen zu gew"ahrleisten,
wird der Zeitrahmen in Form des unendlichen Zeithorizontes idealisiert (Arrow \& Kurz \cite{Arrow}). \\[2mm]
Der erste mathematische Beitrag zu einem Problem mit unbeschr"anktem Zeitintervall besteht in einer Aufgabe der Variationsrechnung,
in der die Frage nach der optimalen Sparquote einer Gesellschaft behandelt wird (Ramsey \cite{Ramsey}).
Bei n"aherer Betrachtung entsteht in diesem Problem nicht nur die Aufgabe die Optimalit"atsbedingungen "uber dem unendlichen Zeithorizont zu formulieren,
sondern insbesondere die Gestalt der Transversalit"atsbedingungen im Unendlichen zu charakterisieren. 
Diese Fragestellung stellt allerdings eine schwerwiegende Herausforderung dar,
denn die bekannten Resultate k"onnen nicht einfach "ubernommen und die bekannten Methoden k"onnen nicht
unmittelbar an das unbeschr"ankte Intervall angepasst werden. \\[2mm]
Von den wenigen uns bekannten Resultaten,
die den unendlichen Zeithorizont nicht auf ein endliches Intervall reduzieren,
z"ahlen wir die Arbeiten von Brodskii \cite{Brodskii} und Pickenhain \cite{Pickenhain} auf.
In diesen Beitr"agen werden auf der Basis von schwachen lokalen Variationen und der Anwendung geeigneter funktionalanalytischer Methoden
notwendige Optimalit"atsbedingungen f"ur Steuerungsprobleme mit unendlichem Zeithorizont erzielt.
In \cite{Brodskii} wird eine sehr allgemeine Aufgabenklasse mit Zustandsbeschr"ankungen und Randbedingungen im Unendlichen betrachtet.
Aufgrund der Wahl des Raumes der messbaren und beschr"ankten Funktionen f"uhrt der funktionalanalytische Rahmen allerdings
zu keiner ``"asthetischen'' Darstellung der Adjungierten.
Demgegen"uber bezieht sich die Arbeit \cite{Pickenhain} auf Aufgaben mit eindimensionalen linearen Nebenbedingungen und mit freiem rechten Endpunkt.
Der innovative Beitrag in \cite{Pickenhain} ist die Wahl des gewichteten Sobolev-Raumes,
in dessen Rahmen die ``nat"urlichen'' Transversalit"atsbedingungen in gewisser Weise Eigenschaften der Elemente dieser R"aume sind.
An der Arbeit \cite{Pickenhain} ist jedoch anzumerken,
dass sich die angewandte Beweismethode auf linear-quadratische Aufgaben fokussiert.
Die Anwendbarkeit der Beweisstrategie im Fall von Aufgaben mit einer nichtlinearen Dynamik ist offen. \\
In dieser Arbeit pr"asentieren wir das Pontrjaginsche Maximumprinzip f"ur Steuerungsprobleme mit unendlichem Zeithorizont.
Dabei erweitern wir die mehrfachen Nadelvariationen nach Ioffe \& Tichomirov \cite{Ioffe} f"ur das unbeschr"ankte Zeitintervall.
Weiterhin gehen wir auf Randbedingungen im Unendlichen und auf Zustandsbeschr"ankungen "uber dem unendlichen Zeithorizont ein.
Zur Aufgabenklasse mit Zustandsbeschr"ankungen existieren bisher nur sehr wenige Beitr"age.
Wir erarbeiten einen Zugang im Rahmen der stetigen Funktionen, die im Unendlichen einen Grenzwert besitzen.
Die Anforderung, die wir dabei an die Aufgabe stellen,
sind verwandt mit denen bei Brodskii \cite{Brodskii}.
Die Darstellung der notwendigen Optimalit"atsbedingungen besitzen wegen diesem Zugang eine enge Verwandtschaft mit dem Maximumprinzip f"ur Aufgaben mit endlichem
Zeitintervall. \\[2mm]
Ein sehr wichtiges und zentrales Anliegen dieser Arbeit ist die Einordnung der Aufgabe mit unendlichem Zeithorizont in die Theorie der Optimalen Steuerung.
Wir gehen dazu auf die h"aufig angewandte Approximation mit endlichem Zeithorizont und auf die R"uckf"uhrung auf ein endliches Zeitintervall mittels der
Substitution der Zeit ein.
Es zeigt sich,
dass einerseits bei der Approximation der vollst"andige Satz an notwendigen Optimalit"atsbedingungen verloren geht
und andererseits die Substitution der Zeit den unendlichen Zeithorizont in eine Singularit"at in der Aufgabe mit endlichem Zeitintervall "uberf"uhrt.
Umgekehrt k"onnen wir wiederum zeigen,
dass das Pontrjaginsche Maximumprinzip im Kapitel zu starken lokalen Minimalstellen "uber einem endlichen Zeitintervall eine direkte Konsequenz der Ergebnisse
dieses Kapitels darstellt.
Somit sind die Steuerungsprobleme "uber endlichem Zeithorizont eine echte Unterklasse der Steuerungsprobleme mit unendlichem Zeithorizont.
Die Untersuchungen liefern damit nicht nur technisch bez"uglich der Nadelvariationsmethode und der Behandlung der Zustandsbeschr"ankungen innovative Beitr"age,
sondern es zeigt sich viel mehr die besondere Qualit"at der Aufgabe mit unendlichem Zeithorizont als eine echte Verallgemeinerung der Aufgabenklasse mit
endlichem Zeithorizont.

       \rhead[]{Aufgabenstellung \hspace*{1mm} \thepage}
       \section{Die Aufgabenstellung}
In diesem Kapitel betrachten wir starke lokale Minimalstellen der Aufgabe
\begin{eqnarray}
&& \label{PAUH1} J\big(x(\cdot),u(\cdot)\big) = \int_0^\infty \omega(t)f\big(t,x(t),u(t)\big) \, dt \to \inf, \\
&& \label{PAUH2} \dot{x}(t) = \varphi\big(t,x(t),u(t)\big), \\
&& \label{PAUH3} h_0\big(x(0)\big)=0, \qquad \lim_{t \to \infty} h_1\big(t,x(t)\big)=0, \\
&& \label{PAUH4} u(t) \in U \subseteq \R^m, \quad U \not= \emptyset, \\
&& \label{PAUH5} g_j\big(t,x(t)\big) \leq 0 \quad \mbox{f"ur alle } t \in \R_+, \quad j=1,...,l.
\end{eqnarray}
Dabei ist $\omega(\cdot) \in L_1(\R_+,\R_+)$ und es gelten f"ur die eingehenden Abbildungen
$$f:\R \times \R^n \times \R^m \to \R, \qquad \varphi:\R \times \R^n \times \R^m \to \R^n, \qquad g_j: \R \times \R^n \to \R,$$
sowie f"ur die Randbedingungen
$$h_0:\R^n \to \R^{s_0}, \qquad h_1:\R \times \R^n \to \R^{s_1}.$$
Wir nennen die Trajektorie $x(\cdot)$ eine L"osung des dynamischen Systems (\ref{PAUH2}),
falls $x(\cdot)$ auf $\R_+$ definiert ist und auf jedem endlichen Intervall die Dynamik mit Steuerung $u(\cdot)$
im Sinn von Carath\'eodory l"ost. \\[2mm]
Zu $x(\cdot)$ bezeichne $V_\gamma$ die Menge
$V_\gamma= \{ (t,x) \in \overline{\R}_+ \times \R^n \,|\, \|x-x(t)\| \leq \gamma\}$.
Dann geh"oren zur Menge $\mathscr{B}_{\rm Lip}$ diejenigen
$x(\cdot) \in W^1_\infty(\R_+,\R^n)$,
f"ur die es zu jeder kompakten Menge $U_1 \subseteq \R^m$ eine Zahl $\gamma>0$ derart gibt,
dass auf der Menge $V_\gamma \times U_1$ die Abbildungen
\begin{enumerate}
\item[(B)] $f(t,x,u)$, $\varphi(t,x,u)$, $g_j(t,x)$ und $h_0(x)$, $h_1(t,x)$ gleichm"a"sig stetig und
           gleichm"a"sig stetig differenzierbar bez"uglich $x$ sind.
\end{enumerate}

Der Steuerungsprozess $\big(x(\cdot),u(\cdot)\big) \in W^1_\infty(\R_+,\R^n;\nu) \times L_\infty(\R_+,U)$
hei"st zul"assig in der Aufgabe (\ref{PAUH1})--(\ref{PAUH5}),
falls $\big(x(\cdot),u(\cdot)\big)$ dem System (\ref{PAUH2}) gen"ugt,
die Randbedingungen (\ref{PAUH3}) und Restriktionen (\ref{PAUH4}) erf"ullt,
sowie das Lebesgue-Integral in (\ref{PAUH1}) endlich ist.
Die Menge $\mathscr{B}_{\rm adm}$ bezeichnet die Menge der zul"assigen Steuerungsprozesse der Aufgabe (\ref{PAUH1})--(\ref{PAUH5}). \\
Der zul"assige Steuerungprozess $\big(x_*(\cdot),u_*(\cdot)\big)$ ist eine starke lokale Minimalstelle\index{Minimum, schwaches lokales!--, starkes lokales} in
der Aufgabe (\ref{PAUH1})--(\ref{PAUH5}),
falls eine Zahl $\varepsilon > 0$ derart existiert, dass die Ungleichung
$$J\big(x(\cdot),u(\cdot)\big) \geq J\big(x_*(\cdot),u_*(\cdot)\big)$$
f"ur alle $\big(x(\cdot),u(\cdot)\big) \in \mathscr{B}_{\rm adm}$ mit 
$\|x(\cdot)-x_*(\cdot)\|_\infty \leq \varepsilon$ gilt.

       \rhead[]{Pontrjaginsches Maximumprinzip \hspace*{1mm} \thepage}
       \section{Das Pontrjaginsche Maximumprinzip} \label{AbschnittPMPUH}
\subsection{Notwendige Optimalit"atsbedingungen}
Im Weiteren bezeichnet $H: \R \times \R^n \times \R^m \times \R^n \times \R \to \R$ die Pontrjagin-Funktion
$$H(t,x,u,p,\lambda_0) = \langle p, \varphi(t,x,u) \rangle-\lambda_0 \omega(t)f(t,x,u).$$

\begin{theorem}[Pontrjaginsches Maximumprinzip] \label{SatzPAUH}\index{Pontrjaginsches Maximumprinzip}
Es sei $\big(x_*(\cdot),u_*(\cdot)\big) \in \mathscr{B}_{\rm adm} \cap \mathscr{B}_{\rm Lip}$.
Weiterhin nehmen wir an,
dass
\begin{equation} \label{PMPBedingung}
\int_0^\infty \big\|\varphi\big(t,x_*(t),u_*(t)\big)\big\| \, dt < \infty, \qquad \int_0^\infty \big\|\varphi_x\big(t,x_*(t),u_*(t)\big)\big\| \, dt < \infty
\end{equation}
ausfallen und es m"oge zu jedem $\delta>0$ ein $T>0$ existieren mit
\begin{eqnarray}
&& \int_T^\infty \big\| \varphi\big(t,x(t),u_*(t)\big)-\varphi\big(t,x'(t),u_*(t)\big) - \varphi_x\big(t,x_*(t),u_*(t)\big)\big(x(t)-x'(t)\big) \big\| \, dt
   \nonumber \\
&& \label{PMPBedingung2} \hspace*{20mm} \leq \delta \|x(\cdot)-x'(\cdot)\|_\infty
\end{eqnarray}
f"ur alle $x(\cdot), x'(\cdot) \in W^1_\infty(\R_+,\R^n)$ mit $\|x(\cdot)-x_*(\cdot)\|_\infty < \gamma$, $\|x'(\cdot)-x_*(\cdot)\|_\infty < \gamma$. \\[2mm]
Ist $\big(x_*(\cdot),u_*(\cdot)\big)$ ein starkes lokales Minimum der Aufgabe (\ref{PAUH1})--(\ref{PAUH4}),
dann existieren nicht gleichzeitig verschwindende Multiplikatoren $\lambda_0 \geq 0$,
$p(\cdot) \in W^1_\infty(\R_+,\R^n)$ und $l_i \in \R^{s_i}$, $i=0,1$, derart, dass
\begin{enumerate}
\item[(a)] die Funktion $p(\cdot)$ fast "uberall der adjungierten Gleichung\index{adjungierte Gleichung}
           \begin{equation}\label{SatzPAUH1}
           \dot{p}(t) = -\varphi_x^T\big(t,x_*(t),u_*(t)\big) p(t) + \lambda_0 \omega(t)f_x\big(t,x_*(t),u_*(t)\big)
           \end{equation}
           gen"ugt und die Transversalit"atsbedingungen\index{Transversalit"atsbedingungen}
           \begin{equation}\label{SatzPAUH2}
           p(0) = {h_0'}^T\big(x_*(0)\big)l_0, \qquad \lim_{t \to \infty} p(t)= - \lim_{t \to \infty}h_{1x}^T\big(t,x_*(t)\big) l_1
           \end{equation}
           erf"ullt;
\item[(b)] in fast allen Punkten $t \in \R_+$ die Maximumbedingung gilt:
           \begin{equation}\label{SatzPAUH3}
           H\big(t,x_*(t),u_*(t),p(t),\lambda_0\big) = \max_{u \in U} H\big(t,x_*(t),u,p(t),\lambda_0\big).
           \end{equation}
\end{enumerate}
\end{theorem}

Im Pontrjaginsche Maximumprinzip ergibt sich aus der ersten Bedingung in (\ref{PMPBedingung}),
dass die Trajektorie $x_*(\cdot)$ im Unendlichen einen Grenzwert besitzt und damit dem Raum $C_{\lim}(\R_+,\R^n)$ angeh"ort.
Die zweite Bedingung in (\ref{PMPBedingung}) stellt sicher,
dass die Dynamik als Abbildung im Rahmen des Raumes $C_{\lim}(\R_+,\R^n)$ den Anforderungen des Extremalprinzips \ref{SatzExtremalprinzipStark} gen"ugt.
F"ur die Einf"uhrung mehrfacher Nadelvariationen "uber dem unendlichen Zeithorizont liefert wiederum die Bedingung (\ref{PMPBedingung2}) die Grundlage. \\[2mm]
Die strengen Einschr"ankungen (\ref{PMPBedingung}) und (\ref{PMPBedingung2}) haben zur Folge,
dass das Pontrjaginsche Maximumprinzip nicht auf lineare Systeme mit konstanten Koeffizienten angewendet werden kann.
An dieser Stelle liefert der gewichtete Rahmen mit dem Schwachen Optimalit"tsprinzip in der Arbeit \cite{TauchnitzWMPIHOC}
einen wesentlich besser angepassten Zugang.

\begin{beispiel} {\rm Wir betrachten die Aufgabe
\begin{eqnarray*}
&& J\big(x(\cdot),u(\cdot)\big) = \int_0^\infty e^{-\varrho t} \big(1-u(t)\big) x(t) \, dt \to \sup,\\
&& \dot{x}(t)=u(t)x(t), \quad x(0)=1, \quad \lim_{t \to \infty} x(t)=x_1>1, \quad u(t) \in [0,1], \quad \varrho \in (0,1).
\end{eqnarray*}
Wir stellen zun"achst die Pontrjagin-Funktion auf:
$$H(t,x,z,u,p,q,1) = pux+\lambda_0 e^{-\varrho t}(1-u)x.$$
Mit Hilfe der Bedingungen (\ref{SatzPAUH1})--(\ref{SatzPAUH3}) k"onnen wir den Fall $\lambda_0=0$ ausschlie"sen.
Weiterhin ergeben sich der Steuerungsprozess
$$x_*(t) = \left\{ \begin{array}{ll} e^t, & t \in [0,\tau), \\ x_1, & t \in [\tau, \infty), \end{array}\right. \quad
  u_*(t)= \left\{ \begin{array}{ll} 1, & t \in [0,\tau), \\ 0, & t \in [\tau, \infty), \end{array}\right. \quad \tau=\ln x_1$$
und die Adjungierte
$$p(t) = \left\{ \begin{array}{ll}
          e^{(1-\varrho)\tau} e^{-t}, & t \in [0,\tau), \\
          \frac{\varrho-1}{\varrho}e^{-\varrho \tau} + \frac{1}{\varrho} e^{-\varrho t}, & t \in [\tau, \infty). \end{array}\right.$$
F"ur die Adjungierte gilt dabei im Unendlichen $\displaystyle\lim_{t \to \infty} p(t)= \frac{\varrho-1}{\varrho}e^{-\varrho \tau} \not=0$.
Damit konnten wir aus den notwendigen Bedingungen (\ref{SatzPAUH1})--(\ref{SatzPAUH3}) des Maximumprinzips einen eindeutigen
Kandidaten bestimmen. \hfill $\square$}
\end{beispiel}

\begin{beispiel} \label{BeipielRWUnendlich}
{\rm Wir betrachten die Aufgabe
\begin{eqnarray*}
&& J\big(x(\cdot),u(\cdot)\big)=\int_0^\infty e^{-\varrho t}\big(1-u(t)\big)x(t) \, dt \to \sup, \\
&& \dot{x}(t)=u(t)x(t), \qquad x(0)=1, \qquad u \in [0,1], \qquad \varrho \in (0,1)
\end{eqnarray*}
mit der Budgetbeschr"ankung
$$\int_0^\infty e^{-\varrho t}x(t) \, dt = Z, \qquad Z> \frac{1}{\varrho}.$$
Bez"uglich der Budgetbeschr"ankung f"uhren wir die folgende Zustandsgleichung mit Randwert im Unendlichen ein:
$$\dot{z}(t)=e^{-\varrho t}x(t), \qquad z(0)=0, \qquad \lim_{t \to \infty} z(t) = Z.$$
Die Zustandsgleichung und -beschr"ankung f"ur die Trajektorie $z(\cdot)$ ergibt sich aus der isoperimetrischen Nebenbedingung in Form einer Budgetbeschr"ankung
$$\int_0^\infty e^{-\varrho t}x(t) \, dt \leq Z.$$
Offensichtlich ist $\dot{z}(t) >0$ auf $\R_+$ und damit $z(t)$ streng monoton wachsend.
Demzufolge kann die Beschr"ankung $z(t) \leq Z$ erst im Unendlichen aktiv werden und greift nur durch das Verhalten in $t=\infty$ in die
gestellte Aufgabe ein.
Da stets $\dot{x}(t)\geq 0$ ist,
muss f"ur jede zul"assige Trajektorie $e^{-\varrho t}x(t) \to 0$ f"ur $t \to \infty$ gelten,
denn nur dann ist $z(t) \leq Z$ erf"ullt.
Damit erhalten wir f"ur zul"assige Steuerungsprozesse zun"achst im Zielfunktional
\begin{eqnarray*}
    J\big(x(\cdot),z(\cdot),u(\cdot)\big)
&=& \int_0^\infty e^{-\varrho t}\big(1-u(t)\big)x(t) \, dt = \int_0^\infty \dot{z}(t) \, dt - \int_0^\infty e^{-\varrho t} \dot{x}(t) \, dt \\
&=& \int_0^\infty \dot{z}(t) \, dt + 1 - \varrho \int_0^\infty e^{-\varrho t} x(t) \, dt \leq 1 + (1-\varrho) Z.
\end{eqnarray*}
Es ergibt sich also f"ur das Zielfunktional die obere Schranke $1 + (1-\varrho) Z$.
D.\,h., dass jeder Steuerungsprozess $\big(x(\cdot),z(\cdot),u(\cdot)\big)$,
f"ur den die Zustandsbeschr"ankung $z(t) \leq Z$ im Unendlichen aktiv wird,
global optimal ist und $J\big(x(\cdot),z(\cdot),u(\cdot)\big) = 1 + (1-\varrho) Z$ gilt. \\
Die Voraussetzungen des Pontrjaginschen Maximumprinzips an einen zul"assigen Steuerungsprozess sind in dem vorliegenden Beispiel genau dann erf"ullt,
wenn die Steuerung $u(\cdot)$ dem Raum $L_1(\R_+,[0,1])$ angeh"ort.
Denn in diesem Fall gelten $x(\cdot) \in C_{\lim}(\R_+,\R)$ f"ur die korrespondierende Trajektorie und die Bedingungen 
(\ref{PMPBedingung}), (\ref{PMPBedingung2}).
Wir diskutieren zwei optimale Steuerungsprozesse:
\begin{enumerate}
\item[(A)] In diesem Beispiel liefert der Steuerungsprozess
           $$x_*(t) = e^{\alpha t},\qquad z_*(t)= \frac{1}{\alpha - \varrho} (e^{(\alpha-\varrho)t}-1),\qquad
             u_*(t)=\alpha,\qquad \alpha=\varrho-\frac{1}{Z} \in (0, \varrho)$$
           ein globales Maximum.
           Da die vorgeschlagene Steuerung $u_*(\cdot)$ "uber $\R_+$ nicht integrierbar ist,
           gelten weder $x_*(\cdot) \in C_{\lim}(\R_+,\R)$ noch die Bedingungen (\ref{PMPBedingung}). \\
           Das Maximumprinzip ist auf diesen Steuerungsprozess nicht anwendbar.
\item[(B)] Ebenfalls stellt der Steuerungsprozess
           $$y_*(t) = \left\{ \begin{array}{ll} e^t,& t \in [0,\tau), \\ e^\tau, & t \in [\tau,\infty), \end{array} \right. \quad
             w_*(t) = \left\{ \begin{array}{ll} 1,& t \in [0,\tau), \\ 0, & t \in [\tau,\infty) \end{array} \right.$$
           mit dem Umschaltzeitpunkt $\tau >0$, der der Bedingung
           $$e^{(1-\varrho)\tau}\bigg(\frac{1}{\varrho}+\frac{1}{1-\varrho}\bigg) = Z+\frac{1}{1-\varrho}, \qquad Z> \frac{1}{\varrho},$$
           gen"ugt,
           ein globales Maximum dar.
           Die zugeh"orige Trajektorie $z_*(\cdot)$ lautet
           $$z_*(t) = \left\{ \begin{array}{ll} \frac{1}{1-\varrho}\big(e^{(1-\varrho)t} - 1 \big) ,& t \in [0,\tau), \\
                                                z(\tau) + \frac{1}{\varrho}\big(e^{(1-\varrho)\tau} - e^{\tau-\varrho t} \big), & t \in [\tau,\infty).
             \end{array} \right.$$
           Da die Steuerung $w_*(\cdot)$ dem Raum $L_1(\R_+,[0,1])$ angeh"ort,
           gelten s"amtliche Voraussetzungen von Theorem \ref{SatzPAUH}.
           Mit den Multiplikatoren
           $$\lambda_0=1, \qquad p(t)= e^{-\varrho t}, \qquad q(t)= \varrho-1$$
           sind dann die notwendigen Bedingungen (\ref{SatzPAUH1})--(\ref{SatzPAUH3}) erf"ullt. \hfill $\square$
\end{enumerate}}
\end{beispiel}
       \subsection{Der Nachweis der notwendigen Optimalit"atsbedingungen} \label{AbschnittBeweisPMPUH}
Es seien $x_*(\cdot) \in \mathscr{B}_{\rm Lip} \cap C_{\lim}(\R_+,\R^n)$ und $V_\gamma= \{ (t,x) \in \overline{\R}_+ \times \R^n \,|\, \|x-x(t)\| \leq \gamma\}$.
Wir betrachten f"ur $\big(x(\cdot),u(\cdot)\big) \in C_{\lim}(\R_+,\R^n) \times L_\infty(\R_+,\R^m)$ die Abbildungen
\begin{eqnarray*}
J\big(x(\cdot),u(\cdot)\big) &=& \int_0^\infty \omega(t)f\big(t,x(t),u(t)\big) \, dt, \\
F\big(x(\cdot),u(\cdot)\big)(t) &=& x(t) -x(t_0) -\int_0^t \varphi\big(s,x(s),u(s)\big) \, ds, \quad t \in \R_+,\\
H_0\big(x(\cdot)\big) &=& h_0\big(x(0)\big), \qquad H_1\big(x(\cdot)\big) = \lim_{t \to \infty} h_1\big(t,x(t)\big).
\end{eqnarray*}
Die Abbildungen fassen wir in folgenden Funktionenr"aumen auf:
\begin{eqnarray*}
J &:& C_{\lim}(\R_+,\R^n) \times L_\infty(\R_+,\R^m) \to \R, \\
F &:& C_{\lim}(\R_+,\R^n) \times L_\infty(\R_+,\R^m) \to C_{\lim}(\R_+,\R^n), \\
H_i &:& C_{\lim}(\R_+,\R^n) \to \R^{s_i}, \quad i=0,1.
\end{eqnarray*}

Wir setzen $\mathscr{F}=(F,H_0,H_1)$, sowie die Menge $\mathscr{U}$ gem"a"s
\begin{eqnarray*}
\mathscr{U}=\big\{ u(\cdot)  \in L_\infty(\R_+,U) &\big|& u(t)=u_*(t) + \chi_M(t)\big(w(t)-u_*(t)\big), \; w(\cdot) \in L_\infty(\R_+,U), \\
                                                  &     & M \subset \R_+ \mbox{ me"sbar und beschr"ankt} \big\}
\end{eqnarray*}
und pr"ufen f"ur die Extremalaufgabe
\begin{equation} \label{ExtremalaufgabePMPUH}
J\big(x(\cdot),u(\cdot)\big) \to \inf, \qquad \mathscr{F}\big(x(\cdot),u(\cdot)\big)=0, \qquad u(\cdot) \in \mathscr{U}
\end{equation}
im Punkt $\big(x_*(\cdot),u_*(\cdot)\big)$ die Voraussetzungen von Theorem \ref{SatzExtremalprinzipStark}:

\begin{enumerate}
\item[(A$_1$)] F"ur jedes $u(\cdot) \in \mathscr{U}$ ist die Abbildung $x(\cdot) \to J\big(x(\cdot),u(\cdot)\big)$
               nach Beispiel \ref{DiffZielfunktionalS} im Punkt $x_*(\cdot)$ Fr\'echet-differenzierbar.
\item[(A$_2$)] Die Abbildung $F$ ist die Summe der Abbildung $x(\cdot) \to x(t)$ und der Abbildung
               $$\big(x(\cdot),u(\cdot)\big) \to -\int_0^t \varphi\big(s,x(s),u(s)\big) \, ds.$$
               Im Beispiel \ref{DiffDynamikS} ist die Fr\'echet-Differenzierbarkeit der zweiten Abbildung im Punkt $x_*(\cdot)$
               f"ur jedes $u(\cdot) \in \mathscr{U}$ nachgewiesen.
               Da $x_*(\cdot)$ dem Raum $C_{\lim}(\R_+,\R^n)$ angeh"ort und $h_1(t,x)$ in $t=\infty$ stetig und stetig differenzierbar bez"uglich $x$ ist,
               sind die Abbildungen $H_i$ stetig differenzierbar.
\item[(B)] Nach Lemma \ref{FolgerungDGL} besitzt der Operator $\mathscr{F}_x\big(x_*(\cdot),u_*(\cdot)\big)$ eine endliche Kodimension.
\item[(C)] Der Nachweis dieser Voraussetzungen sind in Lemma \ref{LemmaNVUH1} und Lemma \ref{LemmaNVUH2} "uber mehrfache Nadelvariationen
           "uber dem unendlichen Zeithorizont dargestellt.
\end{enumerate}

Zur Extremalaufgabe (\ref{ExtremalaufgabePMPUH}) definieren wir auf
$$C_{\lim}(\R_+,\R^n) \times L_\infty(\R_+,\R^m) \times \R \times C_{\lim}^*(\R_+,\R^n) \times \R^{s_0} \times \R^{s_1}$$
die Lagrange-Funktion $\mathscr{L}=\mathscr{L}\big(x(\cdot),u(\cdot),\lambda_0,y^*,l_0,l_1\big)$,
$$\mathscr{L}= \lambda_0 J\big(x(\cdot),u(\cdot)\big)+ \big\langle y^*, F\big(x(\cdot),u(\cdot)\big) \big\rangle
                         +l_0^T H_0\big(x(\cdot)\big)+l_1^T H_1\big(x(\cdot)\big).$$
Ist $\big(x_*(\cdot),u_*(\cdot)\big)$ eine starke lokale Minimalstelle der Aufgabe (\ref{ExtremalaufgabePMPUH}),
dann existieren nach Theorem \ref{SatzExtremalprinzipStark}
nicht gleichzeitig verschwindende Lagrangesche Multiplikatoren 
$$\lambda_0 \geq 0, \qquad y^* \in C_{\lim}^*(\R_+,\R^n), \qquad l_0 \in \R^{s_0}, \qquad l_1 \in \R^{s_1}$$
derart,
dass gelten:
\begin{enumerate}
\item[(a)] Die Lagrange-Funktion besitzt bez"uglich $x(\cdot)$ in $x_*(\cdot)$ einen station"aren Punkt, d.\,h.
          \begin{equation}\label{SatzPMPLMRUH1}
          \mathscr{L}_x\big(x_*(\cdot),u_*(\cdot),\lambda_0,y^*,l_0,l_1\big)=0;
          \end{equation}         
\item[(b)] Die Lagrange-Funktion erf"ullt bez"uglich $u(\cdot)$ in $u_*(\cdot)$ die Minimumbedingung
           \begin{equation}\label{SatzPMPLMRUH2}
           \mathscr{L}\big(x_*(\cdot),u_*(\cdot),\lambda_0,y^*,l_0,l_1\big)
           = \min_{u(\cdot) \in \mathscr{U}} \mathscr{L}\big(x_*(\cdot),u(\cdot),\lambda_0,y^*,l_0,l_1\big).
           \end{equation}
\end{enumerate}
Aufgrund (\ref{SatzPMPLMRUH1}) ist folgende Variationsgleichung f"ur alle $x(\cdot) \in C_{\lim}(\R_+,\R^n)$ erf"ullt: 
\begin{eqnarray*}
0 &=& \lambda_0 \cdot \int_0^\infty \omega(t) \big\langle f_x\big(t,x_*(t),u_*(t)\big),x(t) \big\rangle\, dt \\
  & & + \int_0^\infty \bigg[ x(t)-x(0) - \int_0^t \varphi_x\big(s,x_*(s),u_*(s)\big) x(s) \,ds \bigg]^T  d\mu(t) \\
  & & + \big\langle l_0, h_0'\big(x_*(0)\big) x(0) \big\rangle + \big\langle l_1, \lim_{t \to \infty} h_{1x}\big(t,x_*(t)\big) x(t) \big\rangle.
\end{eqnarray*}
Dabei ist $\mu \in \mathscr{M}(\overline{\R}_+)$ und besitzt nach Satz \ref{SatzRieszClim} die Darstellung $\mu=\mu_0+\mu_\infty$
mit $\mu_0 \in \mathscr{M}(\R_+)$ und einem in $t = \infty$ konzentrierten signierten Ma"s $\mu_\infty$. \\
Wir "andern die Integrationsreihenfolge und setzen $p(t)= \displaystyle \int_t^\infty \, d\mu(s)$.
Dann gilt
$$\lim_{t \to \infty} p(t) = \mu_\infty.$$
Aus der eindeutigen Darstellung eines stetigen linearen Funktionals im Raum $C_{\lim}(\R_+,\R^n)$ folgen (\ref{SatzPAUH1}) und (\ref{SatzPAUH2}). \\[2mm]
Gem"a"s (\ref{SatzPMPLMRUH2}) gilt f"ur alle $u(\cdot) \in \mathscr{U}$ die Ungleichung
$$\int_0^\infty H\big(t,x_*(t),u_*(t),p(t),\lambda_0\big) \, dt \geq
  \int_0^\infty H\big(t,x_*(t),u(t),p(t),\lambda_0\big) \, dt.$$
Daraus folgt abschlie"send durch Standardtechniken f"ur Lebesguesche Punkte die Maximumbedingung (\ref{SatzPAUH3}).
Der Beweis von Theorem \ref{SatzPAUH} ist abgeschlossen. \hfill $\blacksquare$
       \subsection{Zur normalen Form und zu Transversalit"atsbedingungen}
F"ur die Aufgabe (\ref{PAUH1})--(\ref{PAUH4}) lassen sich Aussagen "uber die Normalform des Pontrjaginschen Maximumprinzips und zu diversen
Transversalit"atsbedingungen ableiten. \\[2mm]
Wir betrachten zun"achst die Aufgabe (\ref{PAUH1})--(\ref{PAUH4}) mit freiem rechten Endpunkt im Unendlichen.
Dann gen"ugt die Adjungierte $p(\cdot)$ nach Theorem \ref{SatzPAUH} dem Randwertproblem
$$\dot{p}(t) = -\varphi_x^T\big(t,x_*(t),u_*(t)\big) p(t) + \lambda_0 \omega(t)f_x\big(t,x_*(t),u_*(t)\big), \qquad
  \lim_{t \to \infty} p(t)=0.$$
  
\begin{folgerung} \label{FolgerungPAUH1}
In der Aufgabe (\ref{PAUH1})--(\ref{PAUH4}) mit freiem rechten Endpunkt im Unendlichen
ergeben sich aus dem Maximumprinzip unmittelbar die ``nat"urlichen'' Transversalit"atsbedingungen\index{Transversalit"atsbedingungen!nat@--, nat"urliche}:
$$\lim_{t \to \infty} p(t) =0, \qquad \lim_{t \to \infty} \langle p(t),x(t) \rangle = 0 \mbox{ f"ur alle } x(\cdot) \in W^1_\infty(\R_+,\R^n).$$
\end{folgerung}

Nach Voraussetzung des Maximumprinzips ist
$$\int_0^\infty \big\|\varphi_x\big(t,x_*(t),u_*(t)\big)\big\| \, dt < \infty.$$
Unter der Annahme $\lambda_0=0$ w"urde die Adjungierte als L"osung der Gleichung
$$\dot{p}(t) = -\varphi_x^T\big(t,x_*(t),u_*(t)\big) p(t), \qquad \lim_{t \to \infty} p(t)=0,$$
nach Folgerung \ref{FolgerungDGL2} im Widerspruch zu Theorem \ref{SatzPAUH} identisch verschwinden.

\begin{folgerung} \label{FolgerungPAUH2}
In der Aufgabe (\ref{PAUH1})--(\ref{PAUH4}) mit freiem rechten Endpunkt im Unendlichen gilt $\lambda_0 \not= 0$ und wir k"onnen ohne Einschr"ankung
$\lambda_0=1$ annehmen.
\end{folgerung}

Wegen der Integrierbarkeit der Abbildung $t \to \varphi_x\big(t,x_*(t),u_*(t)\big)$ "uber $\R_+$ k"onnen wir die 
die in $t=0$ normalisierten Fundamentalmatrizen $Y_*(t)$ bzw. $Z_*(t)$ der homogenen Systeme
$$\dot{y}(t)=\varphi_x\big(t,x_*(t),u_*(t)\big) y(t), \qquad \dot{z}(t)=-\varphi_x^T\big(t,x_*(t),u_*(t)\big) z(t)$$
im Rahmen des Raumes $C_{\lim}(\R_+,\R^n)$ betrachten.
Es bezeichne ferner $y_\xi(\cdot) \in C_{\lim}(\R_+,\R^n)$ die Funktion $y_\xi(t)=Y_*(t) Y^{-1}_*(T) \xi$ mit $\xi \in \R^n$ und $\|\xi\|=1$.
Dann ergibt sich auf die gleiche Weise wie in Tauchnitz \cite{TauchnitzWMPIHOC,TauchnitzOC} die Beziehung
$$\langle p(t), y_\xi(t) \rangle = \Big\langle  p(T) + Z_*(T)\int_T^t \omega(s) Z_*^{-1}(s)f_x\big(s,x_*(s),u_*(s)\big) ds , \xi \Big\rangle.$$
Unter Verwendung der ``nat"urlichen'' Transversalit"atsbedingungen in Folgerung \ref{FolgerungPAUH1} und wegen der Willk"urlichkeit von $\xi$
erhalten wir daraus die Integraldarstellung
der Arbeiten von Aseev \& Kryazhimskii und Aseev \& Veliov \cite{AseKry,AseVel,AseVel2,AseVel3}:

\begin{folgerung} \label{FolgerungPAUH3}
Es gen"uge $\big(x_*(\cdot),u_*(\cdot)\big)$ den Voraussetzungen des Pontrjaginschen Maximumprinzips \ref{SatzPAUH}.
Ist $\big(x_*(\cdot),u_*(\cdot)\big)$ ein starkes lokales Minimum der Aufgabe (\ref{PAUH1})--(\ref{PAUH4}) mit freiem rechten Endpunkt im Unendlichen,
dann besitzt die Adjungierte $p(\cdot)$ die Darstellung\index{Adjungierte!eindeutig@--, eindeutige Darstellung}
$$p(t)= -Z_*(t) \int_t^\infty \omega(s) Z^{-1}_*(s) f_x\big(s,x_*(s),u_*(s)\big) \, ds.$$
Dabei ist $Z_*(t)$ die in $t=0$ normalisierte Fundamentalmatrix des linearen Systems
$$\dot{z}(t)=-\varphi^T_x\big(t,x_*(t),u_*(t)\big) z(t).$$
\end{folgerung}

In der Aufgabe (\ref{PAUH1})--(\ref{PAUH4}) seien nun gewisse Randwerte im Unendlichen explizit gegeben, d.\,h. $h_1\big(t,x(t)\big)=x(t)-x_1$.
Wir schlie"sen dabei nicht aus,
dass dabei gewisse Komponenten von $x_1$ nicht fest vorgegeben, sondern ohne Einschr"ankung sind.
Genauer bedeutet dies f"ur $x(t)=\big(x_1(t),...,x_n(t)\big)$:
$$\lim_{t \to \infty} x_i(t)=x_i \in \R, \qquad \lim_{t \to \infty} x_j(t) \mbox{ frei}, \qquad i=1,...,l,\; j=l+1,...,n.$$
Daher verschwinden einerseits bei expliziten Randwerten im Unendlichen die entsprechenden Komponenten $x_i(t)-x_{i*}(t)$ f"ur $t \to \infty$, $i=1,...,l$.
F"ur diejenigen Komponenten, f"ur die die Randwerte im Unendlichen frei sind,
verschwinden die entsprechenden Komponenten $p_j(t)$ der Adjungierten f"ur $t \to \infty$, $j=l+1,...,n$.
Damit k"onnen wir festhalten:

\begin{folgerung} \label{FolgerungPAUH4}
In der Aufgabe (\ref{PAUH1})--(\ref{PAUH4}) mit expliziten Randbedingungen im Unendlichen
ergibt sich aus dem Maximumprinzip unmittelbar die ``nat"urliche'' Transversalit"atsbedingung\index{Transversalit"atsbedingungen!nat@--, nat"urliche}:
$$\lim_{t \to \infty} \langle p(t),x(t)-x_*(t) \rangle = 0 \quad \mbox{ f"ur alle zul"assigen } x(\cdot) \in W^1_\infty(\R_+,\R^n).$$
\end{folgerung}

In der Aufgabe (\ref{PAUH1})--(\ref{PAUH4}) mit freiem rechten Endpunkt verschwinde die Verteilungsfunktion $\omega(\cdot) \in L_1(\R_+,\R_+)$
im Unendlichen.
Dann folgen aus der gleichm"a"sigen Stetigkeit der Abbildungen $f,\,\varphi$ und aus der Beschr"anktheit des Steuerungsprozesses $\big(x_*(\cdot),u_*(\cdot)\big)$ 
die Grenzwerte
$$\lim_{t \to \infty} \omega(t) f\big(t,x_*(t),u_*(t)\big)=0, \qquad
  \lim_{t \to \infty} \big\langle p(t)\,,\, \varphi\big(t,x_*(t),u_*(t)\big)\big\rangle=0.$$

\begin{folgerung} \label{FolgerungPAUH5}
In der Aufgabe (\ref{PAUH1})--(\ref{PAUH4}) mit freiem Endpunkt im Unendlichen besitze die Verteilungsfunktion $\omega(\cdot)$ einen Grenzwert im
Unendlichen, d.\,h.
$$\lim_{t \to \infty} \omega(t)=0.$$
Dann ergibt sich die Bedingung von Michel \index{Transversalit"atsbedingungen!von@-- von Michel}:
$$\lim_{t \to \infty} H\big(t,x_*(t),u_*(t),p(t),\lambda_0\big)=0.$$
\end{folgerung}

%Abschlie"send betrachten wir die Aufgabe (\ref{PAUH1})--(\ref{PAUH4}) mit freiem Endpunkt im Unendlichen bez"uglich der Verteilungsfunktion
%$\omega(t)=e^{-\varrho t}$ und es seien die eingehenden Abbildungen autonom, d.\,h. es gilt
%$$f=f(x,u), \qquad \varphi=\varphi(x,u).$$

       \subsection{Hinreichende Bedingungen nach Arrow} \label{AbschnittArrowPMP}
Unser\index{hinreichende Bedingungen!Arrow@-- nach Arrow} Vorgehen zur Herleitung der hinreichenden Bedingungen basiert wieder auf der Darstellung in 
Seierstad \& Syds\ae ter \cite{Seierstad} und Aseev \& Kryazhimskii \cite{AseKry}.
Wir betrachten das Steuerungsproblem
\begin{eqnarray}
&& \label{HBPAUH1} J\big(x(\cdot),u(\cdot)\big) = \int_0^\infty \omega(t)f\big(t,x(t),u(t)\big) \, dt \to \inf, \\
&& \label{HBPAUH2} \dot{x}(t) = \varphi\big(t,x(t),u(t)\big), \\
&& \label{HBPAUH3} x(0)=x_0, \qquad \lim_{t \to \infty} x(t)=x_1, \\
&& \label{HBPAUH4} u(t) \in U \subseteq \R^m, \quad U \not= \emptyset.
\end{eqnarray}
In der Aufgabenstellung schlie"sen wir wieder den Fall nicht aus,
dass durch die Randbedingungen (\ref{HBPAUH3}) gewisse Komponenten der Punkte $x_0$ und $x_1$ nicht fest vorgegeben,
sondern ohne Einschr"ankung sind. \\[2mm]
Wir definieren die Mengen $V_\gamma$ und $V_\gamma(t)$:
$$V_\gamma=\{ (t,x) \in \overline{\R}_+ \times \R^n \,|\, \|x-x(t)\| \leq \gamma\}, \quad
  V_\gamma(t)=\{ x \in \R^n \,|\, \|x-x_*(t)\| \leq \gamma\}.$$
Au"serdem bezeichnet $\mathscr{H}(t,x,p) = \sup\limits_{u \in U} H(t,x,u,p,1)$ die Hamilton-Funktion $\mathscr{H}$ im normalen Fall.

\begin{theorem} \label{SatzHBPMPUH}
In der Aufgabe (\ref{HBPAUH1})--(\ref{HBPAUH4}) sei $\big(x_*(\cdot),u_*(\cdot)\big) \in \mathscr{B}_{\rm Lip} \cap \mathscr{B}_{\rm adm}$ 
und es sei $p(\cdot):\R_+ \to \R^n$. Ferner gelte:
\begin{enumerate}
\item[(a)] Das Tripel $\big(x_*(\cdot),u_*(\cdot),p(\cdot)\big)$
           erf"ullt (\ref{SatzPAUH1})--(\ref{SatzPAUH3}) in Theorem \ref{SatzPAUH} mit $\lambda_0=1$.        
\item[(b)] F"ur jedes $t \in \R_+$ ist die Funktion $\mathscr{H}\big(t,x,p(t)\big)$ konkav in $x$ auf $V_\gamma(t)$.
\end{enumerate}
Dann ist $\big(x_*(\cdot),u_*(\cdot)\big)$ ein starkes lokales Minimum der Aufgabe (\ref{HBPAUH1})--(\ref{HBPAUH4}).
\end{theorem}

{\bf Beweis} Es sei $t \in \R_+$ gegeben.
Da die Abbildung $x \to \mathscr{H}\big(t,x,p(t)\big)$ auf $V_\gamma(t)$ konkav ist, ist die Menge
$$Z=\big\{ (\alpha,x) \in \R \times \R^n \,\big|\, x \in V_\gamma(t), \alpha \leq \mathscr{H}\big(t,x,p(t)\big) \big\}$$
konvex und besitzt ein nichtleeres Inneres.
Ferner ist $\big(\alpha_*,x_*(t)\big)$ mit $\alpha_*= \mathscr{H}\big(t,x_*(t),p(t)\big)$ ein Randpunkt der Menge $Z$.
Daher existiert ein nichttrivialer Vektor $\big(a_0(t),a(t)\big) \in \R \times \R^n$ mit
\begin{equation} \label{BeweisHBPMP1}
a_0(t) \alpha + \langle a(t),x\rangle \leq a_0(t) \alpha_* + \langle a(t),x_*(t)\rangle \quad \mbox{ f"ur alle } (\alpha,x) \in Z.
\end{equation}
Es ist $x_*(t)$ ein innerer Punkt der Menge $V_\gamma(t)$.
Weiterhin folgt aus den elementaren Eigenschaften konkaver Funktionen,
dass $x \to \mathscr{H}\big(t,x,p(t)\big)$ in $V_\gamma(t)$ stetig ist,
da sie auf $V_\gamma(t)$ konkav und nach unten durch $H\big(t,x,u_*(t),p(t),1\big)$ beschr"ankt ist. \\
Deswegen existiert ein $\delta>0$ mit $x_*(t)+\xi \in V_\gamma(t)$ und $\big(\alpha_*-1,x_*(t)+\xi\big) \in Z$
f"ur alle $\|\xi\| \leq \delta$.
Aus (\ref{BeweisHBPMP1}) folgt daher $\langle a(t),\xi\rangle - a_0(t) \leq 0$ f"ur alle $\|\xi\| \leq \delta$.
Dies zeigt $a_0(t) >0$ und wir k"onnen ohne Einschr"ankung $a_0(t)=1$ annehmen.
Wiederum (\ref{BeweisHBPMP1}) liefert damit
\begin{equation} \label{BeweisHBPMP2}
\langle a(t),x-x_*(t)\rangle \leq \mathscr{H}\big(t,x_*(t),p(t)\big) - \mathscr{H}\big(t,x,p(t)\big)
  \quad \mbox{ f"ur alle } x \in V_\gamma(t).
\end{equation}
Es sei nun $t \in \R_+$ so gew"ahlt,
dass die Maximumbedingung (\ref{SatzPAUH3}) zu diesem Zeitpunkt erf"ullt ist.
Dann folgt aus (\ref{BeweisHBPMP2}), dass
\begin{eqnarray*}
      -\langle a(t),x-x_*(t)\rangle
&\geq& \mathscr{H}\big(t,x,p(t)\big) - \mathscr{H}\big(t,x_*(t),p(t)\big) \\
&=&    \sup_{u \in U} H\big(t,x,u,p(t),1\big) - \mathscr{H}\big(t,x_*(t),p(t)\big) \\
&\geq& \big\langle p(t), \varphi\big(t,x,u_*(t)\big) \big\rangle - f\big(t,x,u_*(t)\big) \\
&    & -\big[ \big\langle p(t), \varphi\big(t,x_*(t),u_*(t)\big) \big\rangle - f\big(t,x_*(t),u_*(t)\big)
\end{eqnarray*}
f"ur alle $x \in V_\gamma(t)$ gilt.
Wir setzen
\begin{eqnarray*}
\Phi(x) &=& \big\langle p(t), \varphi\big(t,x,u_*(t)\big)-\varphi\big(t,x_*(t),u_*(t)\big) \big\rangle \\
        & & - \big[f\big(t,x,u_*(t)\big)-f\big(t,x_*(t),u_*(t)\big)\big] + \langle a(t),x-x_*(t)\rangle.
\end{eqnarray*}
Die Funktion $\Phi(x)$ ist stetig differenzierbar auf $V_\gamma(t)$.
Ferner gilt $\Phi(x) \leq 0$ f"ur alle $x \in V_\gamma(t)$ und $\Phi(x_*(t))=0$.
Damit nimmt die Funktion $\Phi$ in dem inneren Punkt $x_*(t)$ der Menge $V_\gamma(t)$ ihr absolutes Maximum an.
Also gilt
$$0=\Phi'(x_*(t)) = \varphi_x^T\big(t,x_*(t),u_*(t)\big) p(t) - f_x\big(t,x_*(t),u_*(t)\big) + a(t)$$
bzw.
\begin{equation} \label{BeweisHBPMP3}
a(t)= -\varphi_x^T\big(t,x_*(t),u_*(t)\big) p(t) + f_x\big(t,x_*(t),u_*(t)\big).
\end{equation}
Die Gleichung (\ref{BeweisHBPMP3}) wurde unter der Annahme erzielt,
dass die Maximumbedingung (\ref{SatzPAUH3}) in dem Zeitpunkt $t \in \R_+$ erf"ullt ist.
Da (\ref{SatzPAUH3}) f"ur fast alle $t \in \R_+$ gilt,
stimmt $a(t)$ mit der verallgemeinerten Ableitung $\dot{p}(t)$ "uberein.
Also gilt auf $V_\gamma$ f"ur fast alle $t \in \R_+$ die Ungleichung
\begin{equation} \label{BeweisHBPMP4}
\langle \dot{p}(t),x-x_*(t)\rangle \leq \mathscr{H}\big(t,x_*(t),p(t)\big)- \mathscr{H}\big(t,x,p(t)\big).
\end{equation}
Es ergibt sich also zu $T \in \R_+$ die Beziehung
\begin{eqnarray*}
    \Delta(T)
&=& \int_0^T \omega(t)\big[f\big(t,x(t),u(t)\big)-f\big(t,x_*(t),u_*(t)\big)\big] \, dt \\
&\geq & \int_0^T\big[\mathscr{H}\big(t,x_*(t),p(t)\big)-\mathscr{H}\big(t,x(t),p(t)\big)\big] \, dt 
        + \int_0^T \langle p(t), \dot{x}(t)-\dot{x}_*(t) \rangle \, dt \\
&\geq& \langle p(T),x(T)-x_*(T)\rangle-\langle p(0),x(0)-x_*(0)\rangle.
\end{eqnarray*}

Im Fall fester Anfangs- und Endbedingungen verschwinden die Differenzen $x(0)-x_*(0)$ und $x(T)-x_*(T)$ f"ur $T \to \infty$.
Sind jedoch gewissen Komponenten im Anfangs- oder Endpunkt $x_0$ bzw. $x_1$ frei,
dann liefern die Transversalit"atsbedingungen,
dass die entsprechenden Komponenten der Adjungierten $p(\cdot)$ zum Zeitpunkt $t=0$ bzw. im Unendlichen verschwinden.
Daher folgt die Beziehung
$$\lim_{T \to \infty} \Delta(T) \geq \lim_{T \to \infty} \langle p(T),x(T)-x_*(T)\rangle-\langle p(0),x(0)-x_*(0)\rangle=0$$
f"ur alle zul"assigen $\big(x(\cdot),u(\cdot)\big)$ mit $\|x(\cdot)-x_*(\cdot)\|_\infty \leq \gamma$. \hfill $\blacksquare$

\begin{bemerkung}{\rm
Die Herleitung der hinreichenden Bedingungen in Theorem \ref{SatzHBPMPUH} basiert im Wesentlichen auf der Konkavit"at der Hamilton-Funktion $\mathscr{H}$.
Daher m"ussen bei der Anwendung von Theorem \ref{SatzHBPMPUH} die Einschr"ankungen (\ref{PMPBedingung}) und (\ref{PMPBedingung2}) nicht gelten.
Dies hat weiterhin zur Folge, dass im Fall des freien rechten Endpunktes die Trajektorie $x_*(\cdot)$ keinen Grenzwert im Unendlichen besitzen muss.
\hfill $\square$}
\end{bemerkung}

\begin{beispiel} {\rm Wir betrachten die Aufgabe
\begin{eqnarray*}
&& J\big(x(\cdot),u(\cdot)\big) = \int_0^\infty e^{-\varrho t} \big(1-u(t)\big) x(t) \, dt \to \sup,\\
&& \dot{x}(t)=u(t)x(t), \quad x(0)=1, \quad \lim_{t \to \infty} x(t)=x_1>1, \quad u(t) \in [0,1], \quad \varrho \in (0,1).
\end{eqnarray*}
Aus den notwendigen Bedingungen (\ref{SatzPAUH1})--(\ref{SatzPAUH3}) ergeben sich
\begin{eqnarray*}
x_*(t) &=& \left\{ \begin{array}{ll} e^t, & t \in [0,\tau), \\ x_1, & t \in [\tau, \infty), \end{array}\right. \quad
  u_*(t)= \left\{ \begin{array}{ll} 1, & t \in [0,\tau), \\ 0, & t \in [\tau, \infty), \end{array}\right. \quad \tau=\ln x_1, \\
p(t) &=& \left\{ \begin{array}{ll}
        e^{(1-\varrho)\tau} e^{-t}, & t \in [0,\tau), \\
        \frac{\varrho-1}{\varrho}e^{-\varrho \tau} + \frac{1}{\varrho} e^{-\varrho t}, & t \in [\tau, \infty).
        \end{array}\right.
\end{eqnarray*}
Offenbar ist die Hamilton-Funktion $\mathscr{H}$ konkav in $x$ und damit $\big(x_*(\cdot),u_*(\cdot)\big)$
ein starkes lokales Minimum der Aufgabe.  \hfill $\square$}
\end{beispiel}

\begin{beispiel} {\rm F"ur den linear-quadratische Regler
\begin{eqnarray*}
&& J\big(x(\cdot),u(\cdot)\big) = \int_0^\infty e^{-2t} \cdot \frac{1}{2}\big( x^2(t)+u^2(t)\big) \, dt \to \inf, \\
&& \dot{x}(t) = 2 x(t)+u(t), \qquad x(0)=2, \qquad u(t) \in \R
\end{eqnarray*}
liefern die Bedingungen (\ref{SatzPAUH1})--(\ref{SatzPAUH3}) den Steuerungsprozess und die Adjungierte
$$x_*(t)=2e^{(1-\sqrt{2})t}, \quad u_*(t)=-2(1+\sqrt{2})e^{(1-\sqrt{2})t}, \quad p(t)=e^{-2t}u_*(t).$$
Die Hamilton-Funktion $\mathscr{H}$ ist konkav bez"uglich $x$.
Damit ist $\big(x_*(\cdot),u_*(\cdot)\big)$ ein starkes lokales Minimum. \hfill $\square$}
\end{beispiel}

\begin{beispiel}
{\rm Im Beispiel \ref{BeipielRWUnendlich} mit Budgetbeschr"ankung,
\begin{eqnarray*}
&& J\big(x(\cdot),z(\cdot),u(\cdot)\big)=\int_0^\infty e^{-\varrho t}\big(1-u(t)\big)x(t) \, dt \to \sup, \\
&& \dot{x}(t)=u(t)x(t), \; x(0)=1,\qquad \dot{z}(t)=e^{-\varrho t}x(t), \; z(0)=0, \; \lim_{t \to \infty} z(t) = Z, \\
&& u \in [0,1], \qquad \varrho \in (0,1),
\end{eqnarray*}
ist jeder zul"assige Steuerungsprozess $\big(x(\cdot),z(\cdot),u(\cdot)\big)$ global optimal.
Die Voraussetzungen in Theorem \ref{SatzHBPMPUH} an einen zul"assigen Steuerungsprozess sind genau dann erf"ullt,
wenn die Steuerung $u(\cdot)$ dem Raum $L_1(\R_+,[0,1])$ angeh"ort.
Weiterhin sind f"ur jeden zul"assigen Steuerungprozess mit den Multiplikatoren
$$\lambda_0=1, \qquad p(t)= e^{-\varrho t}, \qquad q(t)= \varrho-1$$
die notwendigen Bedingungen (\ref{SatzPAUH1})--(\ref{SatzPAUH3}) erf"ullt.
Weiterhin ist die Hamilton-Funktion $\mathscr{H}$ offenbar konkav bez"uglich $(x,z)$ und es gilt Theorem \ref{SatzHBPMPUH}.} \hfill $\square$
\end{beispiel}

\begin{beispiel}[Abbau einer erneuerbaren Ressource] \label{ExampleErneuRessource}\index{Ressourcenabbau}
{\rm Wir betrachten die Aufgabe
\begin{eqnarray}
&& \label{ErneuRessource1} J\big(x(\cdot),u(\cdot)\big)
   =\int_0^\infty e^{-\varrho t}\big[\pi\big(x(t)\big)-\kappa\big(x(t)\big)\big] u(t) \, dt \to \sup, \\
&& \label{ErneuRessource2} \dot{x}(t) = G\big(x(t)\big) - u(t), \quad x(0)=x_0>0,\quad u(t) \in [0, u_{\max}].
\end{eqnarray}
"Okonomische Interpretation: Durch die Funktion $G$ wird die nat"urliche dynamische Entwicklung der Ressource (oder einer Population) beschrieben,
wie diese sich ohne externe Einfl"usse entwickelt.
Diese nat"urliche Entwicklung wird durch die Abbaurate $u$ gest"ort.
Entsprechend des Angebotes $u$ wird ein Gewinn in H"ohe des Preises $\pi$ abz"uglich der Kosten $\kappa$ erzielt.
Dabei nehmen wir, dass Preis und Kosten steigen desto seltener die Ressource ist. \\[2mm]
In vielen Aufgaben dieser Form zeigt die L"osung ein station"ares Verhalten bzw. es tendiert die L"osung zu einem station"aren Paar $(\overline{x},\overline{u})$.
Eine solches station"ares Paar in Verbindung mit der Adjungierten $p(\cdot)$ nennen wir eine Gleichgewichtsl"osung\index{Gleichgewichtsl"osung}.
Wir identifizieren eine Gleichgewichtsl"osung $\big(\overline{x},\overline{u},p(\cdot)\big)$ durch folgende Eigenschaften:
$$\varphi(\overline{x},\overline{u}) = 0, \quad H_u\big(t,\overline{x},\overline{u},1,p(t)\big) =0, \quad
  \dot{p}(t) = -G'(\overline{x}) p(t) - e^{-\varrho t}\big(\pi'(\overline{x})-\kappa'(\overline{x})\big)\overline{u}.$$
Aus den ersten beiden Beziehung ergeben sich f"ur die Gleichgewichtsl"osung $\big(\overline{x},\overline{u},p(\cdot)\big)$:
$$G(\overline{x}) = \overline{u}, \qquad p(t)= e^{-\varrho t}\big(\pi(\overline{x})-\kappa(\overline{x})\big), \qquad \dot{p}(t)=-\varrho p(t).$$
Mit der Bezeichung $\gamma(x)= \pi(x)-\kappa(x)$ f"ur den Gewinn pro Einheit erhalten wir ferner in der adjungierten Gleichung:
\begin{eqnarray*}
\dot{p}(t)= -\varrho p(t) &=& -G'(\overline{x}) p(t) - e^{-\varrho t} \gamma'(\overline{x})\overline{u}, \\
-\varrho e^{-\varrho t} \gamma(\overline{x}) &=& -G'(\overline{x}) e^{-\varrho t} \gamma(\overline{x}) - e^{-\varrho t} \gamma'(\overline{x})G(\overline{x}), \\
0 &=& \big(\varrho -G'(\overline{x})\big) \gamma(\overline{x}) - \gamma'(\overline{x})G(\overline{x})
\end{eqnarray*}
oder ausgedr"uckt in der Form der Gewinnelastizit"at
$$\varepsilon_{\gamma,x}:=\overline{x} \frac{\gamma'(\overline{x})}{\gamma(\overline{x})} =\overline{x} \frac{\varrho -G'(\overline{x})}{G(\overline{x})}.$$
F"ur das logistische Wachstum $G(x)=x(r-Kx)$ ergibt sich die Gewinnelastitzit"at
$$\varepsilon_{\gamma,x} = -\frac{(r-K\overline{x})-(\varrho -K\overline{x})}{r-K\overline{x}} \in (-1,0)$$
im Fall einer initialen Wachstumsrate $G'(0)=r > \varrho$. \hfill $\square$}
\end{beispiel}

\begin{bemerkung}{\rm
Die Darstellung der Pontrjagin-Funktion und der notwendigen Bedingungen im Maximumprinzip \ref{SatzPAUH} nennt man auch
Gegenwartswert-Schreibweise (``present value''). \index{Pontrjaginsches Maximumprinzip!oekonomische@--, Gegenwartswert-Schreibweise}
Durch die Einf"uhrung Momentanwert-Schreibweise (``current value'') \index{Pontrjaginsches Maximumprinzip!oekonomische@--, Momentanwert-Schreibweise}
mittels der Funktionen $q=e^{\varrho t}p$ und $\tilde{H}=e^{\varrho t}H$,
\begin{eqnarray*}
e^{\varrho t} H(t,x,u,p,1) &=& e^{\varrho t} \big[\langle p, \varphi(t,x,u) \rangle- e^{- \varrho t} f(t,x,u)\big] \\
                           &=& \langle q, \varphi(t,x,u) \rangle - f(t,x,u) = \tilde{H}(t,x,u,q,1),
\end{eqnarray*}
lassen sich im letzten Beispiel die Gleichgewichtsbedingungen in folgende Form "uberf"uhren:
$$\varphi(\overline{x},\overline{u}) = 0, \quad \tilde{H}_u(t,\overline{x},\overline{u},1,\overline{q}) =0, \quad
  \varrho \overline{q} - \tilde{H}_x(t,\overline{x},\overline{u},1,\overline{q})=0.$$
In der Momentanwert-Schreibweise kommt der Charakter eines Gleichgewichts vollst"andig zum Ausdruck.
Deswegen findet man in den Untersuchungen zu Gleichgewichten in der Literatur h"aufig das Maximumprinzip in der Momentanwert-Schreibweise. \hfill $\square$}
\end{bemerkung}

\begin{beispiel} \index{Kapitalakkumulation}
{\rm Im Folgenden befassen wir uns mit einem neoklassischen Modell der "okonomischen Wachstumstheorie.
Diese Problemklasse geht auf die bereits erw"ahnte Arbeit von Ramsey \cite{Ramsey} zur"uck.
Es bezeichne $K(\cdot)$ das Verm"ogen einer "Okonomie,
$Y(\cdot)$ das Nationaleinkommen und $C(\cdot)$ die Konsumption.
Weiterhin wird zu jedem Zeitpunkt das Nationaleinkommen $Y$ in Konsumption und Investition aufgeteilt, d.\,h. $Y(t)=C(t)+\dot{K}(t)$.
Mit der Nutzenfunktion $U$ ergibt sich die allgemeine Aufgabenstellung
$$J\big(K(\cdot),C(\cdot)\big) = \int_0^\infty e^{-\varrho t} U\big(C(t)\big) \, dt \to \sup, \quad \dot{K}(t) = Y(t)-C(t).$$
Das Nationaleinkommen bestimmt sich durch die produzierten G"utern,
die aus dem eingesetzen Kapital $K$ und den Arbeitsressourcen $L$ gewonnen werden.
Mit der Cobb-Douglas-Produktionsfunktion $F(K,L)$ \index{Funktion, absolutstetige!Cobb@--, Cobb-Douglas-Produktions-} \index{Cobb-Douglas-Produktionsfunktion}
erhalten wir die Darstellung
$$Y(t)=F\big(K(t),L(t)\big), \qquad F(K,L)=K^\alpha L^{1-\alpha}, \qquad \alpha \in (0,1).$$
Die Verteilung des Nationaleinkommens $Y$ in Investition und Konsumption beschrieben mit $u \in [0,1]$ liefert
$$\dot{K}(t) = u(t) F\big(K(t),L(t)\big), \qquad C(t)=\big(1-u(t)\big) F\big(K(t),L(t)\big).$$
Mit einer logarithmischen Nutzenfunktion $U$ stellt sich so die Aufgabe wie folgt dar:
\begin{eqnarray*}
&& J\big(K(\cdot),C(\cdot)\big) = \int_0^\infty e^{-\varrho t} \ln\big[ \big(1-u(t)\big) F\big(K(t),L(t)\big) \big] \, dt \to \sup, \\
&& \dot{K}(t) = u(t) F\big(K(t),L(t)\big), \qquad \dot{L}(t) = \mu L(t), \qquad u(t) \in [0,1].
\end{eqnarray*}
Wir f"uhren die Kapital- und Konsumintensit"at $k = K/L$, $c=C/L$ und die Produktionsfunktion $f(x)=F(x,1)$ ein.
Dann ergeben sich die Beziehungen
\begin{eqnarray*}
\dot{k}(t) &=& \frac{\dot{K}(t) L(t)- K(t)\dot{L}(t)}{L^2(t)} = \frac{u(t) K^\alpha(t) L^{1-\alpha}(t) L(t)- K(t)\dot{L}(t)}{L^2(t)} \\
           &=& F\big(K(t)/L(t),1\big) - \frac{\dot{L}(t)}{L(t)}\frac{K(t)}{L(t)} =u(t)f\big(k(t)\big)-\mu k(t), \\
c(t)       &=& \big(1-u(t)\big) \frac{F\big(K(t),L(t)\big)}{L(t)} = \big(1-u(t)\big) F\big(K(t)/L(t),1\big) = \big(1-u(t)\big)f\big(k(t)\big).
\end{eqnarray*}
Auf diese Weise bekommt die Aufgabe die finale Form
\begin{eqnarray*}
&& J\big(k(\cdot),c(\cdot)\big) = \int_0^\infty e^{-\varrho t} \ln\big[ \big(1-u(t)\big) f\big(k(t)\big) \big] \, dt \to \sup, \\
&& \dot{k}(t) = u(t) f\big(k(t)\big)-\mu k(t), \qquad u(t) \in [0,1].
\end{eqnarray*}
Die Pontrjagin-Funktion der Aufgabe in der Momentanwert-Schreibweise ist
$$\tilde{H}(t,k,u,q,1) = q[u f(k)-\mu k] + \ln[(1-u) f(k)].$$
Dann ergeben die Gleichgewichtsbedingungen in der Momentanwert-Schreibweise:
\begin{enumerate}
\item[(a)] Das dynamische Gleichgewicht liefert
           $$\varphi(\overline{k},\overline{u}) = \overline{u} f(\overline{k})-\mu \overline{k} = 0
             \qquad\Leftrightarrow\qquad \overline{u}  = \mu\frac{\overline{k}}{f(\overline{k})}.$$ 
\item[(b)] Der stabile Verteilungsparameter erf"ullt die Beziehung
           $$\tilde{H}_u(t,\overline{k},\overline{u},1,\overline{q}) = \overline{q} f(\overline{k}) - \frac{1}{1-\overline{u}} = 0
             \qquad\Leftrightarrow\qquad \overline{u}  = 1- \frac{1}{\overline{q}f(\overline{k})}.$$
\item[(c)] F"ur den konstanten Schattenpreis erhalten wir in der adjungierten Gleichung
           $$\varrho \overline{q} - \tilde{H}_k(t,\overline{k},\overline{u},1,\overline{q})
             = \varrho \overline{q} -  \overline{q}[\overline{u} f(\overline{k})-\mu] - \frac{f'(\overline{k})}{f(\overline{k})} = 0.$$
           Die Division durch $\overline{q}$ und das Einsetzen der Beziehungen (a) und (b) liefert
           $$f'(\overline{k}) = \varrho+\mu.$$
\end{enumerate}
Wegen $f'(\overline{k}) > \mu$ ergibt sich $\overline{u} \in (0,1)$.
Damit wird durch die Beziehungen (a)--(c) ein eindeutiges Gleichgewicht $(\overline{k},\overline{u})$ festgelegt.
Abschlie"send bemerken wir,
dass die Hamilton-Funktion $\mathscr{H}$ stets konkav in $k$ ist. \hfill $\square$}
\end{beispiel}  

       \rhead[]{Aufgaben mit Zustandsbeschr"ankungen \hspace*{1mm} \thepage}
       \section{Aufgaben mit Zustandsbeschr"ankungen}
\subsection{Notwendige Optimalit"atsbedingungen}
\begin{theorem}[Pontrjaginsches Maximumprinzip] \label{SatzPAUHZA}\index{Pontrjaginsches Maximumprinzip}
Es sei $\big(x_*(\cdot),u_*(\cdot)\big) \in \mathscr{B}_{\rm adm} \cap \mathscr{B}_{\rm Lip}$.
Weiterhin nehmen wir an,
dass
$$\int_0^\infty \big\|\varphi\big(t,x_*(t),u_*(t)\big)\big\| \, dt < \infty, \qquad \int_0^\infty \big\|\varphi_x\big(t,x_*(t),u_*(t)\big)\big\| \, dt < \infty$$
ausfallen und es m"oge zu jedem $\delta>0$ ein $T>0$ existieren mit
\begin{eqnarray*}
&& \int_T^\infty \big\| \varphi\big(t,x(t),u_*(t)\big)-\varphi\big(t,x'(t),u_*(t)\big) - \varphi_x\big(t,x_*(t),u_*(t)\big)\big(x(t)-x'(t)\big) \big\| \, dt \\
&& \hspace*{20mm} \leq \delta \|x(\cdot)-x'(\cdot)\|_\infty
\end{eqnarray*}
f"ur alle $x(\cdot), x'(\cdot) \in W^1_\infty(\R_+,\R^n)$ mit $\|x(\cdot)-x_*(\cdot)\|_\infty < \gamma$, $\|x'(\cdot)-x_*(\cdot)\|_\infty < \gamma$. \\[2mm]
Ist $\big(x_*(\cdot),u_*(\cdot)\big)$ ein starkes lokales Minimum der Aufgabe (\ref{PAUH1})--(\ref{PAUH5}),
dann existieren eine Zahl $\lambda_0 \geq 0$, eine Vektorfunktion $p(\cdot):\R_+ \to \R^n$
und auf den Mengen
$$T_j=\big\{t \in \overline{\R}_+ \,\big|\, g_j\big(t,x_*(t)\big)=0\big\}, \quad j=1,...,l,$$
konzentrierte nichtnegative regul"are Borelsche Ma"se $\mu_j$ endlicher Totalvariation
(wobei s"amtliche Gr"o"sen nicht gleichzeitig verschwinden) derart, dass
\begin{enumerate}
\item[(a)] die Vektorfunktion $p(\cdot)$ von beschr"ankter Variation ist, der adjungierten Gleichung\index{adjungierte Gleichung}
           \begin{eqnarray}
           p(t)&=&- \lim_{t \to \infty}h_{1x}^T\big(t,x_*(t)\big) l_1 + \int_t^\infty H_x\big(s,x_*(s),u_*(s),p(s),\lambda_0\big) \, ds \nonumber \\
           \label{SatzPAUHZA1} & & -\sum_{j=1}^l \int_t^\infty g_{jx}\big(s,x_*(s)\big)\, d\mu_j(s)
           \end{eqnarray}
           gen"ugt und die Transversalit"atsbedingung\index{Transversalit"atsbedingungen}
           \begin{equation}\label{SatzPAUHZA2}
           p(0) = {h_0'}^T\big(x_*(0)\big)l_0
           \end{equation}
           erf"ullt;
\item[(b)] f"ur fast alle $t\in \R_+$ die Maximumbedingung
           \begin{equation}\label{SatzPAUHZA3}
           H\big(t,x_*(t),u_*(t),p(t),\lambda_0\big) = \max_{u \in U} H\big(t,x_*(t),u,p(t),\lambda_0\big)
           \end{equation}
           gilt.
\end{enumerate}
\end{theorem}

Die Adjungierte $p(\cdot)$ ist linksseitig stetig und es gilt
im rechten Endpunkt die Transversalit"atsbedingung \index{Transversalit"atsbedingungen}
$$\lim_{t \to \infty} p(t) = \lim_{t \to \infty} \bigg[-{h_1'}^T\big(x_*(t)\big)l_1-\sum_{j=1}^l g_{jx}\big(t,x_*(t)\big)\, \mu_j(\{\infty\})\bigg].$$

\begin{beispiel}[Abbau einer nicht erneuerbaren Ressource] \label{ExampleRessource}
{\rm Wir betrachten die Aufgabe
\begin{eqnarray}
&& \label{Ressource1} J\big(x(\cdot),y(\cdot),u(\cdot)\big)
   =\int_0^\infty e^{-\varrho t}\big[pf\big(u(t)\big)-ry(t)-qu(t)\big] \, dt \to \sup, \\
&& \label{Ressource2} \dot{x}(t) = -u(t),\quad \dot{y}(t)=cf\big(u(t)\big), \quad x(0)=x_0>0,\quad y(0)=y_0\geq 0, \\
&& \label{Ressource3} x(t) \geq 0, \qquad u(t) \geq 0, \qquad b,c,\varrho, q,r >0, \qquad \varrho - rc>0.
\end{eqnarray}
Die Funktion $f$ sei zweimal stetig differenzierbar, $f'>0$, $f'(0)<\infty$, $f''<0$ und es sei
$f'(u) \to 0$ f"ur $u \to \infty$.
In der vorliegenden Formulierung der Aufgabe wurde im Vergleich zu Seierstad \& Syds\ae ter \cite{Seierstad} die
Restriktion $\liminf\limits_{t \to \infty} x(t) \geq 0$ durch die Zustandsbeschr"ankung $x(t) \geq 0$ in (\ref{Ressource3})
ersetzt. \\[2mm]
"Okonomische Interpretation:
$x(t)$ bezeichnet die Menge einer nat"urlichen Ressource und $u(t)$ ist die industrielle Abbaurate dieser Ressource.
Auf Basis der Ressource werden G"uter mit der Produktionsrate $f\big(u(t)\big)$ hergestellt.
Die Kosten der Herstellung einer Produktionseinheit ist $q$ und der Preis einer G"utereinheit am Markt betr"agt $p$.
Bei der Herstellung der G"uter entstehen proportional zur Produktion Abf"alle,
deren Gesamtmenge durch $y(t)$ beschrieben wird.
Die Kosten der Beseitigung der negativen Auswirkungen der Abfallmenge sind $ry(t)$.
Im Weiteren gehen wir von dem Preis $p=1$ aus. \\[2mm]
Wegen der Zustandsbeschr"ankung sind f"ur jeden zul"assigen Steuerungprozess die Beschr"ankungen an die Dynamik in Theorem \ref{SatzPAUHZA} erf"ullt.
Wir stellen die Optimalit"atsbedingungen von Theorem \ref{SatzPAUHZA} mit $\lambda_0=1$ auf:
\begin{enumerate}
\item[(a)] Die Pontrjagin-Funktion der Aufgabe (\ref{Ressource1})--(\ref{Ressource3}) lautet
           $$H(t,x,y,u,p_1,p_2,1) = p_1 (-u)+p_2 cf(u) + e^{-\varrho t}[f(u)-ry-qu].$$
\item[(b)] Die Adjungierten gen"ugen den Gleichungen
           $$p_1(t)=\int_t^\infty \, d\mu(s), \qquad \dot{p}_2(t)=r e^{-\varrho t} \Rightarrow p_2(t)=-\frac{r}{\varrho}e^{-\varrho t} + K.$$
           Das auf der Menge $T=\{t \in \overline{\R}_+ \,|\, x_*(t)=0\}$ konzentrierte regul"are Ma"s $\mu$ ist nichtnegativ.
           Daher ist $p_1(t) \geq 0$ "uber $\R_+$ und monoton fallend.
           Ferner erhalten wir $K=0$ aus den Transversalit"atsbedingungen bez"uglich dem Zustand $y$.
\item[(c)] Die Maximumbedingung k"onnen wir auf folgende Aufgabe reduzieren:
           $$\max_{u \geq 0} \Big[ -p_1(t) u +c p_2(t) f(u) + e^{-\varrho t}[f(u)-qu]\Big].$$
           Das Einsetzen der Darstellung f"ur $p_2(t)$ liefert weiterhin mit $d=(\varrho -rc)/\varrho$:
           $$\max_{u \geq 0} \Big( d f(u)e^{-\varrho t}-u\big(p_1(t)+ qe^{-\varrho t}\big)\Big).$$
\end{enumerate}
Die Reduktion der Maximumbedingung f"uhrt f"ur festes $t$ zu der Funktion
$$g(u)= d f(u)e^{-\varrho t}-u\big(p_1(t)+ qe^{-\varrho t}\big).$$
Diese Funktion ist zweimal stetig differenzierbar und es gilt
$$g'(u)= \big(d f'(u)-q\big)e^{-\varrho t}-p_1(t), \quad g''(u)=df''(u)e^{-\varrho t}, \quad d=\frac{\varrho -rc}{\varrho}>0.$$
Daher ist $g$ streng konkav und besitzt auf der Menge $U=\{u\geq 0\}$ ein Maximum,
da $f'(u)>0$ und $f'(u) \to 0$ f"ur $u \to \infty$ gelten.
Wir diskutieren drei F"alle:
\begin{enumerate}
\item[(A)] $df'(0)\leq q$: In diesem Fall ist $g'(0) \leq 0$ und man erh"alt
           $$u_*(t) \equiv 0, \quad x_*(t) \equiv x_0, \quad
             y_*(t)=y_0 + cf(0)t, \quad p_1(t) \equiv 0, \quad p_2(t)=-\frac{r}{\varrho}e^{-\varrho t}.$$
           Da die Zustandsbeschr"ankung nichtaktiv ist, gelten die Voraussetzungen und Optimalit"atsbedingungen in Theorem \ref{SatzPAUHZA}.
\item[(B)] $df'(0)> q$ und $p_1(0)=0$:
           Aus $g'(u)=\big(d f'(u)-q\big)e^{-\varrho t}=0$ erhalten wir die optimale Strategie $u_*(t)=u_0>0$ f"ur alle $t \in \R_+$.
           Also gilt $x_*(t)=x_0-u_0t$ auf $\R_+$, was der Zustandsbeschr"ankung widerspricht.
\item[(C)] $df'(0)> q$ und $p_1(0)>0$:
           Wegen $p_1(0)>0$ wird die Ressource vollst"andig abgebaut.
           Andernfalls w"are $p_1(t)=p_1(0)>0$ "uber $\R_+$, was (\ref{SatzPAUHZA2}) widerspricht.
           Da die Ressource vollst"andig abgebaut wird,
           gibt es ein $t'>0$ mit $x_*(t)>0$ f"ur $t \in [0,t')$ und $x_*(t)=0$ f"ur $t\geq t'$.
           Demnach folgt unmittelbar $u_*(t)=0$ f"ur $t\geq t'$. \\
           F"ur $t\geq t'$ ist $p_1(\cdot)$ monoton fallend.
           Ferner erhalten wir f"ur $t \in \R_+$ die Beziehung
           $$g'(u)=0 \qquad\Rightarrow\qquad f'\big(u(t)\big)=\frac{1}{d}(q+p_1(t)e^{\varrho t}).$$
           W"urde demnach die Adjungierte $p_1(\cdot)$ f"ur $t \geq t'$ eine Unstetigkeitstelle besitzen,
           dann folgt aus der Monotonie von $p_1(\cdot)$, dass die Abbaurate sich wieder sprunghaft vergr"o"sert.
           Diese Steuerung f"uhrt zu einem erneuten Abbau der Ressource, obwohl diese bereits vollst"andig aufgebraucht ist.
           Daher ist die Adjungierte stetig. \\
           F"ur die Adjungierte erhalten wir damit
           $$p_1(t) = \big(df'(0)-q\big)e^{-\varrho t'} \mbox{ f"ur } t \leq t', \qquad
             p_1(t) = \big(df'(0)-q\big)e^{-\varrho t} \mbox{ f"ur } t \geq t'.$$
           Wir zeigen noch, dass der Zeitpunkt $t'$ existiert und eindeutig ist:
           Durch
           $$f'\big(u_\tau(t)\big)=\frac{1}{d}(q+p_1(0)e^{\varrho t})=\frac{1}{d}(q+[df'(0)-q]e^{\varrho(t-\tau)}), \quad t \in [0,\tau],$$
           und $u_\tau(t)=0$ f"ur $t \geq \tau$ wird wegen $f'\big(u_\tau(\tau)\big)=f'(0)$ eine Familie $u_\tau(\cdot)$ stetiger Funktionen definiert.
           Dabei gilt $f'\big(u_\tau(t)\big) < f'\big(u_s(t)\big)$, d.\,h. $u_\tau(t) > u_s(t)$, f"ur alle $t \in [0,\tau]$ und $\tau>s$.
           Damit ist die Familie
           $$U(\tau):= \int_0^\infty u_\tau(t) \, dt$$
           streng monoton wachsend und es gelten $U(0)=0$, $U(\tau) \to \infty$ f"ur $\tau \to \infty$.
           Der Parameter $t'$ ergibt sich dann aus der Bedingung $U(t')=x_0$. \hfill $\square$
\end{enumerate}}
\end{beispiel}
       \subsection{Der Nachweis der notwendigen Optimalit"atsbedingungen}
Seien $x_*(\cdot) \in \mathscr{B}_{\rm Lip} \cap C_{\lim}(\R_+,\R^n)$ und $V_\gamma= \{ (t,x) \in \overline{\R}_+ \times \R^n \,|\, \|x-x(t)\| \leq \gamma\}$. \\
Wir betrachten f"ur $\big(x(\cdot),u(\cdot)\big) \in C_{\lim}(\R_+,\R^n) \times L_\infty(\R_+,\R^m)$ die Abbildungen
\begin{eqnarray*}
J\big(x(\cdot),u(\cdot)\big) &=& \int_0^\infty \omega(t)f\big(t,x(t),u(t)\big) \, dt, \\
F\big(x(\cdot),u(\cdot)\big)(t) &=& x(t) -x(t_0) -\int_0^t \varphi\big(s,x(s),u(s)\big) \, ds, \quad t \in \R_+,\\
H_0\big(x(\cdot)\big) &=& h_0\big(x(0)\big), \qquad H_1\big(x(\cdot)\big) = \lim_{t \to \infty} h_1\big(t,x(t)\big), \\
G_j\big(x(\cdot)\big) &=& \max_{t \in \overline{\R}_+} g_j\big(t,x(t)\big), \quad j=1,...,l.
\end{eqnarray*}
Die Abbildungen fassen wir in folgenden Funktionenr"aumen auf:
\begin{eqnarray*}
J &:& C_{\lim}(\R_+,\R^n) \times L_\infty(\R_+,\R^m) \to \R, \\
F &:& C_{\lim}(\R_+,\R^n) \times L_\infty(\R_+,\R^m) \to C_{\lim}(\R_+,\R^n), \\
H_i &:& C_{\lim}(\R_+,\R^n) \to \R^{s_i}, \quad i=0,1, \\
G_j &:& C_{\lim}(\R_+,\R^n) \to \R, \quad j=1,...,l.
\end{eqnarray*}

Wir setzen $\mathscr{F}=(F,H_0,H_1)$, sowie die Menge $\mathscr{U}$ gem"a"s
\begin{eqnarray*}
\mathscr{U}=\big\{ u(\cdot)  \in L_\infty(\R_+,U) &\big|& u(t)=u_*(t) + \chi_M(t)\big(w(t)-u_*(t)\big), \; w(\cdot) \in L_\infty(\R_+,U), \\
                                                  &     & M \subset \R_+ \mbox{ me"sbar und beschr"ankt} \big\}
\end{eqnarray*}
und pr"ufen f"ur die Extremalaufgabe
\begin{equation} \label{ExtremalaufgabePMPUHZA}
J\big(x(\cdot),u(\cdot)\big) \to \inf, \quad \mathscr{F}\big(x(\cdot),u(\cdot)\big)=0, \quad G_j\big(x(\cdot)\big) \leq 0, \quad u(\cdot) \in \mathscr{U}
\end{equation}
im Punkt $\big(x_*(\cdot),u_*(\cdot)\big)$ die Voraussetzungen von Theorem \ref{SatzExtremalprinzipStark}:

\begin{enumerate}
\item[(A$_2$)] Mit Verweis auf Abschnitt \ref{AbschnittBeweisPMPUH} sind nur noch die Abbildungen $G_j$ zu diskutieren.
               In den Beispielen \ref{SubdifferentialMaximum3} und \ref{SubdifferentialMaximum4} wird gezeigt,
               dass die Funktionen $G_j$ Hintereinanderausf"uhrungen einer stetigen, konvexen, eigentlichen Funktion und
               einer Fr\'echet-differenzierbaren Abbildung sind.
               Daher sind nach Lemma \ref{LemmaRichtungsableitung} die Funktionen $G_j$ in $x_*(\cdot)$ lokalkonvex und bez"uglich jeder Richtung
               gleichm"a"sig differenzierbar.              
\end{enumerate}

Zur Extremalaufgabe (\ref{ExtremalaufgabePMPUHZA}) definieren wir auf
$$C_{\lim}(\R_+,\R^n) \times L_\infty(\R_+,\R^m) \times \R \times C_{\lim}^*(\R_+,\R^n) \times \R^n \times \R^{s_0} \times  \R^{s_1}$$
die Lagrange-Funktion $\mathscr{L}=\mathscr{L}\big(x(\cdot),u(\cdot),\lambda_0,y^*,l_0,l_1,\lambda\big)$,
$$\mathscr{L}= \lambda_0 J\big(x(\cdot),u(\cdot)\big)+ \big\langle y^*, F\big(x(\cdot),u(\cdot)\big) \big\rangle
                         +l_0^T H_0\big(x(\cdot)\big)+l_1^T H_1\big(x(\cdot)\big) + \sum_{j=1}^l \lambda_j G_j\big(x(\cdot)\big).$$
Ist $\big(x_*(\cdot),u_*(\cdot)\big)$ eine starke lokale Minimalstelle der Extremalaufgabe (\ref{ExtremalaufgabePMPUHZA}),
dann existieren nach Theorem \ref{SatzExtremalprinzipStark}
nicht gleichzeitig verschwindende Lagrangesche Multiplikatoren
$$\lambda_0 \geq 0, \qquad y^* \in C_{\lim}^*(\R_+,\R^n), \qquad l_0 \in  \R^{s_0}, \qquad l_1 \in  \R^{s_1}$$
und $\lambda_1 \geq 0,...,\lambda_l \geq 0$ derart,
dass gelten:
\begin{enumerate}
\item[(a)] Die Lagrange-Funktion besitzt bez"uglich $x(\cdot)$ in $x_*(\cdot)$ einen station"aren Punkt, d.\,h.
           \begin{equation}\label{SatzPMPUHZALMR1}
           0 \in \partial_x \mathscr{L}\big(x_*(\cdot),u_*(\cdot),\lambda_0,y^*,l_0,l_1,\lambda\big);
           \end{equation}         
\item[(b)] Die Lagrange-Funktion erf"ullt bez"uglich $u(\cdot)$ in $u_*(\cdot)$ die Minimumbedingung
           \begin{equation}\label{SatzPMPUHZALMR2}
           \hspace*{-3mm} \mathscr{L}\big(x_*(\cdot),u_*(\cdot),\lambda_0,y^*,l_0,l_1,\lambda\big)
           = \min_{u(\cdot) \in \mathscr{U}} \mathscr{L}\big(x_*(\cdot),u(\cdot),\lambda_0,y^*,l_0,l_1,\lambda\big);
           \end{equation}
\item[(c)] Die komplement"aren Schlupfbedingungen gelten, d.\,h.
           \begin{equation}\label{SatzPMPUHZALMR3}
           0 = \lambda_j G_j\big(x(\cdot)\big), \qquad i=1,...,l.
           \end{equation}
\end{enumerate}
Aufgrund (\ref{SatzPMPUHZALMR1}) ist folgende Variationsgleichung f"ur alle $x(\cdot) \in C_{\lim}(\R_+,\R^n)$ erf"ullt: 
\begin{eqnarray*}
0 &=& \lambda_0 \cdot \int_0^\infty \omega(t) \big\langle f_x\big(t,x_*(t),u_*(t)\big),x(t) \big\rangle\, dt \\
  & & + \int_0^\infty \bigg[ x(t)-x(0) - \int_0^t \varphi_x\big(s,x_*(s),u_*(s)\big) x(s) \,ds \bigg]^T  d\mu(t) \\
  & & + \big\langle l_0, h_0'\big(x_*(0)\big) x(0) \big\rangle + \big\langle l_1, \lim_{t \to \infty} h_{1x}\big(t,x_*(t)\big) x(t) \big\rangle \\
  & & + \sum_{j=1}^l \lambda_j \int_0^\infty \big\langle g_{jx}\big(t,x_*(t)\big),x(t) \big\rangle \,d\tilde{\mu}_j(t).
\end{eqnarray*}
Darin k"onnen wir ohne Einschr"ankung annehmen (vgl. \cite{Ioffe}),
dass alle Ma"se $\mu_j=\lambda_j \tilde{\mu}_j$ auf den Mengen
$$T_j=\big\{t \in \overline{\R}_+ \,\big|\, g_j\big(t,x_*(t)\big)=0\big\}, \quad j=1,...,l,$$
konzentriert sind.
Wir "andern die Integrationsreihenfolge und setzen $p(t)= \displaystyle \int_t^\infty \, d\mu(s)$.
Dann ergeben sich aus der eindeutigen Darstellung eines stetigen linearen Funktionals im Raum $C_{\lim}(\R_+,\R^n)$
die Bedingungen (\ref{SatzPAUHZA1}) und (\ref{SatzPAUHZA2}).
Gem"a"s (\ref{SatzPMPUHZALMR2}) gilt f"ur alle $u(\cdot) \in \mathscr{U}$ die Ungleichung
$$\int_0^\infty H\big(t,x_*(t),u_*(t),p(t),\lambda_0\big) \, dt \geq
  \int_0^\infty H\big(t,x_*(t),u(t),p(t),\lambda_0\big) \, dt.$$
Daraus folgt abschlie"send durch Standardtechniken f"ur Lebesguesche Punkte die Maximumbedingung (\ref{SatzPAUHZA3}).
       \subsection{Hinreichende Bedingungen nach Arrow} \label{AbschnittArrowPMPZB}
Die folgende Herleitung\index{hinreichende Bedingungen!Arrow@-- nach Arrow} basiert wieder auf der Darstellung in Seierstad \& Syds\ae ter \cite{Seierstad}.
Wir betrachten das Steuerungsproblem
\begin{eqnarray}
&& \label{HBPAUHZB1} J\big(x(\cdot),u(\cdot)\big) = \int_0^\infty \omega(t)f\big(t,x(t),u(t)\big) \, dt \to \inf, \\
&& \label{HBPAUHZB2} \dot{x}(t) = \varphi\big(t,x(t),u(t)\big), \\
&& \label{HBPAUHZB3} x(0)=x_0, \qquad \lim_{t \to \infty} x(t)=x_1, \\
&& \label{HBPAUHZB4} u(t) \in U \subseteq \R^m, \quad U \not= \emptyset, \\
&& \label{HBPAUHZB5} g_j\big(t,x(t)\big) \leq 0 \quad \mbox{f"ur alle } t \in \R_+, \quad j=1,...,l.
\end{eqnarray}
In der Aufgabenstellung schlie"sen wir den Fall nicht aus,
dass durch die Randbedingungen (\ref{HBPAUHZB3}) gewisse Komponenten der Punkte $x_0$ und $x_1$ nicht fest vorgegeben,
sondern ohne Einschr"ankung sind. \\[2mm]
Wir definieren die Mengen $V_\gamma$ und $V_\gamma(t)$:
$$V_\gamma=\{ (t,x) \in \overline{\R}_+ \times \R^n \,|\, \|x-x(t)\| \leq \gamma\}, \quad
  V_\gamma(t)=\{ x \in \R^n \,|\, \|x-x_*(t)\| \leq \gamma\}.$$
Au"serdem bezeichnet $\mathscr{H}(t,x,p) = \sup\limits_{u \in U} H(t,x,u,p,1)$ die Hamilton-Funktion $\mathscr{H}$ im normalen Fall.

\begin{theorem} \label{SatzHBPMPUHZB}
In der Aufgabe (\ref{HBPAUHZB1})--(\ref{HBPAUHZB5}) sei $\big(x_*(\cdot),u_*(\cdot)\big) \in \mathscr{B}_{\rm Lip} \cap \mathscr{B}_{\rm adm}$.
Au"serdem sei die Vektorfunktion $p(\cdot):\R_+ \to \R^n$ st"uckweise stetig,
besitze h"ochstens abz"ahlbar viele Sprungstellen $s_k \in (0,\infty)$,
die sich nirgends im Endlichen h"aufen,
und $p(\cdot)$ sei zwischen diesen Spr"ungen stetig differenzierbar. 
Ferne gelte:
\begin{enumerate}
\item[(a)] Das Tripel $\big(x_*(\cdot),u_*(\cdot),p(\cdot)\big)$
           erf"ullt (\ref{SatzPAUH1})--(\ref{SatzPAUH3}) in Theorem \ref{SatzPAUH} mit $\lambda_0=1$.        
\item[(b)] F"ur jedes $t \in \R_+$ ist die Funktion $\mathscr{H}\big(t,x,p(t)\big)$ konkav 
           und es sind die Funktionen $g_j(t,x)$, $j=1,...,l$, konvex bez"uglich $x$ auf $V_\gamma(t)$.
\end{enumerate}
Dann ist $\big(x_*(\cdot),u_*(\cdot)\big)$ ein starkes lokales Minimum der Aufgabe (\ref{HBPAUHZB1})--(\ref{HBPAUHZB5}).
\end{theorem}

\begin{bemerkung}{\rm
Der Teil (a) in Theorem \ref{SatzHBPMPUHZB} bedarf einer detaillierteren Diskussion.
Da wir von einer st"uckweise stetigen und zwischen den Sprungstellen stetig differenzierbaren Adjungierten $p(\cdot)$ ausgehen,
k"onnen wir die adjungierte Gleichung in Integraldarstellung in die Form einer st"uckweise definierten Differentialgleichung mit
Sprungbedingungen "uberf"uhren.
Es bezeichnen $0<s_1<...<s_d<T$ die Unstetigkeitsstellen der Adjungierten $p(\cdot)$ im Intervall $(0,T)$.
Dann gelten die Sprungbedingungen
$$p(s_k^-)-p(s_k^+)= -\sum_{j=1}^l \beta_j^k g_{jx}\big(s_k,x_*(s_k)\big),\qquad \beta_j^k = \mu_j(\{s_k\}) \geq 0, \quad k=1,...,d.$$
Ferner gibt es eine st"uckweise stetige Vektorfunktion $\lambda(\cdot):\R_+ \to \R^l$ derart,
dass die Differentialgleichung
$$\dot{p}(t)=- H_x\big(t,x_*(t),u_*(t),p(t),1\big) + \sum_{j=1}^l \lambda_j(t) g_{jx}\big(t,x_*(t)\big)$$
st"uckweise auf $(s_k,s_{k+1})$, $k=0,...,d$, gilt. Dabei haben wir $s_0=0$, $s_{d+1}=T$ gesetzt. \\
Abschließend halten wir fest,
dass wegen der Positivit"at der Ma"se $\mu_j$ und der Konzentration dieser Ma"se auf den Mengen
$$T_j=\big\{t \in \overline{\R}_+ \,\big|\, g_j\big(t,x_*(t)\big)=0\big\}, \quad j=1,...,l,$$
neben $\beta_j g_j\big(T,x_*(T)\big) =0$ die Bedingungen
$$\lambda_j(t) \geq 0, \qquad \lambda_j(t)g_j\big(t,x_*(t)\big)=0$$
auf $[0,T]$ und f"ur $k=1,...,d$ in den Sprungstellen 
$$\beta_j^k \geq 0, \qquad \beta_j^k g_j\big(s_k,x_*(s_k)\big)=0$$
f"ur $j=1,...,l$ und gelten. \hfill $\square$}
\end{bemerkung}

{\bf Beweis} Wie im Abschnitt \ref{AbschnittArrowPMP} ergibt sich die Ungleichung
\begin{eqnarray*}
    \Delta(T)
&=& \int_0^T \omega(t)\big[f\big(t,x(t),u(t)\big)-f\big(t,x_*(t),u_*(t)\big)\big] \, dt \\
&\geq & \int_0^T\big[\mathscr{H}\big(t,x_*(t),p(t)\big)-\mathscr{H}\big(t,x(t),p(t)\big)\big] \, dt 
        + \int_0^T \langle p(t), \dot{x}(t)-\dot{x}_*(t) \rangle \, dt.
\end{eqnarray*}
F"ur $T \not=s_k$ ergibt sich ferner nach der eingangs getroffenen Bemerkung:
\begin{eqnarray*}
    \Delta(T)
&\geq& \int_0^T \langle \dot{p}(t),x(t)-x_*(t)\rangle + \langle p(t), \dot{x}(t)-\dot{x}_*(t) \rangle \, dt \\
&&    \hspace*{10mm} - \int_0^T \sum_{j=1}^l \lambda_j(t) \big\langle g_{jx}\big(t,x_*(t)\big) , x(t)-x_*(t) \big\rangle \, dt \\
&&    \hspace*{10mm} + \sum_{s_k <T} \langle p(s_k^-)-p(s_k^+),x(s_k)-x_*(s_k)\rangle \\
&\geq& \int_0^T \langle \dot{p}(t),x(t)-x_*(t)\rangle + \langle p(t), \dot{x}(t)-\dot{x}_*(t) \rangle \, dt \\
&=& \langle p(T),x(T)-x_*(T)\rangle-\langle p(0),x(0)-x_*(0)\rangle.
\end{eqnarray*}
Da $p(\cdot)$ eine Funktion beschr"ankter Variation ist, gilt au"serdem
$$\lim_{k \to \infty} \|p(s_k^-)-p(s_k^+)\|=0.$$
Damt ist der Grenz"ubergang $T \to \infty$ gerechtfertigt.
Beachten wir die Transversalit"atsbedingung
$$\lim_{t \to \infty} p(t) = \lim_{t \to \infty} \bigg[l_1-\sum_{j=1}^l g_{jx}\big(t,x_*(t)\big)\, \mu_j(\{\infty\})\bigg],$$
und die Konvexit"at der Abbildungen $x \to g_j(t,x)$ in $t=\infty$,
so erhalten wir die Ungleichung
$$\lim_{T \to \infty} \Delta(T) = \lim_{T \to \infty} \int_0^T \omega(t)\big[f\big(t,x(t),u(t)\big)-f\big(t,x_*(t),u_*(t)\big)\big] \, dt \geq 0$$
f"ur alle zul"assigen $\big(x(\cdot),u(\cdot)\big)$ mit $\|x(\cdot)-x_*(\cdot)\|_\infty \leq \gamma$. \hfill $\blacksquare$

\begin{beispiel}
{\rm Im Beispiel \ref{ExampleRessource} enth"alt die Pontrjagin-Funktion
$$H(t,x,y,u,p_1,p_2,1) = p_1 (-u)+p_2 cf(u) + e^{-\varrho t}[f(u)-ry-qu]$$
bez"uglich der Zustandsvariablen $x$ und $y$ nur den Term $e^{-\varrho t}ry$.
Damit ist die Hamilton-Funktion offenbar konkav bez"uglich $(x,y)$.
Au"serdem ist die Zustandsbeschr"ankung linear in $x$. \\[2mm]
Im Beispiel \ref{ExampleRessource} diskutierten wir die folgenden F"alle:
\begin{enumerate}
\item[(A)] Im Fall $df'(0)\leq q$ gen"ugen
           $$u_*(t) \equiv 0, \quad x_*(t) \equiv x_0, \quad
             y_*(t)=y_0 + cf(0)t, \quad p_1(t) \equiv 0, \quad p_2(t)=-\frac{r}{\varrho}e^{-\varrho t}$$
           den Bedingungen (\ref{SatzPAUH1})--(\ref{SatzPAUH3}) in Theorem \ref{SatzPAUH} mit $\lambda_0=1$ und
           liefert somit ein starkes lokales Maximum.
\item[(B)] Im Fall $df'(0)> q$ und $p_1(0)=0$ erhielten wir keinen Kandidaten.
\item[(C)] Im Fall $df'(0)> q$ und $p_1(0)>0$ wird die Ressource vollst"andig abgebaut.
           Daher gibt es einen eindeutig bestimmten Zeitpunkt $t'>0$ mit $x_*(t)>0$ f"ur $t \in [0,t')$, $x_*(t)=0$ f"ur $t\geq t'$ und
           $u_*(t)=0$ f"ur $t\geq t'$.
           Die Bedingungen (\ref{SatzPAUH1})--(\ref{SatzPAUH3}) in Theorem \ref{SatzPAUH} lieferten mit $\lambda_0=1$ die Adjungierte
           $$p_1(t) = \big(df'(0)-q\big)e^{-\varrho t'} \mbox{ f"ur } t \leq t', \qquad p_1(t) = \big(df'(0)-q\big)e^{-\varrho t} \mbox{ f"ur } t \geq t'$$
           und wir erhielten bez"uglich $u_*(t)$ im Beispiel \ref{ExampleRessource} f"ur $t \leq t'$ den Zusammenhang
           $$f'\big(u_*(t)\big)=\frac{1}{d}(q+p_1(t)e^{\varrho t}),$$
           sowie die Stetigkeit von $u_*(\cdot)$ in $t=t'$.
           Der Vollst"andigkeit halber geben wir noch $y_*(t)$ an:
           $$y_*(t) = y_0 + \int_0^t c f\big(u_*(s)\big) \, ds \mbox{ f"ur } t \leq t', \qquad y_*(t) = y_*(t') \mbox{ f"ur } t \geq t'.$$
           Nach Theorem \ref{SatzHBPMPUHZB} stellt der Kandidat ein starkes lokale Maximum dar. \hfill $\square$
\end{enumerate}}
\end{beispiel}   

       \rhead[]{Aufgaben mit verschiedenen Horizonten \hspace*{1mm} \thepage}
       \section{Aufgaben mit endlichem und unendlichem Horizont}
\subsection{Zur Approximation durch eine Aufgabe "uber endlichem Horizont}
Zur Behandlung von Aufgaben mit unendlichem Zeithorizont findet man in der Literatur verschiedene Methoden und
eine entsprechend gro"se Anzahl an Optimalit"atsbegriffen. \\
In der Aufgabe (\ref{PAUH1})--(\ref{PAUH5}) bezeichnet $\mathscr{A}_{\rm adm}$ die Menge aller zul"assigen
Steuerungsprozesse,
d.\,h. die Menge aller $\big(x(\cdot),u(\cdot)\big)$,
die der Dynamik (\ref{PAUH2}) zur Anfangs- und Randbedingung (\ref{PAUH3}),
den Steuerrestriktionen (\ref{PAUH4}) und den Zustandsbeschr"ankungen (\ref{PAUH5}) gen"ugen,
und f"ur die das Zielfunktional (\ref{PAUH1}) endlich ist.
Zur Gegen"uberstellung von verschiedenen Optimalit"atsbegriffen betrachten wir globale Optimalit"atskriterien.

\begin{enumerate}
\item[(a)] Der Steuerungsprozess $\big(x_*(\cdot),u_*(\cdot)\big) \in \mathscr{A}_{\rm adm}$ hei"st
           global optimal (GO)\index{optimal, catching up!global@--, global}, falls
           $$J\big(x_*(\cdot),u_*(\cdot)\big) \leq J\big(x(\cdot),u(\cdot)\big)$$
           f"ur alle $\big(x(\cdot),u(\cdot)\big) \in \mathscr{A}_{\rm adm}$ gilt.
\item[(b)] Der Steuerungsprozess $\big(x_*(\cdot),u_*(\cdot)\big) \in \mathscr{A}_{\rm adm}$ hei"st
           streng global optimal (SGO)\index{optimal, catching up!streng global@--, streng global}, falls
           $$J\big(x_*(\cdot),u_*(\cdot)\big) < J\big(x(\cdot),u(\cdot)\big)$$
           f"ur alle $\big(x(\cdot),u(\cdot)\big) \in \mathscr{A}_{\rm adm}$ mit
           $\big(x(\cdot),u(\cdot)\big) \not= \big(x_*(\cdot),u_*(\cdot)\big)$ gilt.
\end{enumerate}
Zur Untersuchung des Steuerungsproblems mit unendlichem Zeithorizont wird in den meisten Ans"atzen der Aufgabe (\ref{PAUH1})--(\ref{PAUH5})
ein Problem "uber einem endlichen Zeitintervall $[0,T]$ zugeordnet und der Grenz"ubergang $T \to \infty$ betrachtet.
Wir nennen dies die Approximation durch endlichen Horizont. \\[2mm]
Das zur Aufgabe (\ref{PAUH1})--(\ref{PAUH5}) geh"orende Problem mit endlichem Zeithorizont lautet
\begin{eqnarray}
&& \label{PAUH1a} J_T\big(x(\cdot),u(\cdot)\big) = \int_0^T \omega(t) f\big(t,x(t),u(t)\big)\, dt \to \inf, \\
&& \label{PAUH2a} \dot{x}(t) = \varphi\big(t,x(t),u(t)\big), \\
&& \label{PAUH3a} h_0\big(x(0)\big)=0, \qquad h_1\big(T,x(T)\big)=0, \\
&& \label{PAUH4a} u(t) \in U \subseteq \R^m, \quad U \not= \emptyset, \\
&& \label{PAUH5a} g_j\big(t,x(t)\big) \leq 0 \quad \mbox{f"ur alle } t \in [0,T], \quad j=1,...,l.
\end{eqnarray}
Wegen des "Uberganges zu einer Aufgabe "uber endlichem Zeithorizont und der anschlie"senden Betrachtung des Grenzwertes $T \to \infty$
werden die Optimalit"atsbegriffe angepasst. \\[2mm]
Zur Aufgabe (\ref{PAUH1a})--(\ref{PAUH5a}) definieren wir den Defekt $\Delta(T)$:
$$\Delta(T)= J_T\big(x(\cdot),u(\cdot)\big) - J_T\big(x_*(\cdot),u_*(\cdot)\big).$$
\begin{enumerate}
\item[(c)] Der Steuerungsprozess $\big(x_*(\cdot),u_*(\cdot)\big) \in \mathscr{A}_{\rm adm}$ hei"st
           overtaking optimal (OT)\index{optimal, catching up!overtaking@--, overtaking},
           falls es zu jedem $\big(x(\cdot),u(\cdot)\big) \in \mathscr{A}_{\rm adm}$ eine Zahl $T_0$ gibt mit
           $$\Delta(T)= J_T\big(x(\cdot),u(\cdot)\big) - J_T\big(x_*(\cdot),u_*(\cdot)\big)\geq 0$$
           f"ur alle $T \geq T_0$ (von Weizs"acker \cite{VonWeiz}).
\item[(d)] Der Steuerungsprozess $\big(x_*(\cdot),u_*(\cdot)\big) \in \mathscr{A}_{\rm adm}$ hei"st
           catching up optimal (CU)\index{optimal, catching up}, falls
           $$\liminf_{T \to \infty} \Delta(T)
             = \liminf_{T \to \infty} J_T\big(x(\cdot),u(\cdot)\big) - J_T\big(x_*(\cdot),u_*(\cdot)\big) \geq 0$$
           f"ur jedes $\big(x(\cdot),u(\cdot)\big) \in \mathscr{A}_{\rm adm}$ gilt (Gale \cite{Gale}).
\item[(e)] Das Paar $\big(x_*(\cdot),u_*(\cdot)\big) \in \mathscr{A}_{\rm adm}$ hei"st
           sporadically catching up optimal (SCU)\index{optimal, catching up!sporadically@--, sporadically catching up}, falls
           $$\limsup_{T \to \infty} \Delta(T)
              = \limsup_{T \to \infty} J_T\big(x(\cdot),u(\cdot)\big) - J_T\big(x_*(\cdot),u_*(\cdot)\big)\geq 0$$
           f"ur jedes $\big(x(\cdot),u(\cdot)\big) \in \mathscr{A}_{\rm adm}$ gilt (Halkin \cite{Halkin}).
\end{enumerate}
Durch die Betrachtung auf der Menge $\mathscr{A}_{\rm adm}$ haben wir in den Definitionen bereits eingef"ugt,
dass f"ur einen zul"assigen Steuerungsprozess das Zielfunktional endlich ist.
In diversen Arbeiten zu Steuerungsproblemen mit unendlichem Zeithorizont,
z.\,B. bei Aseev \& Veliov \cite{AseVel,AseVel2,AseVel3}, darf das Integral in (\ref{PAUH1}) divergieren.
Konvergiert das Zielfunktional f"ur jeden zul"assigen Steuerungsprozess,
so sind die Optimalit"atsbegriffe (CU), (SCU) und (GO) "aquivalent, denn es gilt
$$\liminf_{T \to \infty} \Delta(T) \geq 0 \quad\Leftrightarrow\quad
  \limsup_{T \to \infty} \Delta(T) \geq 0 \quad\Leftrightarrow\quad \lim_{T \to \infty} \Delta(T)\geq 0.$$
Zur Einordnung der Begriffe (SGO) und (OT) f"uhren wir folgendes Beispiel an:

\begin{beispiel}
{\rm Wir betrachten die Aufgabe
$$\int_0^\infty -u(t)x(t) \, dt \to \inf, \qquad \dot{x}(t) = -u(t)x(t), \quad x(0)=1,
  \quad u(t) \in [0,1], \quad \varrho >0.$$
Offenbar ist jede zul"assige Trajektorie $x(\cdot)$ monoton fallend und wegen $1 \geq x(t) > 0$ beschr"ankt.
Also besitzt sie deswegen einen Grenzwert f"ur $t \to \infty$.
Damit gilt
$$J\big(x(\cdot),u(\cdot)\big)=\int_0^\infty -u(t)x(t) \, dt = \lim_{t \to \infty} x(t)-1 \geq -1.$$
Daher ist jedes Paar $\big(x(\cdot),u(\cdot)\big) \in \mathscr{A}_{\rm adm}$ mit $x(t) \to 0$ global optimal.
Dabei gilt:
$$\lim_{t \to \infty} x(t) = 0 \quad \Leftrightarrow \quad \int_0^\infty u(t) \, dt = \infty.$$
In diesem Beispiel ist jeder Steuerungsprozess $\big(x(\cdot),u(\cdot)\big) \in \mathscr{A}_{\rm adm}$
mit $x(t) \to 0$ f"ur $t \to \infty$
optimal im Sinn von (GO), (CU) und (SCU).
Wegen $J_T\big(x(\cdot),u(\cdot)\big)\geq x(T)-1$ ist unter diesen Steuerungsprozessen nur derjenige optimal im Sinn von (OT),
f"ur den $u(t) \equiv 1$ auf $\R_+$ gilt.
Eine (SGO)-L"osung existiert nicht. \hfill $\square$}
\end{beispiel}
Dieses Beispiel verdeutlicht,
dass unter der Annahme eines endlichen Zielfunktionals zu den bereits erw"ahnten "Aquivalenzen folgende weitere Relationen
zwischen den verschiedenen Optimalit"atsbegriffen bestehen:
$$\mbox{(SGO) } \Longrightarrow \mbox{ (OT) } \Longrightarrow \mbox{ (CU) } \Longleftrightarrow \mbox{ (SCU) }
  \Longleftrightarrow \mbox{ (GO)}.$$
Ist demnach die globale Optimalstelle eindeutig, so fallen bei endlichem Zielfunktional
die aufgef"uhrten Optimalit"atsbegriffe zusammen. \\[2mm]
Die Betrachtungen dieses Abschnitts ordnen die Begriffe der globalen Optimalit"at in den Kontext 
g"angiger Optimalit"atsbegriffe ein,
die bei der Approximation mit einem endlichen Horizont angewendet werden.
Unbehandelt bleibt dabei allerdings die wesentliche Frage,
wann eine Familie $\big\{\big(x_*^T(\cdot),u_*^T(\cdot)\big)\big\}_{T \in \R_+}$ von optimalen Steuerungsprozessen,
die sich in der Approximation ergibt,
gegen ein globales Optimum konvergiert.
Dazu muss man sich vorab im Klaren sein,
dass f"ur die Approximation die fundamentale Beziehung
$$\lim_{T \to \infty} \int_0^T f(t) \, dt = \int_0^\infty f(t) \, dt$$
im Allgemeinen nur bei vorliegender Existenz des Lebesgue-Integrals im Zielfunktional gilt.
Ein bekanntes Gegenbeispiel ist der Integralsinus:
$$\lim_{T \to \infty} \int_0^T \frac{\sin t}{t} \, dt = \frac{\pi}{2}, \qquad
  \int_0^\infty \frac{\sin t}{t} \, dt \mbox{ existiert nicht.}$$
Noch weniger "uberschaubar wird die Situation in der Aufgabe (\ref{PAUH1})--(\ref{PAUH5}),
in der die Dynamik, sowie Zustands- und Steuerungsbeschr"ankungen vorkommen.
Dann k"onnen noch weitere unerw"unschte Situationen eintreten.
Dazu geben wir das folgende Beispiel an:

\begin{beispiel}
{\rm Wir betrachten die Aufgabe
\begin{eqnarray*}
&& J_T\big(x(\cdot),u(\cdot)\big)=\int_0^T e^{-\varrho t} \big(1-u(t)\big)x(t) \, dt \to \sup, \\
&& \dot{x}(t) = u(t)x(t), \quad x(0)=1, \quad u(t) \in [0,1], \quad \varrho \in (0,1).
\end{eqnarray*}
F"ur jedes feste $T>\tau$ liefert der Steuerungsprozess
$$x^T_*(t)= \left\{ \begin{array}{ll} e^t,& t \in [0,\tau), \\[1mm] e^\tau,& t \in [\tau,T], \end{array} \right. \quad
  u^T_*(t)= \left\{ \begin{array}{ll} 1,& t \in [0,\tau), \\[1mm] 0,& t \in [\tau,T], \end{array} \right. \quad
  \tau = T + \frac{\ln (1-\varrho)}{\varrho}$$
das globale Maximum.
Betrachten wir den Grenz"ubergang $T \to \infty$,
dann konvergiert die Familie $\big\{\big(x_*^T(\cdot),u_*^T(\cdot)\big)\big\}_{T \in \R_+}$ punktweise gegen den
Steuerungsprozess
$$x_*(t) = e^t, \qquad u_*(t)=1, \qquad t \in \R_+.$$
Dieses Paar ist das globale Minimum der Aufgabe mit unendlichem Zeithorizont. \hfill $\square$}
\end{beispiel}
Fordern wir die Konvergenz des Zielfunktionals,
so lassen sich Aussagen zur globalen Optimalit"at treffen.
Sei $\big\{\big(x_*^T(\cdot),u_*^T(\cdot)\big)\big\}_{T \in \R_+}$ eine Familie von globalen Minimalstellen der Aufgabe
(\ref{PAUH1a})--(\ref{PAUH5a}), die folgende Hypothese erf"ullt:
\begin{enumerate}
\item[(H)] Es existiert ein zul"assiger Steuerungsprozess $\big(x_*(\cdot),u_*(\cdot)\big)$ der Aufgabe (\ref{PAUH1})--(\ref{PAUH5}) mit den Eigenschaften
           $$\lim_{T \to \infty} J_T\big(x_*^T(\cdot),u_*^T(\cdot)\big) = J\big(x_*(\cdot),u_*(\cdot)\big),
             \qquad J\big(x_*(\cdot),u_*(\cdot)\big) < \infty.$$
\end{enumerate}

In der Hypothese (H) werden an die Familie 
$\big\{\big(x_*^T(\cdot),u_*^T(\cdot)\big)\big\}_{T \in \R_+}$ keine Konvergenzeigenschaften gefordert.
Wesentlicher Punkt in der Hypothese (H) ist die Existenz des zul"assigen Steuerungsprozesses $\big(x_*(\cdot),u_*(\cdot)\big)$
mit der angegebenen Grenzwertbeziehung.
Dieser Nachweis stellt die eigentliche Herausforderung dar.

\begin{lemma} \label{LemmaKonvergenzFiniteApproximation}
Es sei $\big\{\big(x_*^T(\cdot),u_*^T(\cdot)\big)\big\}_{T \in \R_+}$ eine Familie von globalen Minimalstellen der Aufgabe
(\ref{PAUH1a})--(\ref{PAUH5a}), f"ur die die Hypothese (H) erf"ullt ist.
Dann ist $\big(x_*(\cdot),u_*(\cdot)\big)$ ein globales Minimum der Aufgabe (\ref{PAUH1})--(\ref{PAUH5}).
\end{lemma}

{\bf Beweis} Sei $\big(x(\cdot),u(\cdot)\big)$ zul"assig in der Aufgabe (\ref{PAUH1})--(\ref{PAUH5}).
Insbesondere ist dann das Zielfunktional in (\ref{PAUH1}) endlich.
Dann l"asst sich zu jedem $\varepsilon>0$ eine Zahl $T'>0$ derart angeben,
dass die Einschr"ankung des Zielfunktionals $J\big(x(\cdot),u(\cdot)\big)$
auf das Intervall $[T,\infty)$ vom Betrag kleiner oder gleich $\varepsilon$ f"ur alle $T\geq T'$ ausf"allt.
Ferner kann die Zahl $T'$ so gew"ahlt werden, dass
$\big|J_T\big(x_*^T(\cdot),u_*^T(\cdot)\big) - J\big(x_*(\cdot),u_*(\cdot)\big)\big| \leq \varepsilon$ f"ur alle $T\geq T'$ gilt.
Mit den suggestiven Bezeichnungen, dass $J_T$ bzw. $J$ die Integration des Zielfunktionals "uber $[0,T]$ bzw. "uber $\R_+$ angeben,
k"onnen wir die Differenz
$$J\big(x(\cdot),u(\cdot)\big) - J\big(x_*(\cdot),u_*(\cdot)\big)$$
in die Form
\begin{eqnarray*}
&& \Big[J\big(x(\cdot),u(\cdot)\big) - J_T\big(x(\cdot),u(\cdot)\big)\Big] 
+\Big[J_T\big(x(\cdot),u(\cdot)\big) - J_T\big(x_*^T(\cdot),u_*^T(\cdot)\big)\Big] \\
 &&   + \Big[J_T\big(x_*^T(\cdot),u_*^T(\cdot)\big) - J\big(x_*(\cdot),u_*(\cdot)\big)\Big]
\end{eqnarray*}
bringen.
Aufgrund der Wahl von $T'$ fallen der erste und dritte Summand gr"o"ser oder gleich $-\varepsilon$ f"ur alle $T\geq T'$ aus.
Der zweite Ausdruck ist nichtnegativ,
da $\big(x_*^T(\cdot),u_*^T(\cdot)\big)$ ein globales Minimum "uber $[0,T]$ darstellt.
Daher gilt $J\big(x(\cdot),u(\cdot)\big) - J\big(x_*(\cdot),u_*(\cdot)\big) \geq -2\varepsilon$.
Da $\varepsilon>0$ beliebig ist, folgt damit die Behauptung. \hfill $\blacksquare$ \\[2mm]
Bei Aseev \& Kryazhimskii \cite{AseKry} sind Voraussetzungen angegeben,
unter denen die Hypothese (H) in der Aufgabe (\ref{PAUH1})--(\ref{PAUH4}) mit freiem rechten Endpunkt erf"ullt ist.
Insbesondere wird gezeigt,
dass die Approximation durch die Aufgabe (\ref{PAUH1a})--(\ref{PAUH4a}) mit endlichem Horizont eine Familie von globalen Minimalstellen
liefert,
die auf jedem endlichen Intervall gleichm"a"sig gegen das globale Minimum
der Aufgabe (\ref{PAUH1})--(\ref{PAUH4}) konvergiert. \\[2mm]
Ein wesentlicher Augenmerk bei der Herleitung notwendiger Optimalit"atsbedingungen liegt auf der G"ultigkeit der ``nat"urlichen'' Transversalit"atsbedingungen
$$\lim_{t \to \infty} \|p(t)\| = 0, \qquad \lim_{t \to \infty} \langle p(t),x(t) \rangle = 0 \quad \mbox{f"ur alle zul"assigen } x(\cdot)$$
in der Aufgabe  (\ref{PAUH1})--(\ref{PAUH4}) mit freiem rechten Endpunkt.
In diesem Zusammenhang gibt Halkin \cite{Halkin} folgendes Gegenbeispiel an:
\begin{beispiel} {\rm In der Aufgabe
$$\int_0^\infty u(t)\big(1-x(t)\big)\, dt \to \sup, \quad \dot{x}(t)=u(t)\big(1-x(t)\big), \quad x(0)=0,\quad u \in [0,1],$$
liefert der Steuerungsprozess
$$u_*(t)=\left\{\begin{array}{ll}
  0,& t \in [0,1], \\ 1,& t \in [1,\infty),\end{array} \right. \quad
  x_*(t)=\left\{\begin{array}{ll}
  0,& t \in [0,1], \\ 1-e^{1-t},& \in [1,\infty),\end{array} \right.$$
ein globales Optimum.
F"ur die zugeh"origen Multiplikatoren $\big(\lambda_0,p(\cdot)\big)\not=(0,0)$,
mit denen die Bedingungen des Maximumprinzips erf"ullt sind,
gilt dann im Widerspruch zur ``nat"urlichen'' Transversalit"atsbedingung $p(t) \equiv -\lambda_0$. \hfill $\square$}
\end{beispiel}

Neben den Schwierigkeiten,
die bei der Approximation mit endlichem Zeithorizont entstehen,
r"uckt das Beispiel von Halkin die Frage in den Vordergrund,
welche Transversalit"atsbedingungen in der Aufgabe (\ref{PAUH1})--(\ref{PAUH5}) notwendige Optimalit"atsbedingungen darstellen.
Die Antworten k"onnen durch die Approximation mit endlichem Horizont nicht gegeben werden,
da der vollst"andige Satz an notwendigen Optimalit"atsbedingungen im Allgemeinen verloren geht.
Aus diesem Grund sind wir im letzten Kapitel zu schwachen lokalen Minimalstellen und im vorliegenden Kapitel "uber starke lokale Minimalstellen
in Aufgaben mit unendlichem Zeithorizont stets auf die G"ultigkeit der verschiedenen Transversalit"atsbedingungen eingegangen. \\
Insbesondere konnten wir im Pontrjaginschen Maximumprinzip \ref{SatzPAUH} Transversalit"atsbedingungen in der ``klassischen'' Form
$$p(0) = {h_0'}^T\big(x_*(0)\big)l_0, \qquad \lim_{t \to \infty} p(t)= - \lim_{t \to \infty}h_{1x}^T\big(t,x_*(t)\big) l_1$$
angeben,
aus denen sich die G"ultigkeit der ``nat"urlichen'' Transversalit"atsbedingungen
$$\lim_{t \to \infty} p(t) =0, \qquad \lim_{t \to \infty} \langle p(t),x(t) \rangle = 0 \mbox{ f"ur alle zul"assigen } x(\cdot)$$
und f"ur ``gutartige'' Verteilungsfunktionen $\omega(\cdot)$ die Bedingung von Michel
$$\lim_{t \to \infty} H\big(t,x_*(t),u_*(t),p(t),\lambda_0\big)=0$$
in der Aufgabe mit freiem rechten Endpunkt als direkte Konsequenz ergab. 

%%%%%%%%%%%%%%%%%%%%%%%%%%%%%%%%%%%%%

\subsection{Die Einordnung der Aufgabenklassen mit verschiedenen Horizonten}
In diesem Abschnitt gehen wir der grundlegenden Frage nach,
wie sich die Aufgabenklassen mit endlichem und mit unendlichem Zeithorizont zueinander einordnen.
Dazu gehen wir zuerst auf das Wesen des unendlichen Zeithorizontes ein: \\[2mm]
Eine Herangehensweise an die Aufgaben mit unendlichem Zeithorizont ist die Zur"uckf"uhrung auf eine Aufgabe
"uber endlichem Zeitintervall mit Hilfe der Substitution der Zeit.
Betrachten wir z.\,B. die Transformation $t(s)= -\ln(1-s)$,
so gilt
$$t'(s)=v(s)=\frac{1}{1-s}>0$$
und es wird die Aufgabe (\ref{PAUH1})--(\ref{PAUH5}) eineindeutig
in die Aufgabe (\ref{PAUH1a})--(\ref{PAUH5a}) "uber dem endlichen Intervall $[0,1]$ "uberf"uhrt. 
Schauen wir uns nach der Substitution der Zeit die Aufgabe "uber dem endlichen Horizont an,
dann besitzt sie die Form (vgl. Ioffe \& Tichomirov \cite{Ioffe})
\begin{eqnarray*}
&& \int_0^1 v(s) \cdot \omega(s) f\big(t(s),y(s),w(s)\big) \, ds \to \inf, \\
&& y'(s) = v(s) \cdot \varphi\big(t(s),y(s),w(s)\big), \qquad t'(s) = v(s), \\
&& y(0)=x_0,\quad h\big(t(1),y(1)\big)=0, \qquad t(0)=0, \quad t(1)=\infty, \\
&& w(s) \in U, \qquad v(s) > 0, \\
&& g_j\big(t(s),y(s)\big) \leq 0, \quad s \in [0,1], \quad j=1,...,l.
\end{eqnarray*}
Die Aufgabe enth"alt die Singularit"at $t(1)=\infty$,
die in s"amtlichen Elementen dieser Aufgabe einflie"st,
und au"serdem die "uber $[0,1]$ nicht integrable und unbeschr"ankte Funktion $v(\cdot)$.
Die auf $[0,1]$ "uberf"uhrte Aufgabenstellung geh"ort damit nicht in den Rahmen der klassischen Steuerungsprobleme.
Sondern die Substitution der Zeit verdeutlicht viel mehr,
dass der unendliche Horizont in der Aufgabe (\ref{PAUH1})--(\ref{PAUH5}) selbst eine Singularit"at ist,
die man nicht aus der Betrachtung argumentieren kann.
Also kann die Aufgabe (\ref{PAUH1})--(\ref{PAUH5}) kein klassisches Problem der Optimalen Steuerung sein. \\[2mm]
Damit stellt sich umgekehrt die Frage,
ob sich die Aufgabe mit endlichem Zeitintervall der Aufgabe mit unendlichem Horizont unterordnet:
Mit der Aufgabe "uber endlichem Zeithorizont bezeichnen wir das Steuerungsproblem
\begin{equation} \label{EinordnungEH} \left. \begin{array}{l}
  J\big(x(\cdot),u(\cdot)\big) = \displaystyle \int_{t_0}^{t_1} f\big(t,x(t),u(t)\big) \, dt \to \inf, \\[1mm]
  \dot{x}(t) = \varphi\big(t,x(t),u(t)\big), \\[1mm]
  h_0\big(x(t_0)\big)=0, \qquad h_1\big(x(t_1)\big)=0, \\[1mm]
  u(t) \in U \subseteq \R^m, \quad U\not= \emptyset, \\[1mm]
  g_j\big(t,x(t)\big) \leq 0 \quad \mbox{f"ur alle } t \in [t_0,t_1], \quad j=1,...,l.
  \end{array} \right\}
\end{equation}
Bez"uglich des Integranden $f$ und der rechten Seite $\varphi$ der Dynamik setzen wir:
$$\tilde{f}(t,x,u)= \left\{ \begin{array}{ll} f(t,x,u), & t \in [t_0,t_1], \\ 0, & t \not \in [t_0,t_1], \end{array} \right. \qquad
  \tilde{\varphi}(t,x,u)= \left\{ \begin{array}{ll} \varphi(t,x,u), & t \in [t_0,t_1], \\ 0, & t \not \in [t_0,t_1]. \end{array} \right.$$
Mit den Randbedingungen verkn"upfen wir die Abbildungen
$$\tilde{h}_0(x)= h_0(x), \qquad \tilde{h}_1(t,x)=h_1(x).$$
Au"serdem setzen wir bez"uglich der Zustandsbeschr"ankungen
$$\tilde{g}_j(t,x) = \left\{ \begin{array}{ll} g_j(t_0,x) - (1-e^{(t-t_0)^2}), & t < t_0, \\
                                               g_j(t,x), & t \in [t_0,t_1], \\
                                               g_j(t_1,x) - (1-e^{(t-t_1)^2}), & t > t_1. \end{array}\right.$$
Mit der Verteilungsfunktion $\omega(t) = \chi_{[t_0,t_1]}(t)$ ergibt sich auf diese Weise "uber dem unendlichen Zeithorizont die Aufgabe
\begin{equation} \label{EinordnungUH} \left. \begin{array}{l}
  J\big(x(\cdot),u(\cdot)\big) = \displaystyle \int_0^\infty \omega(t)\tilde{f}\big(t,x(t),u(t)\big) \, dt \to \inf, \\[1mm]
  \dot{x}(t) = \tilde{\varphi}\big(t,x(t),u(t)\big), \\[1mm]
  \tilde{h}_0\big(x(0)\big)=0, \qquad \displaystyle\lim_{t \to \infty} \tilde{h}_1\big(t,x(t)\big)=0, \\[1mm]
  u(t) \in U \subseteq \R^m, \quad U \not= \emptyset, \\[1mm]
  \tilde{g}_j\big(t,x(t)\big) \leq 0 \quad \mbox{f"ur alle } t \in \R_+, \quad j=1,...,l.
  \end{array} \right\}
\end{equation}
Die Aufgabe (\ref{EinordnungUH}) besitzt damit die Gestalt der Aufgabe (\ref{PAUH1})--(\ref{PAUH5}).
Bez"uglich der Verweise auf die Theoreme \ref{SatzPAUH} und \ref{SatzPAUHZA} ist anzumerken,
dass die unstetigen Anschl"usse der Abbildungen $\tilde{f}$ und $\tilde{\varphi}$ in den Stellen $t=t_0$ und $t=t_1$ sich nicht
nachteilig auf die Beweise der Maximumprinzipien auswirken. \\
Aus diesen Festlegungen ist ersichtlich,
dass ein Steuerungsprozess $\big(\tilde{x}(\cdot),\tilde{u}(\cdot)\big)$ genau dann in der Aufgabe (\ref{EinordnungUH}) zul"assig bzw. 
ein starkes lokales Minimum ist,
wenn die Einschr"ankung $\big(x(\cdot),u(\cdot)\big)$ mit 
$\big(x(t),u(t)\big)= \big(\tilde{x}(t),\tilde{u}(t)\big)$ f"ur $t \in [t_0,t_1]$
zul"assig bzw. ein starkes lokales Minimum in der Aufgabe (\ref{EinordnungEH}) darstellt.
Insbesondere stimmen dann bez"uglich der Zustandbeschr"ankungen die Mengen
$$T_j = \big\{t \in [t_0,t_1] \,\big|\, g_j\big(t,x(t)\big)=0\big\}, \qquad 
  \tilde{T}_j = \big\{t \in \overline{\R}_+=[0,\infty] \,\big|\, \tilde{g}_j\big(t,\tilde{x}(t)\big)=0\big\}$$
f"ur $j=1,...,l$ "uberein.
Bei der Anwendung der Theoreme \ref{SatzPAUH} und \ref{SatzPAUHZA} auf die Aufgabe (\ref{EinordnungUH}) sind
die restriktiven Annahmen (\ref{PMPBedingung}) und (\ref{PMPBedingung2}) zu beachten,
die der Setzung nach f"ur $\tilde{\varphi}(t,x,u)$ offensichtlich erf"ullt sind.
Daher existieren nach dem Pontrjaginschen Maximumprinzip (Theorem \ref{SatzPAUHZA}) 
f"ur ein starkes lokales Minimum $\big(x_*(\cdot),u_*(\cdot)\big)$ der Aufgabe (\ref{EinordnungUH})
eine Zahl $\lambda_0 \geq 0$, eine Vektorfunktion $p(\cdot):\R_+ \to \R^n$
und auf den Mengen
$$\tilde{T}_j=\big\{t \in \overline{\R}_+ \,\big|\, \tilde{g}_j\big(t,x_*(t)\big)=0\big\}, \quad j=1,...,l,$$
konzentrierte nichtnegative regul"are Borelsche Ma"se $\mu_j$ endlicher Totalvariation
(wobei s"amtliche Gr"o"sen nicht gleichzeitig verschwinden) derart, dass
\begin{enumerate}
\item[(a)] die Vektorfunktion $p(\cdot)$ von beschr"ankter Variation ist, der adjungierten Gleichung
           \begin{eqnarray*}
           p(t)&=&- \lim_{t \to \infty}\tilde{h}_{1x}^T\big(t,x_*(t)\big) l_1 + \int_t^\infty H_x\big(s,x_*(s),u_*(s),p(s),\lambda_0\big) \, ds \\
               & & -\sum_{j=1}^l \int_t^\infty \tilde{g}_{jx}\big(s,x_*(s)\big)\, d\mu_j(s)
           \end{eqnarray*}
           gen"ugt und die Transversalit"atsbedingungen
           \begin{eqnarray*}
           p(0) &=& \tilde{h}_0{'}^T\big(x_*(0)\big)l_0, \\
           \lim_{t \to \infty} p(t) &=& \lim_{t \to \infty} \bigg[-\tilde{h}_{1x}^T\big(x_*(t)\big)l_1
                                         -\sum_{j=1}^l \tilde{g}_{jx}\big(t,x_*(t)\big)\, \mu_j(\{\infty\})\bigg]
           \end{eqnarray*}
           erf"ullt;
\item[(b)] f"ur fast alle $t\in \R_+$ die Maximumbedingung
           $$H\big(t,x_*(t),u_*(t),p(t),\lambda_0\big) = \max_{u \in U} H\big(t,x_*(t),u,p(t),\lambda_0\big)$$
           gilt.
\end{enumerate}
In der Aufgabe ohne Zustandsbeschr"ankungen ergeben sich wegen $\dot{p}(t)=0$ f"ur $t \not \in [t_0,t_1]$ aus den Transversalit"atsbedingungen unmittelbar
$$p(t_0) = {h_0'}^T\big(x_*(t_0)\big)l_0, \qquad p(t_1)= -\tilde{h}^T_{1x}\big(t_1,x_*(t_1)\big)l_1 =-{h_1'}^T\big(x_*(t_1)\big)l_1.$$
Damit ergibt sich f"ur die Aufgabe \ref{EinordnungEH} ohne Zustandsbeschr"ankungen:

\begin{theorem}[Pontrjaginsches Maximumprinzip] 
Ist $\big(x_*(\cdot),u_*(\cdot)\big)$ ein starkes lokales Minimum der Aufgabe \ref{EinordnungEH} ohne Zustandsbeschr"ankungen,
dann existieren nicht gleichzeitig verschwindende Multiplikatoren $\lambda_0 \geq 0$,
$p(\cdot) \in W^1_\infty([t_0,t_1],\R^n)$ und $l_i \in \R^{s_i}$, $i=0,1$, derart, dass
\begin{enumerate}
\item[(a)] die Funktion $p(\cdot)$ der adjungierten Gleichung
           \begin{equation}\label{SatzPMP1}
           \dot{p}(t)=-\varphi_x^T\big(t,x_*(t),u_*(t)\big) p(t) + \lambda_0 f_x\big(t,x_*(t),u_*(t)\big)
           \end{equation}
           gen"ugt und die Transversalit"atsbedingungen
           \begin{equation}\label{SatzPMP2} 
           p(t_0)= {h_0'}^T\big(x_*(t_0)\big)l_0, \qquad p(t_1)=-{h_1'}^T\big(x_*(t_1)\big)l_1
           \end{equation}
           erf"ullt; 
\item[(b)] in fast allen Punkten $t \in [t_0,t_1]$ die Maximumbedingung
           \begin{equation}\label{SatzPMP3}
           H\big(t,x_*(t),u_*(t),p(t),\lambda_0\big) = \max_{u \in U} H\big(t,x_*(t),u,p(t),\lambda_0\big)
           \end{equation}
           gilt.
\end{enumerate}
\end{theorem}

In der Aufgabe mit Zustandsbeschr"ankungen gelten $g_j\big(t,x_*(t)\big) <0$ f"ur alle $t \not \in [t_0,t_1]$.
Daher haben wir wieder $\dot{p}(t)=0$ f"ur $t \not \in [t_0,t_1]$ und erhalten aus der linksseitigen Stetigkeit
$$p(t_1) = -{h_1'}^T\big(x_*(t_1)\big)l_1 -\sum_{j=1}^l g_{jx}\big(t_1,x_*(t_1)\big)\, \mu_j(\{t_1\}).$$
Enthalten dabei die Ma"se $\mu_j$ auch von Null verschiedene Massen,
die in den Stellen $t=t_0$ bzw. $t=t_1$ konzentriert sind,
so ergeben sich wie in Ioffe \& Tichomirov \cite{Ioffe} die Sprungbeziehungen:
$$p(t_i^+)-p(t_i) = \sum_{j=1}^l g_{jx}\big(t_i,x_*(t_i)\big)\, \mu_j(\{t_i\}), \qquad i=0,1.$$
Also ergibt sich f"ur die Aufgabe \ref{EinordnungEH}:

\begin{theorem}[Pontrjaginsches Maximumprinzip]
Ist $\big(x_*(\cdot),u_*(\cdot)\big)$ ein starkes lokales Minimum der Aufgabe  \ref{EinordnungEH},
dann existieren eine Zahl $\lambda_0 \geq 0$, eine Vektorfunktion $p(\cdot):[t_0,t_1] \to \R^n$
und auf den Mengen
$$T_j=\big\{t \in [t_0,t_1] \,\big|\, g_j\big(t,x_*(t)\big)=0\big\}, \quad j=1,...,l,$$
konzentrierte nichtnegative regul"are Borelsche Ma"se $\mu_j$ endlicher Totalvariation
(wobei s"amtliche Gr"o"sen nicht gleichzeitig verschwinden) derart, dass
\begin{enumerate}
\item[(a)] die Vektorfunktion $p(\cdot)$ von beschr"ankter Variation ist, der adjungierten Gleichung
           \begin{eqnarray}
           p(t)&=&-{h_1'}^T\big(x_*(t_1)\big)l_1 + \int_t^{t_1} H_x\big(s,x_*(s),u_*(s),p(s),\lambda_0\big) \, ds \nonumber \\
           \label{SatzPMPZA1} & & -\sum_{j=1}^l \int_t^{t_1} g_{jx}\big(s,x_*(s)\big)\, d\mu_j(s)
           \end{eqnarray}
           gen"ugt und die Transversalit"atsbedingung
           \begin{equation}
           p(t_0)= {h_0'}^T\big(x_*(t_0)\big)l_0
           \end{equation}
           erf"ullt;
\item[(b)] f"ur fast alle $t\in [t_0,t_1]$ die Maximumbedingung
           \begin{equation}
           H\big(t,x_*(t),u_*(t),p(t),\lambda_0\big) = \max_{u \in U} H\big(t,x_*(t),u,p(t),\lambda_0\big)
           \end{equation}
           gilt.
\end{enumerate}
\end{theorem}

Zusammenfassend k"onnen wir das Pontrjaginsche Maximumprinzip f"ur die Aufgaben "uber endlichem Zeithorizont aus den
Theoremen \ref{SatzPAUH} und \ref{SatzPAUHZA} vollst"andig ableiten.
Die Aufgabe mit unendlichem Zeithorizont und die vorgestellten Ergebnisse stellen somit echte Verallgemeinerungen zu der Aufgabenklasse mit endlichem
Zeithorizont dar.

\newpage       
\begin{appendix}
\lhead[\thepage \hspace*{1mm} Funktionalanalytische Hilfsmittel]{}
\rhead[]{Funktionalanalytische Hilfsmittel \hspace*{1mm} \thepage}
\section{Funktionalanalytische Hilfsmittel}
\subsection{Grundprinzipien der Funktionalanalysis}
\begin{theorem}[Satz von Hahn-Banach; Fortsetzungsversion]  \index{Satz, Darstellungssatz von Riesz!von HahnBanach@-- von Hahn-Banach}
Sei $X$ ein normierter Raum und $U$ ein Untervektorraum.
Zu jedem stetigen linearen Funktional $u^*:U \to \R$ existiert dann ein stetiges lineares Funktional Funktional
$x^*: X \to \R$ mit
$$x^*\big|_U=u^*, \qquad \|x^*\|=\|u^*\|.$$
\end{theorem}

\begin{folgerung} \label{FolgerungAnnulator}
Seien $X$ ein normierter Raum, $U$ ein abgeschlossener Unterraum und $x \in X$, $x \not\in U$.
Dann existiert ein $x^* \in X^*$ mit
$$x^*\big|_U=0,\qquad \langle x^*,x \rangle \not=0.$$
\end{folgerung}

\begin{theorem}[Satz von Hahn-Banach; Trennungsversion]  \index{Satz, Darstellungssatz von Riesz!Trennungssatz@--, Trennungssatz}
Seien $X$ ein normierter Raum, $V_1,V_2 \subseteq X$ konvex und $V_1$ offen.
Es gelte $V_1 \cap V_2 = \emptyset$.
Dann existiert ein $x^* \in X^*$ mit
$$\langle x^*, v_1 \rangle < \langle x^*, v_2 \rangle \qquad \mbox{ f"ur alle } v_1 \in V_1, v_2 \in V_2.$$
\end{theorem}

Eine Abbildung $T$ hei"st offen\index{Abbildung, beschr"ankte!offene@--, offene}, wenn $T$ offene Mengen auf offene Mengen abbildet.

\begin{theorem}[Satz von der offenen Abbildung] \label{SatzOffeneAbbildung}  \index{Satz, Darstellungssatz von Riesz!von der offen@-- von der offenen Abbildung}
Es seien $X$, $Y$ Banachr"aume und $T \in L(X,Y)$ surjektiv.
Dann ist $T$ offen.
\end{theorem}

Seien $X,Y$ normierte R"aume und $T \in L(X,Y)$.
Der adjungierte Operator\index{adjungierter Operator} $T^*:Y^* \to X^*$ ist durch
$\langle T^*y^*, x \rangle = \langle y^*, Tx \rangle$
definiert. Offensichtlich folgt daraus $T^* \in L(Y^*,X^*)$.
Seien nun $U \subseteq X$ und $V \subseteq X^*$.
Wir definieren die Mengen
\begin{eqnarray*}
U^\perp &=& \{ x^* \in X^* \,|\, \langle x^*, x \rangle =0 \mbox{ f"ur alle } x \in U\}, \\
V_\perp &=& \{ x \in X \,|\, \langle x^*, x \rangle =0 \mbox{ f"ur alle } x^* \in V\}.
\end{eqnarray*}

\begin{lemma}[Satz vom abgeschlossenen Bild] \label{SatzAbgeschlossenesBild} \index{Satz, Darstellungssatz von Riesz!vom ab@-- vom abgeschlossenen Bild}
Seien $X$, $Y$ Banachr"aume, und es sei $T \in L(X,Y)$.
Dann gelten die "Aquivalenzen:
\begin{eqnarray*}
{\rm Im\,}T \mbox{ ist abgeschlossen } &\Leftrightarrow& {\rm Im\,}T=({\rm Ker\,}T^*)_\perp \\
\Leftrightarrow \; {\rm Im\,}T^* \mbox{ ist abgeschlossen} &\Leftrightarrow&  {\rm Im\,}T^*=({\rm Ker\,}T)^\perp.
\end{eqnarray*}
\end{lemma}

\begin{satz}[Fixpunktsatz von Weissinger] \label{SatzWeissinger} \index{Satz, Darstellungssatz von Riesz!Fixpunkt@--, Fixpunktsatz von Weissinger}
Es sei $U$ eine nichtleere abgeschlossene Teilmenge des Banachraumes $X$,
ferner $\displaystyle \sum_{n=1}^\infty a_n$ eine konvergente Reihe positiver Zahlen und
$A:U \to U$ eine Selbstabbildung von $U$ mit
$$\|A^n u-A^n v\| \leq a_n\|u-v\| \qquad\mbox{ f"ur alle } u,v \in U,\; n \in \N.$$
Dann besitzt $A$ genau einen Fixpunkt, d.\,h. es gibt genau ein $u \in U$ mit $Au=u$. \\
Dieser Fixpunkt ist Grenzwert der Iterationsfolge $u_n=Au_{n-1}, n=1,2,...,$ bei beliebigem Startwert $u_0 \in U$.
Schlie"slich gilt die Fehlerabsch"atzung
$$\|u-u_n\| \leq \|u_1-u_0\| \cdot \sum_{k=n}^\infty a_k.$$
\end{satz}
%%%%%%%%%%%%%%%%%%%%%%%%%%%%%%%%%%%%%%%%%%%%%%%%%%%%%%%

\subsection{Der Darstellungssatz von Riesz} \label{AnhangRiesz}
Es sei $I \subseteq \R_+$ und es bezeichne $\mathscr{M}(I)$ die Menge der signierten regul"aren Borelschen Ma"se auf der
Borelschen $\sigma$-Algebra auf $I$.
Au"serdem bezeichnet $C_{\lim}(\R_+,\R^n)$ den Raum der stetigen Funktionen $x(\cdot)$, die im Unendlichen einen Grenzwert besitzen.
Als abgeschlossener Unterraum des Raumes $C_b(\R_+,\R^n)$ ist $C_{\lim}(\R_+,\R^n)$ vollst"andig.

\begin{satz}[Rieszscher Darstellungssatz] \label{SatzRieszC0}
Der Dualraum $C_0^*(\R_+,\R)$ ist isometrisch isomorph zu $\mathscr{M}(\R_+)$ unter der Abbildung
$$\Lambda:\mathscr{M}(\R_+) \to C_0^*(\R_+,\R), \qquad \Lambda(\mu)x(\cdot)= \int_{\R_+} x(t) \, d\mu(t).$$
\end{satz}

\begin{satz}[Rieszscher Darstellungssatz] \label{SatzRieszClim}
Der Dualraum $C_{\lim}^*(\R_+,\R^n)$ ist unter der Abbildung
$$\Lambda(\mu)x(\cdot)= \int_{\R_+} \langle x(t) , d\mu_0(t) \rangle + \lim_{t \to \infty} \langle x(t) , \mu_\infty \rangle$$
isometrisch isomorph zu den signierten Vektorma"sen $\mu \in \mathscr{M}(\overline{\R}_+)$.
Dabei besitzt $\mu$ die Darstellung $\mu=\mu_0+\mu_\infty$
mit einem $\mu_0 \in \mathscr{M}(\R_+)$ und einem in $t =\infty$ konzentrierten signierten Ma"s $\mu_\infty$.
\end{satz}

{\bf Beweis} Wir betrachten die stetige lineare Abbildung $T:C_{\lim}(\R_+,\R^n) \to C_0(\R_+,\R^n)$,
$$Tx(\cdot)(t)=x(t) - \lim_{t \to \infty} x(t).$$
Der Setzung nach bildet $T$ auf den gesamten Raum $C_0(\R_+,\R^n)$ ab.
Ferner gilt
$${\rm Ker\,}T = \{ x(\cdot) \in C_{\lim}(\R_+,\R^n;\nu) \,|\,  x(t)= \mbox{konstant}\}.$$
Sei $x^* \in C_{\lim}^*(\R_+,\R^n)$ und sei
$a \in \R^n$ mit den Komponenten $a_i = \langle x^*,e_i(\cdot) \rangle$,
wobei die $i$-te Komponente der Funktion $e_i(\cdot)$ identisch Eins ist und alle weiteren Komponenten identisch gleich Null sind.
Wir definieren das Funktional $x_1^*$ durch
$$\langle x_1^*, x(\cdot) \rangle = \langle x^*, x(\cdot) \rangle - \lim_{t \to \infty} \langle a,x(t) \rangle.$$
Dann gilt $x_1^* \in ({\rm Ker\,}T)^\perp$ und es existiert nach dem Satz vom abgeschlossenen Bild ein $y^* \in C_0^*(\R_+,\R^n)$ mit $x_1^*= T^*y^*$.
Daraus folgt mit dem Rieszschen Darstellungssatz
\begin{eqnarray*}
\langle x^*, x(\cdot) \rangle &=& \langle x_1^*, x(\cdot) \rangle + \lim_{t \to \infty} \langle a,x(t) \rangle
                                  = \langle y^*, Tx(\cdot) \rangle +\lim_{t \to \infty} \langle a,x(t) \rangle \\
              &=& \int_{\R_+}  \big\langle x(t)-\lim_{t \to \infty} x(t) , d\mu_0(t) \big\rangle + \lim_{t \to \infty} \langle a, x(t) \rangle \\
              &=& \int_{\R_+}  \langle x(t) , d\mu_0(t) \rangle + \lim_{t \to \infty} \langle x(t) , \mu_\infty \rangle.
\end{eqnarray*}
Der Darstellungssatz ist damit nachgewiesen. \hfill $\blacksquare$
       
%%%%%%%%%%%%%%%%%%%%%%%%%%%%%%%%%%%%%%%%%%%%%%%%%%%%%%%

\subsection{Der Satz von Ljusternik}
In diesem Abschnitt befassen wir uns mit dem fundamentalen Satz von Ljusternik (Ljusternik \cite{Ljusternik}).
Die vorliegende zusammenfassende Darstellung und die Verallgemeinerung ist Ioffe \& Tichomirov \cite{Ioffe} entnommen.
Eine vollst"andige Beweisf"uhrung ist wiederum in Ioffe \& Tichomirov \cite{Ioffe} zu finden.

\begin{definition}[Lokaler Tangentialkegel]
Sei $X$ ein Banachraum und $x_0 \in M \subseteq X$.
Mit $\mathscr{C}(M,x_0)$ bezeichnen wir die Menge aller Elemente $x \in X$,
zu denen ein $\varepsilon_0 > 0$ und eine Abbildung $r(\varepsilon): [0,\varepsilon_0] \to X$ mit den Eigenschaften
$$\lim_{\varepsilon \to 0^+} \frac{\| r(\varepsilon) \|}{\varepsilon} = 0
  \qquad\mbox{und}\qquad x_0 + \varepsilon x + r(\varepsilon) \in M \quad \mbox{ f"ur alle } \varepsilon \in [0,\varepsilon_0]$$
existieren. 
$\mathscr{C}(M,x_0)$ hei"st der lokale Tangentialkegel an die Menge $M$ im Punkt $x_0$.
\end{definition}

\begin{theorem}[Satz von Ljusternik] \label{SatzLjusternik}\index{Satz, Darstellungssatz von Riesz!von Ljusternik@-- von Ljusternik}
Es seien $X$ und $Y$ Banachr"aume,
$V$ eine Umgebung des Punktes $x_* \in X$ und $F$ eine Fr\'echet-differenzierbare Abbildungen der Menge $V$ in $Y$.
Wir setzen voraus,
die Abbildung $F$ sei regul"ar im Punkt $x_*$, d.\,h., es gelte
$${\rm Im\,} F'(x_*) = Y,$$
au"serdem sei ihre Ableitung in diesem Punkt in der gleichm"a"sigen Operatorentopologie des Raumes $L(X,Y)$ stetig. \\
Unter diesen Voraussetzungen stimmt dann der lokale Tangentialkegel an die Menge
$$M= \big\{ x \in V \,\big|\, F(x) = F(x_*) \big\}$$
im Punkt $x_*$ mit dem Kern des Operators $F'(x_*)$ "uberein:
$$\mathscr{C}(M,x_*) = {\rm Ker\,} F'(x_*).$$
\end{theorem}

\begin{lemma} \label{LemmaVerallgemeinerterSatzLjusternik}
Es seien $X$ und $Y$ Banachr"aume und $\Lambda:X \to Y$ ein stetiger linearer Operator.
Wir setzen
$$C(\Lambda) := \sup_{y \not= 0} \left( \frac{\inf \{ \|x\| \,|\, x \in X, \Lambda x = y\} }{\|y\|} \right).$$
Wenn ${\rm Im\,} \Lambda = Y$ gilt, so ist $C(\Lambda) < \infty$.
\end{lemma}

\begin{theorem}[Der verallgemeinerte Satz von Ljusternik] \label{SatzVerallgemeinerterSatzLjusternik}
Es seien $X$ und $Y$ Banachr"aume,
$\Lambda: X \to Y$ ein stetiger linearer Operator und $F$ eine Abbildung einer gewissen Umgebung $V$ des Punktes $x_* \in X$ in $Y$.
Wir setzen voraus, es sei ${\rm Im\,}\Lambda = Y$ und es gebe eine Zahl $\delta > 0$ derart, dass erstens
$$\delta C(\Lambda) < \frac{1}{2}$$
und zweitens
\begin{equation} \label{SatzVerallgemeinerterSatzLjusternik1}
\| F(x) - F(x') - \Lambda(x-x') \| \leq \delta \|x-x'\|
\end{equation}
f"ur alle $x,x'$ aus $V$ gilt.
Dann existieren eine Umgebung $V' \subseteq V$ des Punktes $x_*$,
eine Zahl $K>0$ und eine Abbildung $\xi \to x(\xi)$ der Umgebung $V'$ in $X$
derart, dass die Beziehungen
$$F\big( \xi + x(\xi) \big) = F(x_*) \qquad\mbox{ und }\qquad \| x(\xi) \| \leq K \| F(\xi) - F(x_*) \|$$
f"ur alle $\xi \in V'$ erf"ullt sind.
\end{theorem}

%%%%%%%%%%%%%%%%%%%%%%%%%%%%%%%%%%%%%%%%%%%%%%%%%%%%%%%

\subsection{Differenzierbarkeit konkreter Abbildungen}
In den folgenden Beispielen verwenden wir nach dem Hauptsatz der Differential- und Integralrechnung
f"ur eine differenzierbare Abbildung $h(y):\R^n \to \R^m$ die Beziehung:
$$h(y+ h) - h(y)= \int_0^1 h'(y+ s h) y \, ds.$$

\begin{beispiel} \label{DiffAbbildung}
{\rm Sei $I \subseteq \overline{\R}$.
Zu der Funktion $y_*(\cdot) \in C_{\lim}(I,\R^n)$ definieren wir die Menge
$$V_\gamma = \{ (t,y) \in I \times \R^n \,|\, \|y-y_*(t)\| \leq \gamma\}.$$
Wir nehmen an,
dass die Abbildung $g(t,y) : \R \times \R^n \to \R^m$ auf der Menge $V_\gamma$
gleichm"a"sig stetig und bez"uglich $y$ gleichm"a"sig stetig differenzierbar ist. \\
Dann ist die Abbildung $G: C_{\lim}(I,\R^n) \to C_{\lim}(I,\R^m)$,
$G\big(y(\cdot)\big)(t) = g\big(t,y(t)\big)$,
im Punkt $y_*(\cdot)$ stetig Fr\'echet-differenzierbar und es gilt
$$\big[G'\big(y_*(\cdot)\big) y(\cdot)\big](t) = g_y\big(t,y_*(t)\big) y(t), \quad t \in I.$$
Denn: F"ur $t \in I$, $\|y(\cdot)-y_*(\cdot)\|_\infty < \gamma$ und $0 < \lambda <\lambda_0$ gilt
$$\bigg[\frac{G\big(y(\cdot) + \lambda x(\cdot)\big) - G\big(y(\cdot)\big)}{\lambda} - G'\big(y(\cdot)\big)x(\cdot)\bigg](t)
  = \hspace*{-1mm}\int_0^1  \hspace*{-1mm} \big[g_y\big(t,y(t)+s\lambda x(t)\big) - g_y\big(t,y(t)\big) \big] x(t) ds.$$
Da die Abbildung $g_y(t,y)$ auf der Menge $V_\gamma$ gleichm"a"sig stetig ist,
gibt es ein $C>0$ mit 
$$\bigg\| \frac{G\big(y(\cdot) + \lambda x(\cdot)\big) - G\big(y(\cdot)\big)}{\lambda}
   - G'\big(y(\cdot)\big)x(\cdot) \bigg\|_\infty
  \leq \sup_{t \in I} \int_0^1 C \|\lambda x(t)\| ds \cdot \|x(\cdot)\|_\infty = \lambda C,$$
d.\,h. der Grenzwert $\lambda \to 0^+$ konvergiert gleichm"a"sig bez"uglich $\|x(\cdot)\|_\infty = 1$.
Also ist die Abbildung $G$ auf einer Umgebung von $y_*(\cdot)$ Fr\'echet-differenzierbar. \\
Weiterhin ergibt sich wegen der gleichm"a"sigen Stetigkeit von $g_y(t,y)$ auf der Menge $V_\gamma$
f"ur die Abbildung $y(\cdot) \to G'\big(y(\cdot)\big)$ bez"uglich der Operatornorm in $y_*(\cdot)$:
\begin{eqnarray*}
    \big\|G'\big(y(\cdot)\big) - G'\big(y_*(\cdot)\big) \big\|
&=& \sup_{\|x(\cdot)\|_\infty=1} \big\| \big[G'\big(y(\cdot)\big) - G'\big(y_*(\cdot)\big)\big] x(\cdot) \big\|_\infty \\
&\leq& \sup_{t \in I} \big\|  g_y\big(t,y(t)\big) - g_y\big(t,y_*(t)\big) \big\|
       \leq C \|y(\cdot)-y_*(\cdot) \|_\infty.
\end{eqnarray*}
Somit ist die stetige Fr\'echet-Differenzierbarkeit nachgewiesen. \hfill $\square$}
\end{beispiel}

\begin{beispiel} \label{DiffZielfunktional2}
{\rm Es seien $[t_0,t_1] \subset \R$, $U \subseteq \R^m$, $x_*(\cdot) \in C([t_0,t_1],\R^n)$ und
$$V_\gamma = \{ (t,x) \in [t_0,t_1] \times \R^n \,|\, \|x-x_*(t)\| \leq \gamma\}.$$
Ferner sei die Funktion $f(t,x,u): [t_0,t_1] \times \R^n \times \R^m \to \R$ f"ur jede kompakte Menge $U_1 \subseteq \R^m$
auf der Menge $V_\gamma \times U_1$ stetig und bez"uglich $x$ stetig differenzierbar.
Dann ist f"ur jedes $u(\cdot) \in L_\infty([t_0,t_1],U)$ die Abbildung
$$x(\cdot) \to J\big(x(\cdot),u(\cdot)\big) = \int_{t_0}^{t_1} f\big(t,x(t),u(t)\big) \, dt, \quad J: C([t_0,t_1],\R^n) \times L_\infty([t_0,t_1],U) \to \R,$$
im Punkt $x_*(\cdot)$ Fr\'echet-differenzierbar und es gilt
$$J_x\big(x_*(\cdot),u(\cdot)\big) x(\cdot) = \int_{t_0}^{t_1} \big\langle f_x\big(t,x_*(t),u(t)\big), x(t) \big\rangle \, dt.$$
Denn: Wegen $u(\cdot) \in L_\infty([t_0,t_1],U)$ l"asst sich eine kompakte Menge $U_1 \subseteq \R^m$ mit $u(t) \in U_1$ f"ur fast alle $t$ angeben.
Nach Voraussetzung ist die Abbildung $f(t,x,u)$ auf der Menge $V_\gamma \times U_1$ stetig und stetig differenzierbar bez"uglich $x$.
Daher ist $J$ wohldefiniert.
Weiterhin ist f"ur jedes $u(\cdot) \in L_\infty([t_0,t_1],U)$ der lineare Operator $J_x\big(x_*(\cdot),u(\cdot)\big)$ stetig. \\
Sei $\varepsilon >0$ gegeben.
Aufgrund der gleichm"a"sigen Stetigkeit von $f_x(t,x,u)$ auf der Menge $V_\gamma \times U_1$ gibt es eine Zahl $\lambda_0>0$ mit
$$\big\|f_x\big(t,x_*(t)+\lambda x,u\big)-f_x\big(t,x_*(t),u\big)\big\| \leq \frac{\varepsilon}{t_1-t_0}$$
f"ur fast alle $t \in [t_0,t_1]$ und f"ur alle $\|x\| \leq 1$, $u \in U_1$, $0 < \lambda \leq \lambda_0$.
Damit erhalten wir
\begin{eqnarray*}
\lefteqn{\bigg|\frac{J\big(x_*(\cdot)+\lambda x(\cdot),u(\cdot)\big) - J\big(x_*(\cdot),u(\cdot)\big)}{\lambda}
               - J_x\big(x_*(\cdot),u(\cdot)\big)x(\cdot)\bigg|} \\
&=& \bigg|\int_{t_0}^{t_1} \bigg[\int_0^1 \langle f_x\big(t,x_*(t)+\lambda s x(t),u(t)\big) - f_x\big(t,x_*(t),u(t)\big), x(t) \rangle \, ds \bigg] \, dt \bigg|
    \leq \varepsilon
\end{eqnarray*}
f"ur alle $\|x(\cdot) \|_\infty \leq 1$ und alle $0 < \lambda \leq \lambda_0$. \hfill $\square$}
\end{beispiel}

\begin{beispiel} \label{DiffZielfunktionalS}
{\rm Es seien $U \subseteq \R^m$, $\omega(\cdot)$ eine Gewichtsfunktion aus dem Raum $L_1(\R_+,\R_+)$,
weiterhin $x_*(\cdot) \in C_{\lim}(\R_+,\R^n)$ und
$$V_\gamma= \{ (t,x) \in \overline{\R}_+ \times \R^n \,|\, \|x-x(t)\| \leq \gamma\}.$$
Ferner sei die Funktion $f(t,x,u): \R \times \R^n \times \R^m \to \R$ f"ur jede kompakte Menge $U_1 \subseteq \R^m$
auf der Menge $V_\gamma \times U_1$ gleichm"a"sig stetig und bez"uglich $x$ gleichm"a"sig stetig differenzierbar.
Dann ist f"ur jedes $u(\cdot) \in L_\infty(\R_+,U)$ die Abbildung 
$$x(\cdot) \to J\big(x(\cdot),u(\cdot)\big) = \int_0^\infty \omega(t) f\big(t,x(t),u(t)\big) \, dt,$$
wobei $J:C_{\lim}(\R_+,\R^n) \times L_\infty(\R_+,U) \to \R$ gilt,
im Punkt $x_*(\cdot)$ Fr\'echet-differenzierbar und
$$J_x\big(x_*(\cdot),u(\cdot)\big) x(\cdot) = \int_0^\infty \omega(t) \big\langle f_x\big(t,x_*(t),u(t)\big), x(t) \big\rangle \, dt.$$
Denn: Wegen $u(\cdot) \in L_\infty(\R_+,U)$, $\omega(\cdot) \in L_1(\R_+,\R_+)$ und wegen der gleichm"a"sigen Differenzierbarkeit lassen sich
zu jedem $\varepsilon>0$ Zahlen $\lambda_0,T>0$ und nach dem Satz von Lusin eine kompakte Menge $K \subseteq [0,T]$ derart angeben,
dass $\omega(\cdot)$ auf $K$ stetig ist, die Ungleichung
$$\int_{\R_+ \setminus K} \omega(t) \,dt \leq \frac{\varepsilon}{2}$$
gilt und die Beziehung
$$\big\|f_x\big(t,x_*(t)+\lambda x,u(t)\big)-f_x\big(t,x_*(t),u(t)\big)\big\| \leq 1$$
f"ur fast alle $t\not \in K$ und f"ur alle $\|x\| \leq 1$, $0 \leq \lambda \leq \lambda_0$ erf"ullt ist.
Benutzen wir "uber der kompakten Menge $K$ die gleiche Argumentation wie in Beispiel \ref{DiffZielfunktional2},
so erhalten wir
\begin{eqnarray*}
\lefteqn{\bigg|\frac{J\big(x_*(\cdot)+\lambda x(\cdot),u(\cdot)\big) - J\big(x_*(\cdot),u(\cdot)\big)}{\lambda}
               - J_x\big(x_*(\cdot),u(\cdot)\big)x(\cdot)\bigg|} \\
&=& \bigg|\int_0^\infty \omega(t)
          \bigg[\int_0^1 \langle f_x\big(t,x_*(t)+\lambda s x(t),u(t)\big) - f_x\big(t,x_*(t),u(t)\big), x(t) \rangle \, ds \bigg] \, dt \bigg|
    \leq \varepsilon
\end{eqnarray*}
f"ur alle $\|x(\cdot) \|_\infty \leq 1$ und alle $0 < \lambda \leq \lambda_0$. \hfill $\square$}
\end{beispiel}

\begin{beispiel} \label{DiffDynamik2}
{\rm Es seien $[t_0,t_1] \subset \R$, $U \subseteq \R^m$, $x_*(\cdot) \in C([t_0,t_1],\R^n)$ und
$$V_\gamma = \{ (t,x) \in [t_0,t_1] \times \R^n \,|\, \|x-x_*(t)\| \leq \gamma\}.$$
Ferner sei die Abbildung $\varphi(t,x,u): [t_0,t_1] \times \R^n \times \R^m \to \R^n$ f"ur jede kompakte Menge $U_1 \subseteq \R^m$
auf der Menge $V_\gamma \times U_1$ stetig und bez"uglich $x$ stetig differenzierbar.
Dann ist f"ur jedes $u(\cdot) \in L_\infty([t_0,t_1],U)$ die Abbildung 
$$x(\cdot) \to F\big(x(\cdot),u(\cdot)\big)(t) = \int_{t_0}^t \varphi\big(s,x(s),u(s)\big) \, ds, \quad t \in [t_0,t_1],$$
$F: C([t_0,t_1],\R^n) \times L_\infty([t_0,t_1],\R^n) \to C([t_0,t_1],\R^n)$,
im Punkt $x_*(\cdot)$ Fr\'echet-differenzierbar und es gilt
$$\big[F_x\big(x_*(\cdot),u(\cdot)\big) x(\cdot)\big](t) = \int_{t_0}^t  \varphi_x\big(s,x_*(s),u(s)\big) x(s) \, ds, \quad t \in [t_0,t_1].$$
Denn es l"asst sich zu $u(\cdot) \in L_\infty([t_0,t_1],U)$ eine kompakte Menge $U_1 \subseteq \R^m$ mit $u(t) \in U_1$ 
f"ur fast alle $t$ angeben,
so dass $\varphi_x(t,x,u)$ auf der Menge $V_\gamma \times U_1$ gleichm"a"sig stetig ist.
Daraus ergibt sich die Fr\'echet-Differenzierbarkeit. \hfill $\square$}
\end{beispiel}

\begin{beispiel} \label{DiffDynamikS}
{\rm Es seien $U \subseteq \R^m$, $x_*(\cdot) \in C_{\lim}(\R_+,\R^n)$, $u_*(\cdot) \in L_\infty(\R_+,U)$, und
$$V_\gamma= \{ (t,x) \in \overline{\R}_+ \times \R^n \,|\, \|x-x(t)\| \leq \gamma\}.$$
Ferner sei die Abbildung $\varphi(t,x,u): \R \times \R^n \times \R^m \to \R$ f"ur jede kompakte Menge $U_1 \subseteq \R^m$
auf der Menge $V_\gamma \times U_1$ gleichm"a"sig stetig und bez"uglich $x$ gleichm"a"sig stetig differenzierbar.
Weiterhin nehmen wir an,
dass 
$$\int_0^\infty \big\|\varphi\big(t,x_*(t),u_*(t)\big)\big\| \, dt < \infty, \qquad \int_0^\infty \big\|\varphi_x\big(t,x_*(t),u_*(t)\big)\big\| \, dt < \infty$$
gelten.
Au"serdem m"oge zu jedem $\delta>0$ ein $T>0$ existieren, dass
\begin{eqnarray*}
&& \int_T^\infty \big\| \varphi\big(t,x(t),u_*(t)\big)-\varphi\big(t,x'(t),u_*(t)\big) \nonumber \\
&& \hspace*{20mm} - \varphi_x\big(t,x_*(t),u_*(t)\big)\big(x(t)-x'(t)\big) \big\| \, dt \leq \delta \|x(\cdot)-x'(\cdot)\|_\infty
\end{eqnarray*}
f"ur alle $x(\cdot), x'(\cdot) \in C_{\lim}(\R_+,\R^n)$ mit $\|x(\cdot)-x_*(\cdot)\|_\infty < \gamma$, $\|x'(\cdot)-x_*(\cdot)\|_\infty < \gamma$ erf"ullt ist.
Au"serdem bezeichnen wir mit $\mathscr{U}$ die Menge:
\begin{eqnarray*}
\mathscr{U}=\big\{ u(\cdot)  \in L_\infty(\R_+,U) &\big|& u(t)=u_*(t) + \chi_M(t)\big(w(t)-u_*(t)\big), \; w(\cdot) \in L_\infty(\R_+,U), \\
                                                  &     & M \subset \R_+ \mbox{ me"sbar und beschr"ankt} \big\}.
\end{eqnarray*}
Die Menge ist verbunden mit Nadelvariationen des Steuerungsprozesses $\big(x_*(\cdot),u_*(\cdot)\big)$,
die nur "uber beschr"ankten Mengen erfolgen. \\[2mm]
Unter diesen Voraussetzungen ist f"ur jedes $u(\cdot) \in \mathscr{U}$ die Abbildung 
$$x(\cdot) \to F\big(x(\cdot),u(\cdot)\big) = \int_0^t \varphi\big(s,x(s),u(s)\big) \, ds, \quad t \in \R_+,$$
$F:C_{\lim}(\R_+,\R^n) \times L_\infty(\R_+,U) \to C_{\lim}(\R_+,\R^n)$,
in $x_*(\cdot)$ Fr\'echet-differenzierbar und
$$\big[F_x\big(x_*(\cdot),u(\cdot)\big) x(\cdot)\big](t) = \int_0^t  \varphi_x\big(s,x_*(s),u(s)\big) x(s) \, ds, \quad t \in \R_+.$$
Denn: Wegen $u(\cdot) \in \mathscr{U}$ k"onnen wir eine Zahl $T>0$ angeben mit $u(t)=u_*(t)$ f"ur alle $t\geq T$.
Daher ergibt sich f"ur alle $x(\cdot) \in C_{\lim}(\R_+,\R^n)$ mit $\|x(\cdot)-x_*(\cdot)\|_\infty < \gamma$
$$\int_0^\infty \big\|\varphi\big(t,x(t),u(t)\big)\big\| \, dt
   \leq \int_0^T \big\|\varphi\big(t,x(t),u(t)\big)\big\| \, dt + \int_T^\infty \big\|\varphi\big(t,x(t),u_*(t)\big)\big\| \, dt.$$
Darin ergibt sich nach Voraussetzung au"serdem
$$\int_T^\infty \big\|\varphi\big(t,x(t),u_*(t)\big)\big\| \, dt
   \leq \int_T^\infty \Big[\big\|\varphi\big(t,x_*(t),u_*(t)\big)\big\| + \gamma\big\|\varphi_x\big(t,x_*(t),u_*(t)\big)\big\|\Big] \, dt
   + \delta\gamma.$$
Damit bildet der Operator $F$ in den Raum $C_{\lim}(\R_+,\R^n)$ ab.
Ferner k"onnen wir zu jedem $\varepsilon >0$ Zahlen $\lambda_0, T>0$ angeben mit $u(t)=u_*(t)$ f"ur alle $t\geq T$ und 
$$\int_T^\infty \big\| \varphi\big(t,x_*(t)+\lambda x(t),u_*(t)\big)-\varphi\big(t,x_*(t),u_*(t)\big)
                       - \varphi_x\big(t,x_*(t),u_*(t)\big)\big(\lambda x(t)\big) \big\| \, dt \leq \frac{\varepsilon}{2}$$
f"ur alle $0<\lambda \leq \lambda_0$ und alle $\|x(\cdot)\|_\infty < 1$ ausf"allt. \\
Nutzen wir "uber dem Intervall $[0,T]$ die Argumentation aus Beispiel \ref{DiffDynamik2},
so erhalten wir
$$\bigg\| \frac{F\big(x_*(\cdot)+\lambda x(\cdot),u(\cdot)\big) - F\big(x_*(\cdot),u(\cdot)\big)}{\lambda}
          - F_x\big(x_*(\cdot),u(\cdot)\big) x(\cdot)\bigg\|_\infty \leq \varepsilon$$
f"ur alle $\|x(\cdot) \|_\infty \leq 1$ und alle $0 < \lambda \leq \lambda_0$. \hfill $\square$}
\end{beispiel}

\lhead[\thepage \hspace*{1mm} Differentialgleichungen]{}
\rhead[]{Differentialgleichungen \hspace*{1mm} \thepage}
\section{Differentialgleichungen}
Wir betrachten f"ur $A(t) \in \R^{n \times n}$ und $a(t) \in \R^n$ das lineare Differentialgleichungssystem
$$\dot{x}(t)=A(t)x(t)+a(t).$$

\begin{lemma} \label{LemmaDGL}
Es sei vorausgesetzt,
dass die Abbildung $t \to A(t)$ und die Vektorfunktion $a(\cdot)$ "uber $[t_0,t_1]$ integrierbar.
Dann existiert zu jeder Vektorfunktion $\zeta(\cdot) \in C([t_0,t_1],\R^n)$ und jedem $\tau \in [t_0,t_1]$
eine eindeutig bestimmte Vektorfunktion $x(\cdot) \in C([t_0,t_1],\R^n)$ derart,
dass f"ur alle $t \in [t_0,t_1]$ die folgende Gleichung erf"ullt ist:
$$\zeta(t)=x(t) - \int_\tau^t \big[A(s)x(s)+a(s)\big] \, ds.$$
\end{lemma}

{\bf Beweis} Wir werden im Folgenden zeigen,
dass die Fixpunktgleichung $x(\cdot) = T\big(x(\cdot)\big)$, wobei der Operator $T$ durch
$$x(\cdot) \to T\big(x(\cdot)\big), \quad
  T\big(x(\cdot)\big)(t) = \zeta(t) + \int_\tau^t \big[A(s)x(s)+a(s)\big] \, ds, \quad t \in [t_0,t_1],$$
gegeben wird, stets eine eindeutige L"osung besitzt.
Der Operator $T$ bildet den Raum $C([t_0,t_1],\R^n)$ in sich ab.
Zur abk"urzenden Schreibweise seien 
$$c(t) = \|A(t)\|, \qquad C(t) = \int_\tau^t c(s) \, ds, \qquad c_0=\int_{t_0}^{t_1} c(s) \, ds.$$

Bei mehrfacher Anwendung des Operators $T$ ergeben sich f"ur $x_1(\cdot),x_2(\cdot) \in C([t_0,t_1],\R^n)$ die Beziehungen
\begin{eqnarray*}
\lefteqn{\big\| \big[T\big(x_1(\cdot) - x_2(\cdot)\big)\big](t) \big\|
         \leq \int_\tau^t c(s) \, ds \cdot \| x_1(\cdot) - x_2(\cdot) \|_\infty,} \\
\lefteqn{\big\| \big[T^2\big(x_1(\cdot) - x_2(\cdot)\big)\big](t) \big\|
         \leq \int_\tau^t c(s) \big\| \big[T\big(x_1(\cdot) - x_2(\cdot)\big)\big](s) \big\| \, ds} \\
&& \hspace*{10mm} \leq \int_\tau^t c(s) C(s) \, ds \cdot \big\| x_1(\cdot) - x_2(\cdot) \big\|_\infty
   = \frac{1}{2} C^2(t) \cdot \| x_1(\cdot) - x_2(\cdot) \|_\infty, \\
&& \hspace*{30mm}\vdots \\
\lefteqn{\big\| \big[T^m \big(x_1(\cdot) - x_2(\cdot)\big)\big](t) \big\|
         \leq \int_\tau^t c(s) \big\| \big[T^{m-1}\big(x_1(\cdot) - x_2(\cdot)\big)\big](s) \big\| \, ds} \\
&& \hspace*{10mm} \leq \int_\tau^t c(s) \frac{C^{m-1}(t)}{(m-1)!} \, ds \cdot \big\| x_1(\cdot) - x_2(\cdot) \big\|_\infty
       =\frac{C^m(t)}{m!} \cdot \| x_1(\cdot) - x_2(\cdot) \|_\infty.
\end{eqnarray*}
In der Topologie des Raumes $C([t_0,t_1],\R^n)$ gilt daher
$$\big\|T^m \big(x_1(\cdot) - x_2(\cdot)\big) \big\|_\infty \leq \frac{c_0^m}{m!} \cdot \| x_1(\cdot) - x_2(\cdot) \|_\infty.$$
Die Zahlen $a_m = \frac{c_0^m}{m!}$ liefern eine Folge, deren Reihe konvergiert.
Nach dem Fixpunktsatz von Weissinger (Satz \ref{SatzWeissinger}) existiert daher genau ein $x(\cdot)$ mit $x(\cdot) = T x(\cdot)$.
\hfill $\blacksquare$

\begin{lemma} \label{FolgerungDGL}
Es seien die Abbildung $t \to A(t)$ und die Vektorfunktion $a(\cdot)$ "uber $\R_+$ integrierbar.
Dann existiert zu jeder Vektorfunktion $\zeta(\cdot) \in C_{\lim}(\R_+,\R^n)$ und jedem $\tau \in \overline{\R}_+$
eine eindeutig bestimmte Vektorfunktion $x(\cdot) \in C_{\lim}(\R_+,\R^n)$ derart,
dass f"ur alle $t \in \R_+$ die folgende Gleichung erf"ullt ist:
$$\zeta(t)=x(t) - \int_\tau^t \big[A(s)x(s)+a(s)\big] \, ds.$$
\end{lemma}

{\bf Beweis} Wir ersetzen im Beweis von Lemma \ref{LemmaDGL} das Intervall $[t_0,t_1]$ durch $\R_+$.
Ferner seien der Operator $T$ wie dort definiert und
$$c(t) = \|A(t)\|, \qquad C(t) = \int_\tau^t c(s) \, ds, \qquad c_0 = \int_0^\infty c(s) \, ds.$$
Bei mehrfacher Anwendung des Operators $T$ ergibt sich f"ur $x_1(\cdot),x_2(\cdot) \in C_{\lim}(\R_+,\R^n)$:
$$\big\|T^m \big(x_1(\cdot) - x_2(\cdot)\big) \big\|_\infty \leq \frac{c_0^m}{m!} \cdot \| x_1(\cdot) - x_2(\cdot) \|_\infty.$$
Die Zahlen $a_m = \frac{c_0^m}{m!}$ liefern eine Folge, deren Reihe konvergiert.
Nach dem Fixpunktsatz von Weissinger existiert daher genau ein $x(\cdot)$ mit $x(\cdot) = T x(\cdot)$.
\hfill $\blacksquare$

\begin{folgerung} \label{FolgerungDGL2}
Es seien die Abbildung $A(\cdot)$ und die Vektorfunktion $a(\cdot)$ auf $\R_+$ integrierbar.
Dann existiert zu jedem $\zeta \in \R^n$ und jedem $\tau \in \R_+$ bzw. f"ur $t \to \infty$
eine eindeutig bestimmte L"osung $x(\cdot)$ der Gleichung
$$\dot{x}(t)=A(t)x(t)+a(t), \qquad x(\tau)=\zeta \mbox{ bzw.} \lim_{t \to \infty} x(t)=\zeta.$$
\end{folgerung}

{\bf Beweis} Die Behauptung folgt unmittelbar aus Lemma \ref{FolgerungDGL},
wenn wir $\zeta(t) \equiv \zeta$ setzen. \hfill $\blacksquare$
\lhead[\thepage \hspace*{1mm} Konvexe Analysis]{}
\rhead[]{Konvexe Analysis \hspace*{1mm} \thepage}
\section{Elemente der Konvexen Analysis} \label{AbschnittKonvexeAnalysis}
Bei der Zusammenstellung der grundlegenden Ergebnisse beschr"anken wir uns auf die Eigenschaften konvexer und lokalkonvexer Funktionen
nach Clarke \cite{Clarke}, Ioffe \& Tichomirov \cite{Ioffe} und Rockafellar \cite{Rockafellar}.
Im vorliegenden Rahmen stimmen die klassische Richtungsableitung und der Clarkesche Gradient "uberein.
Deswegen verweisen wir bez"uglich Lemma \ref{LemmaRichtungsableitung} und bez"uglich der Kettenregel \ref{SatzKettenregel}
auf \cite{Ioffe}.

\subsection{Das Subdifferential konvexer Funktionen}
Es seien $X,Y$ Banachr"aume.
Eine Funktion $f$ auf $X$ ist in der Konvexen Analysis eine Abbildung in die erweiterte reelle Zahlengerade, d.\,h.
$f:X \to \overline{\R} = [-\infty,\infty]$.
Der effektive Definitionsbereich
der Abbildung $f$ ist die Menge ${\rm dom\,}f = \{ x \in X | f(x) < \infty\}$. 
Die Funktion $f$ hei"st eigentlich\index{Funktion, absolutstetige!eigentliche@--, eigentliche},
falls ${\rm dom\,}f \not= \emptyset$ und $f(x) > -\infty$ f"ur alle $x \in X$ gelten. \\
Die eigentliche Funktion $f$ hei"st konvex\index{Funktion, absolutstetige!konvexe@--, konvexe},
wenn f"ur alle $x_1,x_2 \in X$ und alle $0 \leq \alpha \leq 1$ folgende Ungleichung gilt:
$$f\big(\alpha x_1 + (1-\alpha) x_2 \big) \leq \alpha f(x_1) + (1-\alpha) f(x_2).$$
Die Funktion $f$ hei"st homogen\index{Funktion, absolutstetige!homogene@--, homogene},
falls $f(0) = 0$ und $f(\lambda x) = \lambda f(x)$ f"ur alle $x \in X,\, \lambda > 0$ ist.
Eine eigentliche konvexe Funktion ist genau dann in einem Punkt stetig,
wenn sie auf einer Umgebung dieses Punktes nach oben beschr"ankt ist.
In diesem Fall ist das Innere des effektiven Definitionsbereichs nichtleer.
Ist andererseits eine homogene Funktion auf einer Umgebung des Nullpunktes stetig, so ist sie auf $X$ stetig. \\
Ist $f$ eine eigentliche konvexe Funktion auf $X$,
dann existiert in jedem Punkt der Menge ${\rm dom\,}f$ die klassische
Richtungsableitung\index{Ableitung, Fr\'echet-!Richtungs@--, Richtungs-}\label{Richtungsableitung},
d.\,h. f"ur alle $z \in X$ der Grenzwert
$$f'(x;z)=\lim_{\lambda \to 0^+} \frac{f(x + \lambda z) - f(x)}{\lambda}.$$
Sei $f$ eigentlich, konvex und in $x$ stetig.
Dann ist $f$ auf einer Umgebung des Punktes $x$ nach oben beschr"ankt,
in $x$ lokal Lipschitz-stetig und es existiert der Clarkesche Gradient\index{Clarkescher Gradient}\label{ClarkescherGradient}
$$f^\circ(x;z) = \limsup_{\substack{y \to x \\ \lambda \to 0^+}} \frac{f(y + \lambda z) - f(y)}{\lambda}.$$
Unter diesen Voraussetzungen ist au"serdem die Funktion $f$ ist im Punkt $x$ regul"ar\index{Funktion, absolutstetige!regul"are@--, regul"are}
im Sinn der Konvexen Analysis,
d.\,h. $f'(x;\cdot)=f^\circ(x;\cdot)$.
Das Subdifferential\index{Subdifferential} der eigentlichen konvexen Funktion $f$ besteht im Punkt $x$ aus allen Subgradienten
$x^* \in X^*$, d.\,h. \label{Subdifferential}
$$\partial f(x) = \{ x^* \in X^* | f(z) - f(x) \geq \langle x^*, z-x \rangle \mbox{ f"ur alle } z \in X \}.$$
F"ur eine eigentliche konvexe Funktion $f$ gilt $\partial f(x)=\partial f'(x;0)$ f"ur alle $x \in {\rm dom\,}f$.
Ist $f$ eine eigentliche homogene konvexe Funktion und $x \not= 0$, dann ist
$$\partial f(x) = \{ x^* \in \partial f(0) | f(x) = \langle x^*, x \rangle \}.$$
Mit $\partial_x f(x,y)$ bezeichnen wir das Subdifferential der Abbildung $x \to f(x,y)$.\label{partSubdifferential}

\subsection{Lokalkonvexe Funktionen}
Es seien $X,Y$ Banachr"aume.
Eine auf $X$ definierte Funktion $G$ hei"st im Punkt $x_0$ lokalkonvex\index{Funktion, absolutstetige!lokalkonvexe@--, lokalkonvexe},
wenn ihre Richtungsableitung in diesem Punkt existiert und $x \to G'(x_0;x)$ konvex ist.
Im Folgenden seien $g: X \to Y$ im Punkt $x_0 \in X$ Fr\'echet-differenzierbar und
$f: Y \to \overline{\R}$ eigentlich, konvex und im Punkt $g(x_0)$ stetig.

\begin{lemma} \label{LemmaRichtungsableitung}
Die Funktion $G: X \to \overline{\R}$, $G(x) = f\big( g(x)\big)$, besitzt in $x_0$ eine klassische Richtungsableitung, es gilt
$$G'(x_0;x) = f'\big(g(x_0);g'(x_0)x\big)$$
und die Richtungsableitung konvergiert bez"uglich jeder Richtung $x$
gleichm"a"sig\index{Funktion, absolutstetige!gleichmassig@--, gleichm"a"sig differenzierbare}:
$$\bigg| \frac{G(x_0+\lambda z) - G(x_0)}{\lambda} - G'(x_0;x) \bigg| < \varepsilon \quad
  \mbox{f"ur alle } z \in U(x), \; 0<\lambda<\lambda_0.$$
Insbesondere folgt aus der gleichm"a"sigen Differenzierbarkeit bez"uglich aller Richtungen,
dass die Richtungsableitung der Abbildung $G$ eine stetige Funktion ist.
\end{lemma}

\begin{lemma} \label{LemmaClarke}
Die Funktion $G(x) = f\big(g(x)\big)$ ist in $x_0$ regul"ar.
\end{lemma}

\begin{satz}[Kettenregel] \label{SatzKettenregel}
Es seien $g: X \to Y$ im Punkt $x_0 \in X$ Fr\'echet-differenzierbar und $f: Y \to \overline{\R}$ eigentlich,
konvex und im Punkt $g(x_0)$ stetig.
Dann ist die Funktion $G(x) = f\big( g(x)\big)$ im Punkt $x_0$ regul"ar und es gilt
$$\partial G(x_0) = g'^*(x_0) \partial f\big(g(x_0)\big).$$
\end{satz}

%%%%%%%%%%%%%%%%%%%%%%%%%%%%%%%%%%%%%%%%%%%

\subsection{Das Subdifferential konkreter Funktionen}
\begin{beispiel} \label{SubdifferentialMaximum1}
{\rm Es sei $[t_0,t_1] \subset \R$ ein kompaktes Intervall.
Im Raum $C([t_0,t_1],\R)$ betrachten wir die Funktion
$$f\big(x(\cdot)\big) = \max_{t \in [t_0,t_1]} x(t).$$
Diese Funktion ist stetig, konvex und homogen.
Das Subdifferential $\partial f(0)$ besteht nach Definition aus denjenigen signierten regul"aren Borelschen Ma"sen $\mu$ auf $[t_0,t_1]$,
die der Bedingung
$$\max_{t \in [t_0,t_1]} x(t) \geq \int_{t_0}^{t_1} x(t) \, d\mu(t) \qquad \mbox{f"ur alle } x(\cdot) \in C([t_0,t_1],\R)$$
gen"ugen.
Hieraus folgt, dass das Ma"s $\mu$ nichtnegativ ist und $\|\mu\| = 1$ gilt.
Man erh"alt n"amlich f"ur $-x(t)$ nach obiger Ungleichung
$$\int_{t_0}^{t_1} x(t) \, d\mu(t) \geq - \max_{t \in [t_0,t_1]} \big( -x(t)\big) = \min_{t \in [t_0,t_1]} x(t).$$
Ist daher $x(t) \geq 0$ f"ur alle $t \in [t_0,t_1]$,
so ist auch $\displaystyle\int_{t_0}^{t_1} x(t) \, d\mu(t) \geq 0$ und $\mu$ ist positiv.
Weiterhin besitzen alle $\mu \in \partial f(0)$ die Norm $\|\mu\| = 1$,
denn nur dann gilt die Relation
$$\max_{t \in [t_0,t_1]} x(t) \geq \int_{t_0}^{t_1} x(t) \,d\mu(t) \geq \min_{t \in [t_0,t_1]} x(t)$$
f"ur konstante Funktionen $x(\cdot)$. \\
Offenbar ist auch umgekehrt, wenn $\mu \geq 0$ ist und Norm Eins hat, die Ungleichung
$$\max_{t \in [t_0,t_1]} x(t) \geq \int_{t_0}^{t_1} x(t) \, d\mu(t)$$
richtig.
Somit haben wir das Subdifferential der Funktion $f$ im Nullpunkt berechnet. \\[2mm]
F"ur $x(\cdot) \not= 0$ gilt f"ur die eigentliche homogene konvexe Funktion $f\big(x(\cdot)\big)$:
$$\partial f\big(x(\cdot)\big) = \big\{ x^* \in \partial f(0) \big| \langle x^*,x(\cdot) \rangle = f\big(x(\cdot)\big) \big\},$$
wobei $\partial f(0)$ die regul"aren Borelschen Ma"se mit $\|\mu\|=1$ enth"alt.
Wir zeigen, dass die Ma"se $\mu \in \partial f\big(x(\cdot)\big)$ auf den Mengen
$$T = \bigg\{ \tau \in [t_0,t_1] \,\bigg|\, x(\tau) = \max_{t \in [t_0,t_1]} x(t) \bigg\}$$
konzentriert sind:
Der K"urze halber bezeichne $M$ das Maximum von $x(t)$ auf $[t_0,t_1]$.
Die Menge $[t_0,t_1] \setminus T$ ist offen, da $T$ abgeschlossen ist.
Angenommen, es existiert eine me"sbare Menge $B \subseteq [t_0,t_1]$ mit $B \cap T = \emptyset$ und $\mu(B) > 0.$
Dann gilt
$$\int_B x(t)\, d\mu(t) < \int_B M\, d\mu(t) \Rightarrow \int_{[t_0,t_1] \setminus T} x(t)\, d\mu(t) < \int_{[t_0,t_1] \setminus T} M\, d\mu(t).$$
Damit folgt der Widerspruch
$$\max_{t\in [t_0,t_1]} x(t) = \int_{t_0}^{t_1} M\, d\mu(t) > \int_T M\, d\mu(t) + \int_{[t_0,t_1] \setminus T} x(t) \, d\mu(t)
   = \int_{t_0}^{t_1} x(t)\, d\mu(t).$$
Ist umgekehrt $\mu \geq 0$ und besitzt Norm Eins, dann geh"ort $\mu$ der Menge $\partial f(0)$ an.
Wenn $\mu$ zus"atzlich auf der Menge $T$ konzentriert ist, dann ist auch
$$\max_{t \in [t_0,t_1]} x(t) = \int_{t_0}^{t_1} x(t)\, d\mu(t)$$
erf"ullt.
Damit besteht das Subdifferential der Funktion $f$ in einem vom Nullpunkt verschiedenen Punkt $x(\cdot)$ aus den regul"aren
Borelschen Ma"sen auf $[t_0,t_1]$, deren Norm gleich Eins ist und die auf der Menge $T$ konzentriert sind. \hfill $\square$}
\end{beispiel}

\begin{beispiel} \label{SubdifferentialMaximum2}
{\rm Es sei $g(t,x)$ eine Funktion auf $[t_0,t_1]\times \R^n$,
die bez"uglich beider Ver"anderlicher stetig und f"ur jedes $t \in [t_0,t_1]$ nach $x$ stetig differenzierbar ist.
D.\,h., dass die Abbildung $g_x(t,x)$ in der Gesamtheit der Variablen stetig ist.
Unter diesen Voraussetzungen ist nach Beispiel \ref{DiffAbbildung} die Abbildung
$$\tilde{g}: C([t_0,t_1],\R^n) \to C([t_0,t_1],\R), \qquad \big[\tilde{g}\big(x(\cdot)\big)\big] (t) = g\big(t,x(t)\big),\quad t \in [t_0,t_1],$$
Fr\'echet-differenzierbar und es gilt
$$\big[\tilde{g}'\big(x(\cdot)\big) z(\cdot)\big] (t) = \big\langle g_x \big(t,x(t)\big), z(t) \big\rangle, \qquad t \in [t_0,t_1].$$
Weiterhin ist die Funktion $f\big(x(\cdot)\big)$ im Beispiel \ref{SubdifferentialMaximum1} auf $C([t_0,t_1],\R)$ stetig. \\[2mm]
Wir bestimmen nun das Subdifferential der Funktion
$$G\big(x(\cdot)\big) = f\Big(g\big(x(\cdot)\big)\Big) =\max_{t \in [t_0,t_1]} g\big(t,x(t)\big).$$
Dazu k"onnen wir die Kettenregel (Satz \ref{SatzKettenregel}) anwenden:
$$\partial G\big(x(\cdot)\big) = \tilde{g}'^*\big(x(\cdot)\big) \partial f \big(\tilde{g}\big(x(\cdot)\big)\big).$$
F"ur $x^* \in \tilde{g}'^*\big(x(\cdot)\big) \partial f \big(\tilde{g}\big(x(\cdot)\big)\big)$ und $z(\cdot) \in C([t_0,t_1],\R^n)$ gilt dann:
$$\big\langle x^*,z(\cdot) \big\rangle
   = \big\langle \tilde{g}'^*\big(x(\cdot)\big)\, \mu, z(\cdot) \big\rangle
   = \big\langle \mu , \tilde{g}'\big(x(\cdot)\big)z(\cdot) \big\rangle
   = \int_{t_0}^{t_1} \big\langle g_x \big(t,x(t)\big), z(t) \big\rangle \,d\mu(t).$$
Daher besteht das Subdifferential der Funktion $G$ im Punkt $x(\cdot)$ genau aus denjenigen stetigen linearen Funktionalen $x^*$,
die die Darstellung
$$\big\langle x^*,z(\cdot) \big\rangle = \int_{t_0}^{t_1} \big\langle g_x \big(t,x(t)\big), z(t) \big\rangle \,d\mu(t)$$
besitzen,
wobei das regul"are Borelsches Ma"s $\mu$ auf $T= \big\{ t \in [t_0,t_1] \,\big|\, g\big(t,x(t)\big) = G\big(x(\cdot)\big) \big\}$
konzentriert ist und $\|\mu\| =1$ gilt. \hfill $\square$}
\end{beispiel}

\begin{beispiel} \label{SubdifferentialMaximum3}
{\rm Wir betrachten im Raum $C_{\lim}(\R_+,\R)$ die Funktion
$$f\big(x(\cdot)\big) = \max_{t \in \overline{\R}_+} x(t).$$
Mit den gleichen Argumenten wie in Beispiel \ref{SubdifferentialMaximum1}
ist diese Funktion stetig, konvex und homogen.
Ferner besteht das Subdifferential $\partial f(0)$ aus denjenigen signierten regul"aren Borelschen Ma"sen $\mu$ auf $\overline{\R}_+$,
die nichtnegativ sind und die Totalvariation $\|\mu\| = 1$ besitzen. \\
Weiterhin besteht das Subdifferential der Funktion $f$ in einem vom Nullpunkt verschiedenen Punkt $x(\cdot)$ aus den regul"aren
Borelschen Ma"sen auf $\overline{\R}_+$, deren Norm gleich Eins ist und die auf der Menge
$T = \big\{ t \in \overline{\R}_+ \,\big|\, x(t) = f\big(x(\cdot)\big) \big\}$
konzentriert sind. \hfill $\square$}
\end{beispiel}

\begin{beispiel} \label{SubdifferentialMaximum4}
{\rm Es sei $g(t,x)$ eine Funktion auf $\R_+ \times \R^n$,
die bez"uglich beider Ver"anderlicher stetig und f"ur jedes $t \in \R_+$ nach $x$ stetig differenzierbar ist.
Unter diesen Voraussetzungen ist nach Beispiel \ref{DiffAbbildung} die Abbildung
$$\tilde{g}: C_{\lim}(\R_+,\R^n) \to C_{\lim}(\R_+,\R), \qquad \big[\tilde{g}\big(x(\cdot)\big)\big] (t) = g\big(t,x(t)\big),\quad t \in \R_+,$$
Fr\'echet-differenzierbar und es gilt
$$\big[\tilde{g}'\big(x(\cdot)\big) z(\cdot)\big] (t) = \big\langle g_x \big(t,x(t)\big), z(t) \big\rangle, \qquad t \in \R_+.$$
Weiterhin ist die Funktion $f$ im Beispiel \ref{SubdifferentialMaximum3} auf $C_{\lim}(\R_+,\R)$ stetig.
Auf die Funktion
$$G\big(x(\cdot)\big) = f\Big(g\big(x(\cdot)\big)\Big) =\max_{t \in \overline{\R}_+} g\big(t,x(t)\big).$$
wir die Kettenregel an:
Wie im Beispiel \ref{SubdifferentialMaximum2} besteht 
das Subdifferential der Funktion $G$ im Punkt $x(\cdot)$ genau aus denjenigen stetigen linearen Funktionalen $x^*$,
die die Darstellung
$$\big\langle x^*,z(\cdot) \big\rangle = \int_0^\infty \big\langle g_x \big(t,x(t)\big), z(t) \big\rangle \,d\mu(t) + 
                                         \lim_{t \to \infty} \big\langle g_x \big(t,x(t)\big), z(t) \big\rangle \, \mu(\{\infty\})$$
besitzen,
wobei das regul"are Borelsches Ma"s $\mu$ auf $T= \big\{ t \in \overline{\R}_+ \,\big|\, g\big(t,x(t)\big) = G\big(x(\cdot)\big) \big\}$
konzentriert ist und $\|\mu\| =1$ gilt. \hfill $\square$}
\end{beispiel}
\lhead[\thepage \hspace*{1mm} Mehrfache Nadelvariationen]{}
\rhead[]{Mehrfache Nadelvariationen \hspace*{1mm} \thepage}
\section{Mehrfache Nadelvariationen} \label{AnhangNV}
Wir folgen bei der Darstellung mehrfacher Nadelvariationen Ioffe \& Tichomirov \cite{Ioffe}.
Als erstes geben wir eine Konstruktion geeigneter Tr"agerfamilien an,
auf denen die Nadelvariationen durchgef"uhrt werden.
Danach definieren wir die mehrfachen Nadelvariationen und untersuchen abschlie"send die Eigenschaften,
die die mehrfache Nadelvariation dem Steuerungsproblem "ubergibt.

\subsection{Die Konstruktion nach Ioffe \& Tchomirov} \label{AbschnittMengenfamilie}
\begin{lemma} \label{LemmaMengenfamilie3}
Es seien $y_i(\cdot),\, y_i : [t_0,t_1] \to \R^{n_i},\, i=1,...,d,$ beschr"ankte messbare Vektorfunktionen.
Dann gibt es zu jedem $\delta > 0$ einparametrige Mengenfamilien $M_1(\alpha),...,M_d(\alpha)$
messbarer Teilmengen des Intervalls $[t_0,t_1]$,
wobei der Parameter $\alpha$ Werte zwischen $0$ und $1/d$ annimmt derart, dass f"ur alle
$i=1,...,d,$ $0 \leq \alpha' \leq \alpha \leq 1/d$ und $i \not= i'$ gilt:
\begin{eqnarray}
&& \label{Nadelvariation1} 
   \hspace*{-15mm} |M_i(\alpha)| = \alpha (t_1-t_0), \qquad M_i(\alpha') \subseteq M_i(\alpha),
   \qquad M_i(\alpha) \cap M_{i'}(\alpha') = \emptyset,\\
&& \label{Nadelvariation2}
   \hspace*{-15mm} \max_{t \in [t_0,t_1]}
   \bigg\| \int_{t_0}^t \big( \chi_{M_i(\alpha)}(\tau) - \chi_{M_i(\alpha')}(\tau) \big) y_i(\tau) \, d\tau
   - (\alpha-\alpha') \int_{t_0}^t y_i(\tau) \, d\tau \bigg\| \leq \delta |\alpha-\alpha'|.
\end{eqnarray}
\end{lemma}

Es seien $x_*(\cdot) \in C([t_0,t_1],\R^n)$, $u_*(\cdot),u_1(\cdot),...,u_d(\cdot) \in L_\infty\big([t_0,t_1],U\big)$
und es bezeichne $V_\gamma$ die Menge
$$V_\gamma= \{ (t,x) \in [t_0,t_1] \times \R^n\,|\, \|x-x_*(t)\|\leq \gamma\}.$$
Wir nehmen an, dass die Abbildungen $f(t,x,u)$, $\varphi(t,x,u)$ auf $V_\gamma \times \R^m$ stetig in der Gesamtheit der
Variablen und stetig differenzierbar bez"uglich $x$ sind.
D.\,h., dass die partiellen Ableitungen $f_x(t,x,u)$, $\varphi_x(t,x,u)$ auf $V_\gamma \times \R^m$ stetig in der 
Gesamtheit der Variablen sind.
Unter diesen Annahmen sind die Vektorfunktionen $y_i(\cdot)$,
\begin{eqnarray*}
y_i (t) &=& \Big( \varphi\big(t,x_*(t),u_i(t) \big) - \varphi\big(t,x_*(t),u_*(t) \big), \\
             & & \hspace*{20mm} f\big(t,x_*(t),u_i(t) \big) - f\big(t,x_*(t),u_*(t) \big) \Big), \quad i = 1,...,d,
\end{eqnarray*}
auf $[t_0,t_1]$ messbar und beschr"ankt, weil die Steuerungen $u_*(\cdot),u_1(\cdot),...,u_{d}(\cdot)$ messbar, beschr"ankt und $f, \varphi$ stetig sind.
Daher existieren nach Lemma \ref{LemmaMengenfamilie3} auf dem Intervall $[t_0,t_1]$
einparametrige Mengenfamilien $\big\{ M_i(\alpha) \big\}$ mit
\begin{eqnarray*}
&& |M_i(\alpha)| = \alpha (t_1-t_0), \qquad M_i(\alpha') \subseteq M_i(\alpha),
   \qquad M_i(\alpha) \cap M_{i'}(\alpha') = \emptyset,\\
&& \max_{t \in [t_0,t_1]}
   \bigg\| \int_{t_0}^t \big( \chi_{M_i(\alpha)}(\tau) - \chi_{M_i(\alpha')}(\tau) \big) y_i(\tau) \, d\tau
   - (\alpha-\alpha') \int_{t_0}^t y_i(\tau) \, d\tau \bigg\| \leq \delta |\alpha-\alpha'|
\end{eqnarray*}
f"ur alle $i=1,...,d,$ $0 \leq \alpha' \leq \alpha \leq 1/d$ und $i \not= i'$.
Auf dem Quader
$$Q^{d} = \Big\{ \alpha = (\alpha_1,...,\alpha_{d}) \in \R^{d} \,\Big|\,
                                0 \leq \alpha_i  \leq 1/d, \, i =1,...,d \Big\}$$
ist die Abbildung $\alpha \to u_\alpha(\cdot) \in L_\infty\big([t_0,t_1],U\big),$
$$u_\alpha(t) = u_*(t) + \sum_{i=1}^{d} \chi_{M_i(\alpha_i)}(t) \cdot \big( u_i(t)-u_*(t) \big),$$
wohldefiniert und $u_\alpha(\cdot)$ nennen wir mehrfache Nadelvariation von $u_*(\cdot)$. \\[2mm]
Wir definieren die Mengen
\begin{eqnarray*}
\Sigma(\Delta) &=& \left\{ \alpha = (\alpha_1,...,\alpha_{d}) \in \R^{d} \,\bigg|\,
                                \alpha_1,...,\alpha_{d} \geq 0,\, \sum_{i=1}^{d} \alpha_i \leq \Delta \right\}, \\
V &=& \{ x(\cdot) \in C([t_0,t_1],\R^n) \,\big|\, \|x(\cdot) - x_*(\cdot)\|_{\infty} \leq \sigma \}.
\end{eqnarray*}
Dann betrachten wir f"ur $t \in [t_0,t_1]$ die Abbildungen
\begin{eqnarray*}
\Phi_1\big(x(\cdot),\alpha\big)(t) &=& \int_{t_0}^t \big[\varphi\big(\tau,x(\tau),u_\alpha(\tau)\big)
                                                 - \varphi\big(\tau,x_*(\tau),u_*(\tau)\big)\big] \, d\tau, \\
\Phi_2\big(x(\cdot),\alpha\big) &=&  \int_{t_0}^{t_1} \big[f\big(t,x(t),u_\alpha(t)\big)-f\big(t,x(t),u_*(t)\big)\big] \, dt
\end{eqnarray*}
und die linearen Operatoren
\begin{eqnarray*}
\Lambda_1\big(x(\cdot),\alpha\big)(t) &=& \int_{t_0}^t \Big[\varphi_x\big(\tau,x_*(\tau),u_*(\tau)\big) \big(x(\tau)-x_*(\tau)\big) \\
&&  \hspace*{5mm}+ \sum_{i=1}^{d} \alpha_i \cdot 
   \Big( \varphi\big(\tau,x_*(\tau),u_i(\tau)\big) - \varphi\big(\tau,x_*(\tau),u_*(\tau)\big)\Big)\Big] \, d\tau, \\
\Lambda_2\big(x(\cdot),\alpha\big) &=& \sum_{i=1}^{d} \int_{t_0}^{t_1} \alpha_i \cdot \Big( f\big(t,x(t),u_i(t)\big)
                                                    - f\big(t,x(t),u_*(t)\big)\Big) \, dt.
\end{eqnarray*}
Dann ergeben sich folgende Resultate (vgl. \cite{Ioffe}):
\begin{lemma} \label{LemmaEigenschaftNadelvariation1}
Es sei $\Delta \in (0,1/d]$ hinreichend klein. 
\begin{enumerate}
\item[(a)] F"ur alle $\alpha, \alpha' \in \Sigma(\Delta)$ und alle $x(\cdot), x'(\cdot) \in V$ gilt:
           \begin{eqnarray}
           && \Big\| \big[\Phi_1\big(x(\cdot),\alpha\big)-\Phi_1\big(x'(\cdot),\alpha'\big)
                    -\Lambda_1\big(x(\cdot),\alpha\big)+\Lambda_1\big(x'(\cdot),\alpha'\big)\big] (\cdot) \Big\|_\infty \nonumber \\
           && \label{EigenschaftNadelvariation7}
              \hspace*{20mm} \leq\delta \bigg( \|x(\cdot)-x'(\cdot)\|_\infty + \sum_{i=1}^d |\alpha_i - \alpha_i'| \bigg).
           \end{eqnarray}
\item[(b)] F"ur alle $\alpha \in \Sigma(\Delta)$ und alle $x(\cdot) \in V$ gilt
           \begin{equation} \label{EigenschaftNadelvariation8}
           \Phi_2\big(x(\cdot),\alpha\big) - \Lambda_2\big(x(\cdot),\alpha\big) \leq \delta \sum_{i=1}^d \alpha_i.
           \end{equation}
\end{enumerate}
\end{lemma}

%%%%%%%%%%%%%%%%%%%%%%%%%%%%%%%%%%%%%%%%%%%

\subsection{Mehrfache Nadelvariationen "uber dem unendlichen Zeithorizont} \label{AnhangNVUH}
Wir werden nun zeigen,
wie wir f"ur das unbeschr"ankte Zeitintervall und f"ur einen lokal unbeschr"ankten Integranden
die Konstruktion mehrfacher Nadelvariationen auf eine kompakte Menge $K$ zur"uckf"uhren k"onnen. \\[2mm]
Die Konstruktion der Tr"agermengen l"asst sich unmittelbar auf den Fall anwenden,
dass statt dem Intervall $[t_0,t_1]$ eine kompakte Menge $K \subset \R_+$ betrachtet wird.
Dann ergibt sich:

\begin{lemma} \label{LemmaMengenfamilie4}
Es sei $K \subset \R_+$ kompakt und es seien $y_i(\cdot),\, y_i : K \to \R^{n_i},\, i=1,...,d,$ beschr"ankte messbare Vektorfunktionen.
Dann existieren zu jedem $\delta > 0$ einparametrige Mengenfamilien $M_1(\alpha),...,M_d(\alpha)$, $0\leq \alpha \leq 1/d$,
messbarer Teilmengen der Menge $K$ derart, dass f"ur alle
$i=1,...,d,$ $0 \leq \alpha' \leq \alpha \leq 1/d$ und $i \not= i'$ gilt:
\begin{eqnarray*}
&& \hspace*{-5mm} |M_i(\alpha)| = \alpha |K|, \qquad M_i(\alpha') \subseteq M_i(\alpha),
   \qquad M_i(\alpha) \cap M_{i'}(\alpha') = \emptyset,\\
&& \hspace*{-5mm} \max_{t \in K}
   \bigg\| \int_{[0,t] \cap K} \big( \chi_{M_i(\alpha)}(\tau) - \chi_{M_i(\alpha')}(\tau) \big) y_i(\tau) \, d\tau 
            - (\alpha-\alpha') \int_{[0,t] \cap K} y_i(\tau) \, d\tau \bigg\| \leq \delta |\alpha-\alpha'|.
\end{eqnarray*}
\end{lemma}

F"ur die weiteren Betrachtungen geben wir einen "Uberblick "uber die ben"otigten Eigenschaften an die Aufgabe (\ref{PAUH1})--(\ref{PAUH5}),
die wir im Kapitel \ref{KapitelStrong} getroffen haben: \\
Es seien $x_*(\cdot) \in C_{\lim}(\R_+,\R^n)$ und $u_*(\cdot) \in L_\infty(\R_+,U)$.
Zu $x_*(\cdot)$ bezeichnet $V_\gamma$ die Menge
$$V_\gamma= \{ (t,x) \in \overline{\R}_+ \times \R^n \,|\, \|x-x*(t)\| \leq \gamma\}.$$
Wir nehmen an,
dass es zu jeder kompakten Menge $U_1 \subseteq \R^m$ eine Zahl $\gamma>0$ derart gibt,
dass auf der Menge $V_\gamma \times U_1$ die Abbildungen
$f(t,x,u)$ und $\varphi(t,x,u)$ gleichm"a"sig stetig und gleichm"a"sig stetig differenzierbar bez"uglich $x$ sind.
Ferner seien f"ur das Paar $\big(x_*(\cdot),u_*(\cdot)\big)$ folgende Eigenschaften erf"ullt:
Es sind die Lebesgue-Integrale
$$\int_0^\infty \big\|\varphi\big(t,x_*(t),u_*(t)\big)\big\| \, dt, \qquad \int_0^\infty \big\|\varphi_x\big(t,x_*(t),u_*(t)\big)\big\| \, dt$$
endlich.
Au"serdem nehmen wir an, es existiert zu jedem $\delta>0$ ein $T>0$ mit
\begin{eqnarray}
&& \int_T^\infty \big\| \varphi\big(t,x(t),u_*(t)\big)-\varphi\big(t,x'(t),u_*(t)\big) - \varphi_x\big(t,x_*(t),u_*(t)\big)\big(x(t)-x'(t)\big) \big\| \, dt
   \nonumber \\
&& \label{AnhangNVUH1} \hspace*{20mm} \leq \delta \|x(\cdot)-x'(\cdot)\|_\infty
\end{eqnarray}
f"ur alle $x(\cdot), x'(\cdot) \in C_{\lim}(\R_+,\R^n)$ mit $\|x(\cdot)-x_*(\cdot)\|_\infty < \gamma$, $\|x'(\cdot)-x_*(\cdot)\|_\infty < \gamma$. \\
Abschlie"send bezeichnet $\mathscr{U}$ die Menge aller $u(\cdot) \in L_\infty(\R_+,U)$,
die die Darstellung
$$u(t)=u_*(t) + \chi_M(t)\big(w(t)-u_*(t)\big)$$
mit $w(\cdot) \in L_\infty(\R_+,U)$ und einer me"sbar und beschr"ankten Menge $M \subset \R_+$ besitzen. \\[2mm]
Es seien $u_1(\cdot),...,u_d(\cdot) \in \mathscr{U}$ und $\delta >0$ gegeben.
Dann l"asst sich eine Zahl $T>0$ derart w"ahlen, dass
die Mengen $M_i$, die in den Darstellungen der Steuerungen $u_i(\cdot) \in \mathscr{U}$ auftreten, im Intervall $[0,T]$ enthalten sind und
die Relation (\ref{AnhangNVUH1}) mit $\delta/3$ erf"ullt ist. \\
Ferner existiert nach dem Satz von Lusin eine kompakte Menge $K \subseteq [0,T]$,
auf der die Funktion $\omega(\cdot)$ stetig ist und zudem f"ur $i=1,...,d$ folgende Relationen gelten:
\begin{eqnarray*}
&& \int_{[0,T] \setminus K} \big\|\varphi_x\big(t,x_*(t),u_*(t)\big)\big\| \, dt \leq \frac{\delta}{3}, \\
&& \int_{[0,T] \setminus K} \big\|\varphi\big(t,x_*(t),u_i(t)\big)-\varphi\big(t,x_*(t),u_*(t)\big)\big\| \, dt \leq \frac{\delta}{3}, \\
&& \int_{[0,T] \setminus K} \omega(t)\big|f\big(t,x_*(t),u_i(t)\big)-f\big(t,x_*(t),u_*(t)\big)\big| \, dt \leq \frac{\delta}{2}.
\end{eqnarray*}
Wir betrachten im Folgenden die Funktionen $y_i(\cdot)$,
\begin{eqnarray*}
y_i (t) &=& \Big( \varphi\big(t,x_*(t),u_i(t) \big) - \varphi\big(t,x_*(t),u_*(t) \big), \\
             & & \hspace*{20mm} \omega(t)\big[ f\big(t,x_*(t),u_i(t) \big) - f\big(t,x_*(t),u_*(t) \big)\big] \Big), \quad i = 1,...,d,
\end{eqnarray*}
die nach Wahl von $K$ "uber dieser Menge messbar und beschr"ankt sind.
Mit den zugeh"origen Tr"agermengen $M_1(\alpha),...,M_d(\alpha) \subseteq K$ aus Lemma \ref{LemmaMengenfamilie4} definieren wir "uber $\R_+$ 
die Abbildung $\alpha \to u_\alpha(\cdot) \in \mathscr{U}$:
$$u_\alpha(t) = u_*(t) + \sum_{i=1}^{d} \chi_{M_i(\alpha_i)}(t) \cdot \big( u_i(t)-u_*(t) \big).$$
Diese ist auf dem Quader $Q^{d}$,
$$Q^{d} = \Big\{ \alpha = (\alpha_1,...,\alpha_{d}) \in \R^{d} \,\Big|\,
                                0 \leq \alpha_i  \leq 1/d, \, i =1,...,d \Big\},$$
wohldefiniert.
Im Weiteren bezeichnen $\Sigma(\Delta)$ und $V$ die Mengen
\begin{eqnarray*}
\Sigma(\Delta) &=& \left\{ \alpha = (\alpha_1,...,\alpha_{d}) \in \R^{d} \,\bigg|\,
                                \alpha_1,...,\alpha_{d} \geq 0,\, \sum_{i=1}^{d} \alpha_i \leq \Delta \right\}, \\
V &=& \{ x(\cdot) \in C_{\lim}(\R_+,\R^n) \,\big|\, \|x(\cdot) - x_*(\cdot)\|_{\infty} \leq \sigma \}.
\end{eqnarray*}
Auf diesen Mengen betrachten wir die Abbildungen
\begin{eqnarray*}
    \Phi_1\big(x(\cdot),\alpha\big)(t)
&=& \int_0^t \big[\varphi\big(\tau,x(\tau),u_\alpha(\tau)\big) - \varphi\big(\tau,x_*(\tau),u_*(\tau)\big)\big] \, d\tau, \quad t \in \R_+, \\
    \Lambda_1\big(x(\cdot),\alpha\big)(t)
&=& \int_0^t \bigg[\varphi_x\big(\tau,x_*(\tau),u_*(\tau)\big) \big(x(\tau)-x_*(\tau)\big) \\
& & \hspace*{7mm} + \sum_{i=1}^{d} \alpha_i \cdot 
    \Big( \varphi\big(\tau,x_*(\tau),u_i(\tau)\big) - \varphi\big(\tau,x_*(\tau),u_*(\tau)\big)\Big)\bigg] \, d\tau, \quad t \in \R_+
\end{eqnarray*}
und die Funktionale
\begin{eqnarray*}
    \Phi_2\big(x(\cdot),\alpha\big)
&=& \int_0^\infty \omega(t)\big[f\big(t,x(t),u_\alpha(t)\big)-f\big(t,x(t),u_*(t)\big)\big] \, dt, \\
    \Lambda_2\big(x(\cdot),\alpha\big)
&=& \sum_{i=1}^{d} \int_0^\infty \alpha_i \cdot \omega(t) \big[ f\big(t,x(t),u_i(t)\big)
                                                    - f\big(t,x(t),u_*(t)\big)\big] \, dt.
\end{eqnarray*}

\begin{lemma} \label{LemmaNVUH1}
Es existieren $\Delta \in (0,1/d]$ und $\sigma>0$ derart,
dass f"ur alle $\alpha, \alpha' \in \Sigma(\Delta)$ und alle $x(\cdot), x'(\cdot) \in V$ gelten:
\begin{eqnarray}
&& \Big\| \big[\Phi_1\big(x(\cdot),\alpha\big)-\Phi_1\big(x'(\cdot),\alpha'\big)
          -\Lambda_1\big(x(\cdot),\alpha\big)-\Lambda_1\big(x'(\cdot),\alpha'\big)\big] (\cdot) \Big\|_\infty \nonumber \\
&& \label{NVUH1}
   \hspace*{20mm} \leq\delta \bigg( \|x(\cdot)-x'(\cdot)\|_\infty + \sum_{i=1}^d |\alpha_i - \alpha_i'| \bigg).
\end{eqnarray}  
\end{lemma}

{\bf Beweis} Beachten wir, dass die Abbildung $\alpha \to u_\alpha(\cdot)$ nur f"ur $t \in K$ zum Tragen kommt,
so k"onnen wir die linke Seite in (\ref{NVUH1}) gegen folgenden Ausdruck nach oben absch"atzen:

\begin{eqnarray*}
& & \int_{[0,T] \setminus K}  \big\|\varphi_x\big(t,x_*(t),u_*(t)\big)\big(x(t)-x'(t)\big)\big\| \, dt  \\
& & \hspace*{15mm} + \sum_{i=1}^d (\alpha_i- \alpha'_i) \cdot
    \int_{[0,T] \setminus K} \big\|\varphi\big(t,x_*(t),u_i(t)\big)-\varphi\big(t,x_*(t),u_*(t)\big)\big\| \, dt \\
&+& \max_{t \in K} \bigg\| \int_{[0,t] \cap K} \bigg[\varphi \big( \tau, x(\tau), u_\alpha(\tau) \big)
                            - \varphi \big( \tau, x'(\tau), u_{\alpha'}(\tau) \big) \\
& & \hspace*{30mm} - \varphi_x \big( \tau, x_*(\tau), u_*(\tau) \big) \big( x(\tau) - x'(\tau) \big) \\
& & \hspace*{30mm} - \sum_{i=1}^d (\alpha_i- \alpha'_i) \cdot 
    \Big( \varphi \big( \tau, x_*(\tau), u_i(\tau) \big) - \varphi \big( \tau, x_*(\tau), u_*(\tau) \big) \Big) \bigg] d\tau \bigg\| \\
&+& \int_T^\infty \big\| \varphi\big(t,x(t),u_*(t)\big)-\varphi\big(t,x'(t),u_*(t)\big)
                         - \varphi_x\big(t,x_*(t),u_*(t)\big)\big(x(t)-x'(t)\big) \big\| \, dt.
\end{eqnarray*}
Nach obiger Wahl der Menge $K$ f"allt der erste Summand kleiner oder gleich
$$\frac{\delta}{3} \bigg( \|x(\cdot)-x'(\cdot)\|_\infty + \sum_{i=1}^d |\alpha_i - \alpha_i'| \bigg)$$
aus.
Nach Wahl der Zahl $T$ ist der letzte Summand kleiner gleich
$$\frac{\delta}{3} \|x(\cdot)-x'(\cdot)\|_\infty.$$
Bez"uglich der kompakten Menge $K$ folgt diese Absch"atzung direkt aus Lemma \ref{LemmaEigenschaftNadelvariation1},
da im Beweis das Intervall $[t_0,t_1]$ durch die kompakte Menge $K$ ersetzt werden darf,
auf der die Dichtefunktion $\omega(\cdot)$ stetig ist. \hfill $\blacksquare$

\begin{lemma} \label{LemmaNVUH2}
Es existieren $\Delta \in (0,1/d]$ und $\sigma>0$ derart,
dass die Ungleichung
\begin{equation} \label{NVUH2}
\Phi_2\big(x(\cdot),\alpha\big) - \Lambda_2\big(x(\cdot),\alpha\big) \leq \delta \sum_{i=1}^d \alpha_i
\end{equation}
f"ur alle $\alpha \in \Sigma(\Delta)$ und alle $x(\cdot) \in V$ gilt.
\end{lemma}

{\bf Beweis} Die linke Seite in (\ref{NVUH2}) ist gleich
$$\sum_{i=1}^{d} \int_0^\infty \Big(\chi_{M_i(\alpha_i)}(t) -\alpha_i \Big) \cdot
                \Big( \omega(t) \big[f \big(t,x(t),u_i(t)\big) - f\big(t,x(t),u_*(t)\big)\big]\Big) \, dt.$$
Beachten wir $M_i(\alpha_i) \subseteq K$, so folgt nach Wahl der Menge $K$:
$$\sum_{i=1}^{d} \int_{[0,T] \setminus K} \alpha_i \cdot
                \Big( \omega(t) \big[f \big(t,x(t),u_i(t)\big) - f\big(t,x(t),u_*(t)\big)\big]\Big) \, dt
    \leq \frac{\delta}{2} \sum_{i=1}^d \alpha_i.$$
Au"serdem gilt $u_i(t)=u_*(t)$ f"ur alle $t >T$ und $i=1,...,d$.
Daher ist
$$\sum_{i=1}^{d} \int_T^\infty \alpha_i \cdot \Big( \omega(t) \big[f \big(t,x(t),u_i(t)\big) - f\big(t,x(t),u_*(t)\big)\big]\Big) \, dt = 0.$$
Ferner erhalten wir "uber der kompakten Menge $K$
$$\sum_{i=1}^{d} \int_K \Big(\chi_{M_i(\alpha_i)}(t) -\alpha_i \Big) \cdot
                \Big( \omega(t) \big[f \big(t,x(t),u_i(t)\big) - f\big(t,x(t),u_*(t)\big)\big]\Big) \, dt$$
unmittelbar aus Lemma \ref{LemmaEigenschaftNadelvariation1},
da wieder das Intervall $[t_0,t_1]$ durch die kompakte Menge $K$ ersetzt werden darf,
auf der die Dichtefunktion $\omega(\cdot)$ stetig ist. \hfill $\blacksquare$
\lhead[\thepage \hspace*{1mm} Theorie der Extremalaufgaben]{}
\rhead[]{Theorie der Extremalaufgaben \hspace*{1mm} \thepage}
\section{Ein Extremalprinzip f"ur ein starkes lokales Minimum}
Es seien $X$, $Y$ Banachr"aume, $\mathscr{U}$ ein normierter Raum und $U \subseteq \mathscr{U}$.
Ferner seien $J$ ein Funktional auf $X \times U$,
$\mathscr{F}$ eine Abbildung des Produktes $X \times U$ in den Raum $Y$ und $G_j:X \to \R$ f"ur $j=1,...,l$. \\
Unter diesen Angaben betrachten wir in diesem Abschnitt die Extremalaufgabe
\begin{equation}\label{EAAnhang3}
J(x,u) \to \inf; \quad \mathscr{F}(x,u)=0, \quad G_j(x) \leq 0, \; j=1,...,l, \quad x \in X,\; u \in U.
\end{equation}
Der Punkt $(x,u)$ ist ein zul"assiges Element der Aufgabe (\ref{EAAnhang3}),
falls s"amtliche Nebenbedingungen erf"ullt sind.
Ein Punkt $(x_*,u_*)$ hei"st ein starkes\index{Minimum, schwaches lokales!--, starkes lokales} lokales Minimum der Extremalaufgabe (\ref{EAAnhang3}),
wenn ein $\varepsilon >0$ derart existiert,
dass f"ur alle zul"assigen Paare $(x,u)$ mit $\|x-x_*\|_X \leq \varepsilon$ die Ungleichung $J(x_*,u_*) \leq J(x,u)$ gilt.
Auf $X \times U \times \R \times Y^* \times \R^l$ definieren wir die Lagrange-Funktion
$$\mathscr{L}(x,u,\lambda_0,y^*,\lambda_1,...,\lambda_l)
  = \lambda_0 J(x,u) + \langle y^*, \mathscr{F}(x,u) \rangle + \sum_{j=1}^l  \lambda_j G_j(x).$$
Au"serdem bezeichnet in diesem Abschnitt $\Delta\Sigma^m$ das folgende $d$-dimensionale Simplex:
$$\Sigma(\Delta) = \bigg\{ \alpha = (\alpha_1,...,\alpha_d) \in \R^d \,\bigg|\, \alpha_1,...,\alpha_d \geq 0,\, \sum_{i=1}^{d} \alpha_i \leq \Delta \bigg\}.$$

\begin{theorem}[Extremalprinzip] \label{SatzExtremalprinzipStark}
Sei $(x_*,u_*)$ ein zul"assiges Element der Aufgabe (\ref{EAAnhang3}).
\begin{enumerate}
\item[(A)] Wir nehmen an, dass der Punkt $x_*$ eine Umgebung $V$ besitzt mit:
           \begin{enumerate}
           \item[(A$_1$)] F"ur jedes $u \in U$ ist $x \to J(x,u)$ im Punkt $x_*$ Fr\'echet-differenzierbar;
           \item[(A$_2$)] F"ur jedes $u \in U$ ist $x \to \mathscr{F}(x,u)$ im Punkt $x_*$ G\^ateaux-differenzierbar;
           \item[(A$_3$)] Die Funktionen $G_j(x)$ sind im Punkt $x_*$ lokalkonvex und bez"uglich jeder Richtung gleichm"a"sig differenzierbar.
           \end{enumerate}
\item[(B)] Weiterhin setzen wir voraus, dass der Operator $\mathscr{F}_x(x_*,u_*)$ eine endliche Kodimension besitzt.
\item[(C)] Au"serdem nehmen wir an,
           dass zu jedem endlichen System von Punkten $u_1,...,u_d$ aus $U$ und zu jedem $\delta>0$ eine Zahl $\Delta >0$ und eine Abbildung
           $u:\Sigma(\Delta) \to U$ derart existieren, dass $u(0)=u_*$ gilt und f"ur alle $x,x'$ aus $V'$ und alle $\alpha,\alpha'$
           aus $\Sigma(\Delta)$ folgende Ungleichungen erf"ullt sind:
           \begin{eqnarray*}
           && \hspace*{-10mm} \bigg\| \mathscr{F}\big(x,u(\alpha)\big) - \mathscr{F}\big(x',u(\alpha')\big)-\mathscr{F}_x(x_*,u_*)(x-x') \\
           && -\sum_{k=1}^d (\alpha_k - \alpha_k') \big( \mathscr{F}(x_*,u_k)-\mathscr{F}(x_*,u_*)\big) \bigg\| 
              \leq \delta \bigg( \|x-x'\| + \sum_{k=1}^d |\alpha_k - \alpha_k'|\bigg), \\
           && \hspace*{-10mm} J\big(x,u(\alpha)\big) - J(x,u_*) - \sum_{k=1}^d \alpha_k \big( J(x,u_k)-J(x,u_*)\big)
              \leq \delta \bigg( \|x-x_*\| + \sum_{k=1}^d \alpha_k \bigg).
           \end{eqnarray*}
\end{enumerate}
Ist dann $(x_*,u_*)$ starke lokale Minimalstelle der Aufgabe (\ref{EAAnhang3}),
so existieren nicht gleichzeitig verschwindende Lagrangesche Multiplikatoren $\lambda_0 \geq 0,...,\lambda_l \geq 0$ und $y^* \in Y^*$ derart,
dass folgende Bedingungen gelten:
\begin{enumerate}
\item[(a)] Die Lagrange-Funktion besitzt bez"uglich $x$ in $x_*$ einen station"aren Punkt, d.\,h.
           \begin{equation}\label{SatzEPstark1}
           0 \in \partial_x \mathscr{L}(x_*,u_*,\lambda_0,y^*,\lambda_1,...,\lambda_l);
           \end{equation}   
\item[(b)] Die Lagrange-Funktion erf"ullt bez"uglich $u$ in $u_*$ die Minimumbedingung
           \begin{equation}\label{SatzEPstark2}
           \mathscr{L}(x_*,u_*,\lambda_0,y^*,\lambda_1,...,\lambda_l) = \min_{u \in U} \mathscr{L}(x_*,u,\lambda_0,y^*,\lambda_1,...,\lambda_l);
           \end{equation}
\item[(c)] Die komplement"aren Schlupfbedingungen gelten, d.\,h.
           \begin{equation}\label{SatzEPstark3}
           0 = \lambda_j G_j(x_*) \qquad\mbox{f"ur } j=1,...,l.
           \end{equation}
\end{enumerate}
\end{theorem}

Das Extremalprinzip \ref{SatzExtremalprinzipStark} ist eine Vereinfachung des Extremalprinzips f"ur regul"ar lokalkonvexe Aufgaben in
Ioffe \& Tichomirov \cite{Ioffe}.
In \cite{Ioffe} h"angt die Steuerungsvariable $u$ zus"atzlich zu $\alpha$ von der Zustandsvariable $x$ ab.   
\end{appendix}       
   
\newpage
\addcontentsline{toc}{section}{Literatur}
\lhead[\thepage \, Literatur]{Optimale Steuerung mit unendlichem Zeithorizont}
\rhead[Optimale Steuerung mit unendlichem Zeithorizont]{Literatur \thepage}


\begin{thebibliography}{...}             
\bibitem{Arrow} Arrow,\,K.J., Kurz,\,M.: Public Investment, the Rate of Return, and Optimal Fiscal Policies.
                The Johns Hopkins Press, Baltimore (1970).
\bibitem{AseKry} Aseev,\,S.M., Kryazhimskii,\,A.V.: The Pontryagin Maximum Principle and Optimal Economic Growth Problems.
                 Proc. Steklov Inst. Math., 257, 1--255 (2007).
\bibitem{AseVel} Aseev,\,S.M., Veliov,\,V.M.: Maximum Principle for Infinite-Horizon Optimal Control Problems with Dominating Discount.
                 Dynamics of Continuous, Discrete and Impulsive Systems, Series B, Vol. 19 1-2 (2012).
\bibitem{AseVel2} Aseev,\,S.M., Veliov,\,V.M.: Needle variations in infinite-horizon optimal control.
                  Variational and Optimal Control Problems on Unbounded Domains,
                  eds. G.\,Wolansky, A.J.\,Zaslavski, Amer. Math. Soc. Contemporary Mathematics, 619 (2014).
\bibitem{AseVel3} Aseev,\,S.M., Veliov,\,V.M.: Maximum principle for infinite-horizon optimal control problems under weak regularity
                  assumptions. Trudy Inst. Mat. i Mekh. UrO RAN, 20, no.3, 41--57 (2014).
\bibitem{Brodskii} Brodskii,\,Yu.I.: Necessary Conditions for a Weak Extremum in Optimal Control Problems on an Infinite Time Interval.
                   Mat. Sb. (N.S.), Vol. 105(147), Number 3, 371--388 (1978).
\bibitem{Clarke} Clarke,\,F.: Optimization and Nonsmooth Analysis. John Wiley \& Sons, New York (1983).
\bibitem{Gale} Gale,\,D.: On Optimal Development in a Multi-Sector Economy. RES 34, 1--18 (1967).
\bibitem{Halkin} Halkin,\,H.: Necessary conditions for optimal control problems with infinite horizons. Econometrica 42, 267--272 (1979).
\bibitem{Ioffe} Ioffe,\,A.D., Tichomirov,\,V.M.: Theorie der Extremalaufgaben. VEB Deutscher Verlag der Wissenschaften Berlin, (1979).
\bibitem{Ljusternik} Ljusternik,\,L.A.: On Conditional Extremums of Functionals. Mat. Sb. 41, 390--401 (1934).
\bibitem{Pickenhain} Pickenhain,\,S.: Hilbert Space Treatment of Optimal Control Problems with Infinite Horizon.
                     Bock,\,H.G., Phu,\,H.X., Rannacher,\,R., Schloeder,\,J.P.:
                     Modeling, Simulation and Optimization of Complex Processes -- HPSC 2012, Springer, 169--182 (2014).
\bibitem{Ramsey} Ramsey\,F.P.: A Mathematical Theory of Saving. Econ. J. 38, 543--559 (1928).
\bibitem{Rockafellar} Rockafellar,\,R.T.: Convex Analysis. Princeton University Press, Princeton, New Jersey (1970).
\bibitem{Seierstad} Seierstad,\,A., Syds\ae ter,\,K.: Optimal Control Theory with Economic Applications.
                    North-Holland Amsterdam-New York-Oxford-Tokyo, (1987).
\bibitem{TauchnitzWMPIHOC} Tauchnitz,\,N.:: Weak Local Extrema in Infinite Horizon Optimal Control Problems.
                           Revised Preprint https://arxiv.org/abs/1510.04544 (2018).
\bibitem{TauchnitzOC} Tauchnitz,\,N.: Necessary Conditions in Optimal Control Theory.
                      Revised Preprint https://arxiv.org/abs/1610.02829 (2018).
\bibitem{VonWeiz} Von Weizsaecker,\,C.C.: Existence of optimal programs of accumulation for an infinite time horizon.
                  Review of Economic Studies 32, 85--92 (1965).
\end{thebibliography}
\end{document}